\documentclass{article}
\usepackage{amsmath}
\usepackage{mathrsfs}
\usepackage{txfonts}
\usepackage{bm}
\usepackage{latexsym,amsfonts,amssymb}
\begin{document}
\setcounter{page}{1}
\newtheorem{thm}{Theorem}[section]
\newtheorem{fthm}[thm]{Fundamental Theorem}
\newtheorem{dfn}[thm]{Definition}
\newtheorem{rem}[thm]{Remark}
\newtheorem{lem}[thm]{Lemma}
\newtheorem{cor}[thm]{Corollary}
\newtheorem{exa}[thm]{Example}
\newtheorem{pro}[thm]{Proposition}
\newtheorem{prob}[thm]{Problem}

\newtheorem{con}[thm]{Conjecture}
\renewcommand{\thefootnote}{\fnsymbol{footnote}}
\newcommand{\qed}{\hfill\Box\medskip}
\newcommand{\proof}{{\it Proof.\quad}}
\newtheorem{ob}[thm]{Observation}
\newcommand{\rmnum}[1]{\romannumeral #1}
\renewcommand{\abovewithdelims}[2]{%
\genfrac{[}{]}{0pt}{}{#1}{#2}}

\renewcommand{\thefootnote}{\fnsymbol{footnote}}

\begin{center}{\LARGE\bf Graph homomorphisms  on rectangular matrices  over
division rings II}\footnote{Project 11371072 supported by National Natural Science Foundation of China}
\vspace{0.5 true cm}
\end{center}
\begin{center}
Li-Ping Huang\footnote{E-mail: \ lipingmath@163.com (L.P. Huang);
zhaokangmath@126.com (K. Zhao)}, \ Kang Zhao
\\
{\small  School of Math. and Statis., {Changsha University of Science and Technology, Changsha, 410004, China}}
\end{center}

\begin{abstract}
Let ${\mathbb{D}}^{m\times n}$  be the set of $m\times n$  matrices  over a division ring $\mathbb{D}$. Two matrices $A,B\in {\mathbb{D}}^{m\times n}$ are adjacent if
${\rm rank}(A-B)=1$. By the adjacency, ${\mathbb{D}}^{m\times n}$ is a connected graph.
Suppose  $\mathbb{D}, \mathbb{D}'$ are division rings and $m,n,m',n'\geq2$ are integers. We determine  additive graph homomorphisms
from ${\mathbb{D}}^{m\times n}$ to ${\mathbb{D}'}^{m'\times n'}$. When $|\mathbb{D}|\geq 4$, we characterize the graph homomorphism
$\varphi:  {\mathbb{D}}^{n\times n}\rightarrow  {\mathbb{D}'}^{m'\times n'}$ if $\varphi(0)=0$ and there exists $A_0\in {\mathbb{D}}^{n\times n}$ such that ${\rm rank}(\varphi(A_0))=n$.
We also discuss  properties and ranges on degenerate  graph homomorphisms.
If $f:{\mathbb{D}}^{m\times n}\rightarrow  {\mathbb{D}'}^{m'\times n'}$ (where ${\rm min}\{m,n\}=2$) is a degenerate graph homomorphism,
we  prove that the image of $f$ is contained in a union of two maximal adjacent sets of different types.
For the case of finite fields, we obtain two better results on degenerate  graph homomorphisms.

\vspace{2mm}


\vspace{2mm}

\noindent{\bf 2010 AMS Classification}:  05C50, 15B33, 05C60,  51A10, 51D20

\end{abstract}

\section{Introduction}
 \setcounter{equation}{0}

\ \ \ \ \
The study of the {\em geometry of matrices} was initiated by Hua in the mid 1940s \cite{Hua8,W3}. The fundamental problem in the geometry of  matrices
is to characterize the transformation group of matrices by as few geometric invariants as possible. In the view of equivalence, the basic problem of the geometry of  matrices
is also to study graph isomorphisms on matrices.
In 1951,  Hua \cite{Hua8} proved the  fundamental theorem of the geometry of rectangular matrices over division rings. Hua's theorem also characterized the graph isomorphism on
 rectangular matrices, and his work was continued by many scholars
 (cf. \cite{Chooi}, \cite{LAA436Huang}-\cite{HWEN}, \cite{Lim-1}, \cite{Pazzis}-\cite{Wan5}).

In the algebraic graph theory, the research of graph homomorphisms is an important subject \cite{Godsil,Oxford}.
Thus, it is of significance to determine  graph homomorphisms on matrices.
Recently, literatures \cite{optimal, Huang-Semrl-2014,  Pazzis, Huang-Li,Huangli-I} discussed the graph homomorphisms on matrices over a division ring, and
 these work is interesting  for the geometry of matrices, the algebraic graph theory and  the preserver problems.
This paper is  continuation of literature \cite{Huangli-I}, and our goal is further to characterize the graph homomorphism on rectangular matrices.
 All mathematical symbols and definitions that are not explained, see  \cite{Huangli-I}.

Throughout this paper, we assume that
 $\mathbb{D}, \mathbb{D}'$ are division rings,  ${\mathbb{D}}^*=\mathbb{D}\setminus \{0\}$, and $m, n, m',n'$ are positive integers.
  Denote by  $\mathbb{F}_q$ the finite field
with $q$ elements where $q$ is a power of a prime.
Let $|X|$ be the cardinality  of a set $X$. Let ${\mathbb{D}}^n$ [resp. $^n{\mathbb{D}}$] denote the left [resp. right] vector space over
$\mathbb{D}$ whose elements are $n$-dimensional row [resp. column] vectors over  $\mathbb{D}$.
On the basic properties of matrices over a division ring, one may see literatures \cite{Thoms-Algebra,W3}.

Let ${\mathbb{D}}^{m\times n}$ and ${\mathbb{D}}^{m\times n}_r$ denote the sets of $m\times n$  matrices and  $m\times n$  matrices of rank $r$ over $\mathbb{D}$,
respectively. A matrix in ${\mathbb{D}}^{m\times n}$ is also called a point. Denote the set of $n\times n$ invertible matrices over $\mathbb{D}$ by $GL_n(\mathbb{D})$.
Let $I_r$ ($I$ for short) be the $r\times r$ identity matrix,  $0_{m,n}$ the $m\times n$ zero matrix ($0$ for short) and $0_n=0_{n,n}$.
Let   $E^{m\times n}_{ij}$ ($E_{ij}$ for short) denote the $m\times n$ matrix whose $(i,j)$-entry is $1$ and all other entries are $0$'s.
Denote by  $^tA$ the transpose matrix of a matrix $A$. If $\sigma: \mathbb{D}\rightarrow \mathbb{D}'$ is a map and  $A=(a_{ij})\in {\mathbb{D}}^{m\times n}$,
we write  $A^{\sigma}=(a_{ij}^{\sigma})$ and $^tA^\sigma=\,^t(A^\sigma)$.

Let $\Gamma(\mathbb{D}^{m\times n})$ be the graph whose vertex set is ${\mathbb{D}}^{m\times n}$ and
two vertices $A,B\in {\mathbb{D}}^{m\times n}$ are adjacent if ${\rm rank}(A-B)=1$. We write $A\sim B$ if ${\rm rank}(A-B)=1$.
The  $\Gamma(\mathbb{D}^{m\times n})$ is called  the graph on $m\times n$ matrices  over $\mathbb{D}$.
 $\Gamma(\mathbb{D}^{m\times n})$ is a connected distance transitive graph \cite{Godsil}.
When $\mathbb{D}=\mathbb{F}_q$, $\Gamma(\mathbb{F}_q^{m\times n})$ is also called a {\em bilinear forms graph}  \cite{Brouwera2}.
For $A, B\in  {\mathbb{D}}^{m\times n}$,  $d(A,B):={\rm rank}(A-B)$ is  the {\em distance} between $A$ and $B$ in $\Gamma(\mathbb{D}^{m\times n})$ (cf. \cite{Thoms-Algebra,W3}).
The  {\em diameter} of a subgraph $G$ of $\Gamma(\mathbb{D}^{m\times n})$, denoted by  ${\rm diam}(G)$, is the maximum distance between two distinct vertices in $G$.

Let $\varphi:  \mathbb{D}^{m\times n} \rightarrow  {\mathbb{D}'}^{m'\times n'}$ be a map.  The  $\varphi$ is called a {\em graph homomorphism}
if $A\sim B$ implies that $\varphi(A)\sim\varphi(B)$.
 The  $\varphi$ is called a {\em  distance preserving map} if $d(A, B)=d(\varphi(A),\varphi(B))$ for all $A,B\in {\mathbb{D}}^{m\times n}$.
The  $\varphi$ is called  a {\em  distance $k$ preserving map} if $d(A, B)=k$ implies that $d(\varphi(A),\varphi(B))=k$ for some fixed $k$.
 In the geometry of matrices, a graph homomorphism [resp. {\em graph isomorphism}] is also called an {\em adjacency preserving map}
[resp.  {\em adjacency preserving bijection in both directions}].
If $\varphi:  {\mathbb{D}}^{m\times n}\rightarrow  {\mathbb{D}'}^{m'\times n'}$ is a graph homomorphism, then
\begin{equation}\label{uy87mbmm}
\mbox{$d(\varphi(A), \varphi(B))\leq d(A, B)$, \ for all $A, B\in  {\mathbb{D}}^{m\times n}$.}
\end{equation}

A nonempty subset $\cal S$ of ${\mathbb{D}}^{m\times n}$ is called an {\em adjacent set} if any two distinct vertices (matrices) in $\cal S$ are adjacent.
 An adjacent set $\cal M$ in ${\mathbb{D}}^{m\times n}$ is called a {\em maximal adjacent set} ({\em maximal set} for short),
 if  there is no adjacent set in ${\mathbb{D}}^{m\times n}$ which properly contains $\cal M$ as a subset.
 In graph theory, a maximal  set is also called a {\em maximal clique} \cite{Brouwera2,Godsil}.
In  ${\mathbb{D}}^{m\times n}$, every adjacent set can be extended to a maximal set, and there are only two type of maximal sets.
For convenience, we think that a maximal set and its vertex set are equal.

Suppose that $\mathbb{D}$ and $\mathbb{D}'$ are division rings and $m,n,m',n'\geq 2$ are integers.  Write $E_{ij}=E_{ij}^{m\times n}$ and  $E_{ij}'=E_{ij}^{m'\times n'}$.
 For $1\leq i\leq m$ and $1\leq j\leq n$ (or $1\leq i\leq m'$ and $1\leq j\leq n'$), we let
\begin{equation}\label{maximalsets56463gd}
{\cal M}_i=\left\{\sum_{j=1}^{n}x_jE_{ij}: x_j\in \mathbb{D}\right \}, \ \ {\cal N}_j=\left\{\sum_{i=1}^{m} y_iE_{ij}: y_i\in \mathbb{D}\right \}\subset {\mathbb{D}}^{m\times n};
\end{equation}
\begin{equation}\label{maxiwq42hfh463}
{\cal M}_i'=\left\{\sum_{j=1}^{n'}x_jE_{ij}': x_j\in \mathbb{D}'\right \}, \ \ {\cal N}_j'=\left\{\sum_{i=1}^{m'} y_iE_{ij}': y_i\in \mathbb{D}'\right \}\subset {\mathbb{D}'}^{m'\times n'}.
\end{equation}

\begin{lem}\label{Rectangular-PID2-4}{\rm (cf. \cite[Prpopsition 3.8]{W3})} \
 Every maximal set  ${\cal M}$ of ${\mathbb{D}}^{m\times n}$ is one of  the following  forms.

{Type one.} \  ${\cal M}=P{\cal M}_1Q+A=P{\cal M}_1+A$, where $P\in GL_m(\mathbb{D})$, $Q\in GL_n(\mathbb{D})$ and $A\in {\mathbb{D}}^{m\times n}$.

{Type two.} \ ${\cal M}=P{\cal N}_1Q+A={\cal N}_1Q+A$, where $P\in GL_m(\mathbb{D})$, $Q\in GL_n(\mathbb{D})$ and  $A\in {\mathbb{D}}^{m\times n}$.
\end{lem}

For $1\leq k\leq {\rm min}\{m,n\}$ and $A\in {\mathbb{D}}^{m\times n}$, we let
$$\mathbb{D}_{\leq \,k}^{m\times n}=\left\{X\in\mathbb{D}^{m\times n}: {\rm rank}(X)\leq k \right\},$$
$$\mathbb{B}_{A}=\mathbb{D}_{\leq 1}^{m\times n}+A,$$
$$\mathbb{N}_{A}=\mathbb{D}_{1}^{m\times n}+A\subset\mathbb{B}_{A}.$$
The  $\mathbb{B}_{A}$ is called  the {\em unit ball} with a central point $A$, and the $\mathbb{N}_{A}$ is called
the {\em neighbourhood} of $A$.
 If $\varphi:  {\mathbb{D}}^{m\times n} \rightarrow  {\mathbb{D}'}^{m'\times n'}$ is a graph homomorphism, then
$\varphi(\mathbb{N}_{A})\subseteq \mathbb{N}_{\varphi(A)}$ and  $\varphi(\mathbb{B}_{A})\subseteq \mathbb{B}_{\varphi(A)}$  for every $A\in {\mathbb{D}}^{m\times n}$.

Let $\varphi: {\mathbb{D}}^{m\times n}\rightarrow{\mathbb{D}'}^{m'\times n'}$ be a graph homomorphism. The homomorphism $\varphi$  is called  {\em degenerate},
 if there exists a matrix $A\in\mathbb{D}^{m\times n}_{\leq 1}$ and there are two  maximal sets $\mathcal{M}$ and $\mathcal{N}$ of
 different types in ${\mathbb{D}'}^{m'\times n'}$, such that $\varphi\left(\mathbb{B}_{A}\right)\subseteq \mathcal{M}\cup \mathcal{N}$
 with $\varphi(A)\in\mathcal{M}\cap\mathcal{N}$.
 The homomorphism $\varphi$ is called {\em non-degenerate} if it is not degenerate.
The homomorphism $\varphi$ is called a (vertex) {\em colouring} if $\varphi(\mathbb{D}^{m\times n})$ is an adjacent set in ${\mathbb{D}'}^{m'\times n'}$.

Every (vertex) colouring is a degenerate graph homomorphism.
It is easy to see that $\varphi$ is a (vertex) colouring if and only if ${\rm diam}(\varphi(\mathbb{D}^{m\times n}))=1$.
The word ``colouring" is derived from graph theory \cite{Godsil}.
Any distance preserving map is a non-degenerate graph homomorphism but not vice versa.

This  paper is organized as follows. In Section 2, we  introduce  some Lemmas on maximal sets.
In Section 3, We determine  additive graph homomorphisms from ${\mathbb{D}}^{m\times n}$ to ${\mathbb{D}'}^{m'\times n'}$.
 In Section 4, we will characterize the graph homomorphism $\varphi$ from ${\mathbb{D}}^{n\times n}$ to ${\mathbb{D}'}^{m'\times n'}$ (where $|\mathbb{D}|\geq 4$)
if there exists $A_0\in {\mathbb{D}}^{n\times n}$ such that ${\rm rank}(\varphi(A_0))=n$, this result extends \cite[Theorem 4.1]{optimal} to
the cases of two division rings and $n=2$.
In Section 5,  we  discuss the degenerate graph homomorphisms.  Let $f:  {\mathbb{D}}^{m\times n}\rightarrow  {\mathbb{D}'}^{m'\times n'}$ (where $|\mathbb{D}|\geq 4$)
be a  degenerate graph homomorphism. We prove that $f\left(\mathbb{D}_{\leq 1}^{m\times n} \right)$ and $f(\mathbb{B}_{A_0})$ are two adjacent sets,
if there exists $A_0\in {\mathbb{D}}^{m\times n}$ such that ${\rm rank}(f(A_0))={\rm min}\{m,n\}$. Moreover, when ${\rm min}\{m,n\}=2$, the image of $f$ is contained in
a union of two maximal  sets of different types. For the case of finite fields, we obtain two better results on the degenerate graph homomorphisms.

\section{Lemmas on maximal sets}

\ \ \ \ \  In this section, we will introduce maximal sets on ${\mathbb{D}}^{m\times n}$ and their affine geometries.

\begin{lem}\label{Matrix-PID5-4bb}{\rm (cf. \cite[Lemma 3.2]{LAA436Huang})} \ Let  $m,n,r,s$ be integers with $1\leq r, s<{\rm min}\{m, n\}$.
Assume that $\alpha=\{i_1, \ldots, i_r\}$, $\beta=\{j_1, \ldots,
j_s\}$, where $1\leq i_1<\cdots <i_r\leq m$ and $1\leq j_1<\cdots <j_s\leq n$. Let $A=(a_{ij})\in \mathbb{D}^{m \times n}$,
$B_i=\sum_{t=1}^r\sum_{k=1}^sb^{(i)}_{i_tj_k}E_{i_tj_k}\in \mathbb{D}^{m \times n}$ (where $b^{(i)}_{i_tj_k}\in \mathbb{D}$), $i=1, 2$, and $B_1\neq B_2$.
If $A\sim B_i$, $i=1, 2$, then either $a_{ij}=0$ for all $i\notin \alpha$, or  $a_{ij}=0$ for all $j\notin \beta$.
\end{lem}

Using Lemmas \ref{Matrix-PID5-4bb} and \ref{Rectangular-PID2-4}, we can prove the following results.

\begin{cor}\label{Rectangular-PID2-7}{\rm (cf. \cite[Corollary 3.10]{W3})}
 Let $A$ and $B$ be two  adjacent points in ${\mathbb{D}}^{m\times n}$. Then there are exactly two maximal sets
 $\mathcal{M}$ and $\mathcal{M}'$ containing $A$ and $B$. Moreover, $\mathcal{M}$ and $\mathcal{M}'$ are of different types.
\end{cor}

\begin{cor}\label{Rectangular-1-13}{\rm (cf. \cite{LAA436Huang,Huang-book,W3})}
If $\mathcal{M}$ and $\mathcal{N}$ are two distinct maximal sets of the same type $[$resp.  different types$]$ in ${\mathbb{D}}^{m\times n}$ with
$\mathcal{M}\cap\mathcal{N}\neq \emptyset$, then $|\mathcal{M}\cap\mathcal{N}|=1$ $[$resp. $|\mathcal{M}\cap \mathcal{N}|\geq 2$$]$.
\end{cor}

\begin{lem}\label{maximalset022}{\rm (cf. \cite[Lemma 3.4]{Huangli-I})} Suppose that $\mathcal{M}$, $\mathcal{N}$ are two distinct maximal sets in  ${\mathbb{D}}^{m\times n}$
with $\mathcal{M}\cap\mathcal{N}\neq \emptyset$. Then:

{\rm (i)} \  if $\mathcal{M}$, $\mathcal{N}$ are of different types, then for any $A\in\mathcal{M}\cap\mathcal{N}$,
there are $P\in GL_m(\mathbb{D})$ and $Q\in GL_n(\mathbb{D})$
such that ${\cal M}=P{\cal M}_1Q+A=P{\cal M}_1+A$ and  ${\cal N}=P{\cal N}_1Q+A={\cal N}_1Q+A$;

{\rm (ii)} \  if both $\mathcal{M}$ and $\mathcal{N}$ are of type one $[$resp. type two$]$, then there exists an invertible  matrix $P$ $[$resp. $Q$$]$
such that ${\cal M}=P{\cal M}_1+A$ and ${\cal N}=P{\cal M}_2+A$
$[$resp. ${\cal M}={\cal N}_1Q+A$ and ${\cal N}={\cal N}_2Q+A$$]$, where $\mathcal{M}\cap\mathcal{N} = \{A \}$.
\end{lem}

 There is an  axiomatic definition of the {\em affine geometry} (cf. \cite{Bennett, ModernProjective}). Let $V$ be
 an $r$-dimensional  left vector subspace  of ${\mathbb{D}}^n$ and $a\in\mathbb{D}^n$. Then $V+a$ is called an $r$-dimensional {\em left affine flat} ({\em flat} for short)
over $\mathbb{D}$.  When $r\geq 2$, the set of all flats in $V+a$ is called the {\em left  affine geometry} on $V+a$, which is denoted by $AG(V+a)$.  The {\em dimension} of $AG(V+a)$ is $r$,
 denoted by ${\rm dim}(AG(V+a))=r$. The flats of dimensions $0, 1, 2$  are called {\em points, lines, planes} in $AG(V+a)$.
Similarly, we have the {\em right  affine geometry} on a right affine flat over $\mathbb{D}$.

Let ${\cal M}=P{\cal M}_1+A$ be a maximal set of type one in $\mathbb{D}^{m\times n}$, where $P\in GL_m(\mathbb{D})$ and $A\in\mathbb{D}^{m\times n}$. Then we have
a left affine geometry $AG(P{\cal M}_1+A)$ such that $AG(P{\cal M}_1+A)$ and $AG({\mathbb{D}}^n)$ are affine isomorphic.
Similarly,  we have a right affine geometry $AG({\cal N}_1Q+A)$ such that $AG({\cal N}_1Q+A)$ and $AG(^m{\mathbb{D}})$ are affine isomorphic,
where $Q\in GL_n(\mathbb{D})$ (cf. \cite{Huangli-I}).
In $AG(P{\cal M}_1+A)$, the parametric equation of a line $\ell$ is
$$\mbox{$\ell=P\left(
\begin{array}{cc}
\mathbb{D}\alpha+\beta \\
0\\
\end{array}\right)+A$, \ where $\alpha, \beta\in \mathbb{D}^n$ with $\alpha\neq 0$.}
$$

\begin{lem}\label{Rectangular-PID2-11}{\rm (cf. \cite[p.95]{W3})} \ Let ${\cal M}$
be a maximal set in $\mathbb{D}^{m\times n}$. Then $\ell$ is a line in $AG({\cal M})$ if and only if
$\ell=\mathcal{M}\cap \mathcal{N}$, where $\mathcal{M}$ and $\mathcal{N}$ are two  maximal sets of different types with $\mathcal{M}\cap\mathcal{N}\neq\emptyset$.
Moreover,  $|\ell|=|\mathcal{M}\cap\mathcal{N}|=|\mathbb{D}|$.
\end{lem}

\begin{lem}\label{non-degeneratelemma00a}{\rm(see \cite[Lemma 4.6]{Huangli-I})} \ Let $m,n,m', n'\geq 2$ be integers. Suppose that
$\varphi:  {\mathbb{D}}^{m\times n}\rightarrow  {\mathbb{D}'}^{m'\times n'}$ is a non-degenerate graph homomorphism with $\varphi(0)=0$. Then
\begin{itemize}
\item[{\rm (a)}]   if ${\cal M}$ is a maximal set of type one [resp. type two]  containing $0$ in ${\mathbb{D}}^{m\times n}$ and $\varphi({\cal M})\subseteq {\cal M}'$, where
 ${\cal M}'$ is a maximal set  in ${\mathbb{D}'}^{m'\times n'}$, then for any $A\in{\cal M}$,  there exists a  maximal set ${\cal R}$
 of type one [resp. type two]  in ${\mathbb{D}}^{m\times n}$, such that $\varphi({\cal R})\subseteq {\cal R}'$ and ${\cal R}\cap {\cal M}=\{A\}$,
where ${\cal R}'$ is a maximal set  in ${\mathbb{D}'}^{m'\times n'}$, ${\cal R}'$ and ${\cal M}'$ are of the same type with ${\cal R}'\neq {\cal M}'$;

\item[{\rm (b)}] if ${\cal M}$ and ${\cal N}$ are two distinct maximal sets of the same type  [resp. different types]
in ${\mathbb{D}}^{m\times n}$ such that $0\in \mathcal{M}$ and ${\cal M}\cap {\cal N}\neq \emptyset$, then  $\varphi({\cal M})\subseteq {\cal M}'$
and $\varphi({\cal N})\subseteq {\cal N}'$,
where ${\cal M}'$ and ${\cal N}'$ are two maximal sets of the same type  [resp. different types]
in ${\mathbb{D}'}^{m'\times n'}$;

\item[{\rm (c)}] if ${\cal M}$ is a maximal set containing $0$ in ${\mathbb{D}}^{m\times n}$  and  $\varphi(\mathcal{M})\subseteq \mathcal{M}'$
where ${\cal M}'$ is a maximal set in ${\mathbb{D}'}^{m'\times n'}$,
then ${\cal M}'$ is the unique  maximal set containing $\varphi({\cal M})$, and $\varphi({\cal M})$ is not contained in any line of $AG({\cal M}')$.
\end{itemize}
\end{lem}

Let $\mathcal{S}$ be an adjacent set in ${\mathbb{D}}^{m\times n}$ with $0\in \mathcal{S}$ and $|\mathcal{S}|\geq 2$. Then there exists a maximal set $\mathcal{M}$  containing $\mathcal{S}$.
By Lemma \ref{Rectangular-PID2-4}, $\mathcal{M}=P\mathcal{M}_1$ or $\mathcal{M}=\mathcal{N}_1Q$,  where $P$ and $Q$ are invertible matrices. Hence
$P^{-1}\mathcal{S}\subseteq\mathcal{M}_1$ or $\mathcal{S}Q^{-1}\subseteq\mathcal{N}_1$. When $P^{-1}\mathcal{S}\subseteq\mathcal{M}_1$
[resp. $\mathcal{S}Q^{-1}\subseteq\mathcal{N}_1$],
the {\em dimension} of $\mathcal{S}$,  denoted by ${\rm dim}(\mathcal{S})$,  is the
number of matrices in a maximal left [resp. right] linear independent subset  of $P^{-1}\mathcal{S}$ [resp. $\mathcal{S}Q^{-1}$].
The ${\rm dim}(\mathcal{S})$ is uniquely determined by $\mathcal{S}$ and ${\rm dim}(\mathcal{S})\leq {\rm max}\{m,n\}$.

\begin{lem}\label{degenerate-2}{\rm (see \cite[Lemma 4.13]{Huangli-I})} \
Let $\mathbb{D}, \mathbb{D}'$ be division rings with $|\mathbb{D}|\geq 4$. Suppose  $\varphi:  {\mathbb{D}}^{m\times n}\rightarrow  {\mathbb{D}'}^{m'\times n'}$
is a non-degenerate graph homomorphism with $\varphi(0)=0$. Let ${\cal M}$ be a maximal set containing $0$ in ${\mathbb{D}}^{m\times n}$.
If $\varphi({\cal M})\subseteq {\cal M}'$ where ${\cal M}'$ is a maximal set in ${\mathbb{D}'}^{m'\times n'}$,
then the restriction map $\varphi\mid_{\mathcal{M}}: \mathcal{M}\rightarrow \mathcal{M}'$
 is an injective weighted semi-affine map. Moreover, if $\mathcal{S}$ is an adjacent set in ${\mathbb{D}}^{m\times n}$ with  $0\in \mathcal{S}$
 and $|\mathcal{S}|\geq 2$, then $${\rm dim}(\varphi(\mathcal{S}))\leq {\rm dim}(\mathcal{S}).$$
\end{lem}

\begin{lem}\label{gdkj7ssb33} Let $\mathcal{M}$  be a maximal set in $\mathbb{D}^{m\times n}$ $(m,n\geq 2)$ and $A\in \mathbb{D}^{m\times n}$. Then
$A\in \mathcal{M}$ if and only if $A$ is adjacent with three  noncollinear points in $AG(\mathcal{M})$.
\end{lem}
\proof
If $A\in \mathcal{M}$, it is clear that $A$ is adjacent with three  noncollinear points in $AG(\mathcal{M})$.
Now,  suppose that  $A$ is adjacent with three  noncollinear points $B_1,B_2,B_3$ in $AG(\mathcal{M})$.
 Then, there are $P\in GL_m(\mathbb{D})$ and $Q\in GL_n(\mathbb{D})$ such that
$B_2-B_1=PE_{11}Q$. After by the map $X\mapsto P^{-1}(X-B_1)Q^{-1}$, we may assume with no loss of  generality that $B_1=0$ and $B_2=E_{11}$.
By Corollary \ref{Rectangular-PID2-7}, we have either $\mathcal{M}=\mathcal{M}_1$ or $\mathcal{M}=\mathcal{N}_1$.
Without loss of generality  we assume that $\mathcal{M}=\mathcal{M}_1$.
Clearly, $\ell=\{xE_{11}: x\in \mathbb{D}\}$ is a line in
$AG(\mathcal{M}_1)$ containing $0$ and $E_{11}$. Since  $0,E_{11},B_3$ are  noncollinear points in $AG(\mathcal{M}_1)$, $B_3$ is of the form
$\scriptsize B_3=\left(
  \begin{array}{cc}
    b & \beta\\
    0 & 0 \\
  \end{array}
\right)$ where $0\neq \beta\in \mathbb{D}^{n-1}$.  Write $\scriptsize A=\left(
  \begin{array}{cc}
    a_{11} & A_{12}\\
    A_{21} & A_{22} \\
  \end{array}
\right)$ where $a_{11}\in \mathbb{D}$. By $A\sim 0$, $A\sim E_{11}$ and
$\scriptsize A\sim \left(
  \begin{array}{cc}
    b & \beta\\
    0 & 0 \\
  \end{array}
\right)$, we get
$${\rm rank}\left(
  \begin{array}{cc}
    a_{11} & A_{12}\\
    A_{21} & A_{22} \\
  \end{array}
\right)={\rm rank}\left(
  \begin{array}{cc}
    a_{11}-1 & A_{12}\\
    A_{21} & A_{22} \\
  \end{array}
\right)={\rm rank}\left(
  \begin{array}{cc}
    a_{11}-b & A_{12}-\beta\\
    A_{21} & A_{22} \\
  \end{array}
\right)=1.$$
It follows from Lemma \ref{Matrix-PID5-4bb} that $(A_{21},A_{22})=0$, and hence $A\in \mathcal{M}_1$.
$\qed$

\section{Additive graph homomorphisms on rectangular matrices }

\ \ \ \ \ \
Let $\phi: {\mathbb{D}}^{m\times n}\rightarrow  {\mathbb{D}'}^{m'\times n'}$ be a map. The  $\phi$ is called an {\em additive map}  if $\phi(A+B)=\phi(A)+\phi(B)$
for all $A,B\in \mathbb{D}^{m\times n}$. The  $\phi$ is called an {\em additive graph homomorphism} if $\phi$
is a graph homomorphism and an additive map. If $\phi$ is an  additive graph homomorphism, then  $\phi(0)=0$ and
$\phi(A-B)=\phi(A)-\phi(B)$ for all $A,B\in\mathbb{D}^{m\times n}$.

The additive preserver problems on matrices is an active research area in linear algebra and the theory of matrices \cite{X.Zhang}.
In the additive preserver problems,  an additive graph homomorphism is also called an {\em additive rank-$1$ preserving map}.
When $\mathbb{D}=\mathbb{D}'$ is a field, additive rank-$1$ preserving maps from  ${\mathbb{D}}^{m\times n}$ to ${\mathbb{D}}^{m'\times n'}$ is determined by \cite{X.Zhang}.
In this section, we discuss the case of division rings, and our main result is the following theorem.

\begin{thm}\label{additivehomomorphism0a}   Suppose $\mathbb{D}, \mathbb{D}'$ are division rings
and $m,n,m',n'\geq 2$ are integers. Let
$\phi:  {\mathbb{D}}^{m\times n}\rightarrow  {\mathbb{D}'}^{m'\times n'}$
be an additive graph homomorphism. Then either there exist  matrices  $P\in {\mathbb{D}'}^{m'\times m}$ and $Q\in {\mathbb{D}'}^{n\times n'}$ with their ranks $\geq 2$,
and a nonzero ring homomorphism $\tau: \mathbb{D} \rightarrow\mathbb{D}'$,  such that
\begin{equation}\label{adhf57ni665fdc0}
 \phi(X)=P X^\tau Q, \ \, X\in {\mathbb{D}}^{m\times n};
\end{equation}
or there exist matrices  $P\in {\mathbb{D}'}^{m'\times n}$ and $Q\in {\mathbb{D}'}^{m\times n'}$ with their ranks $\geq 2$,
and a nonzero ring anti-homomorphism $\sigma: \mathbb{D}\rightarrow\mathbb{D}'$,  such that
\begin{equation}\label{addn6gkkglad53ss75}
\phi(X)=P \,{^tX^\sigma}Q, \ \,  X\in {\mathbb{D}}^{m\times n};
\end{equation}
or $\phi$ is  a (vertex) colouring.
\end{thm}

In order to prove Theorem \ref{additivehomomorphism0a},  we need  the following results.

\begin{thm}\label{additiverank}
Suppose $\mathbb{D}, \mathbb{D}'$ are division rings  and  $m,n,m', n'\geq 2$ are integers.
 Let  $\phi: {\mathbb{D}}^{m\times n}\rightarrow  {\mathbb{D}'}^{m'\times n'}$
be an additive graph homomorphism. If $\phi$ is  degenerate, then $\phi$ is  a (vertex) colouring.
 \end{thm}
\proof
Assume that $\phi$ is  degenerate. Since $\phi$ is additive,  there are two maximal sets $\mathcal{M}, \mathcal{N}$ of
 different types in ${\mathbb{D}'}^{m'\times n'}$, such that $\phi\left(\mathbb{D}_{\leq 1}^{m\times n}\right)\subseteq \mathcal{M}\cup \mathcal{N}$
 with $0\in\mathcal{M}\cap\mathcal{N}$. We show that $\phi\left(\mathbb{D}_{\leq 1}^{m\times n}\right)$ is an adjacent set by contradiction as follows.

Suppose  $\phi\left(\mathbb{D}_{\leq 1}^{m\times n}\right)$ is not an adjacent set.
Then there are matrices $A,B\in \mathbb{D}^{m\times n}_1$ such that ${\rm rank}(A-B)={\rm rank}(\phi(A)-\phi(B))=2$. Let $E_{ij}=E_{ij}^{m\times n}$ and $E_{ij}'=E_{ij}^{m'\times n'}$.
There are invertible matrices $P_1, Q_1$ over $\mathbb{D}$ such that $A=P_1E_{11}Q_1$ and $B=P_1E_{22}Q_1$.
Also, there are invertible matrices $P_2, Q_2$ over $\mathbb{D}'$ such that $P_2\phi(A)Q_2=E_{11}'$ and $P_2\phi(B)Q_2=E_{22}'$.
Replacing $\phi$ by the map $X\longmapsto P_2\phi(P_1XQ_1)Q_2$, we have
\begin{equation}\label{nbvwer42ftr6sdg}
\phi(E_{ii})=E_{ii}', \ \ i=1,2.
\end{equation}
By Corollary \ref{Rectangular-PID2-7}, there are exactly two maximal sets containing $E_{ii}$ and $0$ [resp. $E_{ii}'$ and $0$], they are $\mathcal{M}_i$ and  $\mathcal{N}_i$
[resp. $\mathcal{M}_i'$ and  $\mathcal{N}_i'$] ($i=1,2$). Thus, by (\ref{nbvwer42ftr6sdg}) and $\phi(0)=0$, we get that
 $\phi(\mathcal{M}_1)\subseteq \mathcal{M}_1'$ or $\phi(\mathcal{M}_1)\subseteq \mathcal{N}_1'$. Without loss of generality, we assume  $\phi(\mathcal{M}_1)\subseteq \mathcal{M}_1'$.
 From the above discussion we have  $\phi(\mathcal{N}_2)\subseteq \mathcal{N}_2'$ or  $\phi(\mathcal{N}_2)\subseteq \mathcal{M}_2'$.
 Since  $\mathcal{M}_1\cap\mathcal{N}_2=\mathbb{D}E_{12}$ and $\mathcal{M}_1'\cap\mathcal{M}_2'=\{0\}$, we must have  $\phi(\mathcal{N}_2)\subseteq \mathcal{N}_2'$.
 By   $\mathcal{M}_1'\cap\mathcal{N}_2'=\mathbb{D}'E_{12}'$, we obtain $\phi(\lambda E_{12})=\lambda^*E_{12}'$ for all $\lambda\in \mathbb{D}$.
 Since $\phi$ is additive,
 $$\phi(-E_{11}-E_{12}-E_{22})= -E_{11}'+(-1)^*E_{12}'-E_{22}'.$$
Recall that  $\phi\left(\mathbb{D}_{\leq 1}^{m\times n}\right)\subseteq \mathcal{M}\cup \mathcal{N}$.
Using Corollary \ref{Rectangular-PID2-7}, it is clear that
$\phi\left(\mathbb{D}_{\leq 1}^{m\times n}\right)\subseteq \mathcal{M}_1'\cup \mathcal{N}_2'$.
Since $E_{21}\sim (-E_{11}-E_{12}-E_{22})$,  $\phi(E_{21})\sim (-E_{11}'+(-1)^*E_{12}'-E_{22}')$.  On the other hand, by $\phi(E_{21})\in \mathcal{M}_1'\cup \mathcal{N}_2'$ and
$\phi(E_{21})\sim E_{ii}'$ ($i=1,2$), we must have $\phi(E_{21})=bE_{12}'$ where $b\in{\mathbb{D}'}^*$. Hence
$${\rm rank}\left(\phi(E_{21})- (-E_{11}'+(-1)^*E_{12}'-E_{22}')\right)=2,$$
a contradiction.

Thus, $\phi\left(\mathbb{D}_{\leq 1}^{m\times n}\right)$ is an adjacent set. Then there is a fixed maximal set $\mathcal{M}'$ containing $0$ in ${\mathbb{D}'}^{m'\times n'}$,
such that $\phi\left(\mathbb{D}_{\leq 1}^{m\times n}\right)\subseteq \mathcal{M}'$. Note that $Y_1,Y_2\in\mathcal{M}'$ implies that
 $Y_1+ Y_2 \in\mathcal{M}'$.
Since $\phi$ is additive, it is clear that
$\phi\left(\mathbb{D}^{m\times n}\right)\subseteq \mathcal{M}'$, and hence $\phi$ is  a (vertex) colouring.
$\qed$

\begin{lem}\label{s4245549bcnfte}
Suppose $\mathbb{D}, \mathbb{D}'$ are division rings  and  $m,n, m',n'\geq 2$ are
integers. Let  $\varphi:  {\mathbb{D}}^{m\times n}\rightarrow  {\mathbb{D}'}^{m'\times n'}$
be a non-degenerate graph homomorphism with $\varphi(0)=0$.
Then  there exist two invertible matrices  $T_1$ and $T_2$ over $\mathbb{D}'$, and two invertible matrices  $\scriptsize T_3^{-1}=\left(
\begin{array}{cc}
 1 & * \\
 0 & * \\
 \end{array}
 \right)\in GL_m(\mathbb{D})$ and  $\scriptsize T_4^{-1}=\left(
\begin{array}{cc}
 1 & 0 \\
 * & * \\
 \end{array}
 \right)\in GL_n(\mathbb{D})$,  such that
either
$$\mbox{$\varphi(T_3^{-1}\mathcal{M}_iT_4^{-1})\subseteq T_1\mathcal{M}_i'T_2$ \ with
\ $\varphi(T_3^{-1}\mathcal{N}_iT_4^{-1})\subseteq T_1\mathcal{N}_i'T_2$ $(i=1,2)$,}$$
or
$$\mbox{$\varphi(T_3^{-1}\mathcal{M}_iT_4^{-1})\subseteq T_1\mathcal{N}_i'T_2$ \ with \ $\varphi(T_3^{-1}\mathcal{N}_iT_4^{-1})\subseteq T_1\mathcal{M}_i'T_2$ $(i=1,2)$.}$$
\end{lem}
\proof
{\em Case 1}. \ $\varphi(\mathcal{M}_1)$ is contained in a maximal set of type one. Then there exists an invertible matrix $T_0$ over $\mathbb{D}'$
such that $\varphi(\mathcal{M}_1)\subseteq T_0\mathcal{M}_1'$.
By Lemma \ref{non-degeneratelemma00a}(a), there exists a type one maximal set ${\cal M}$ in ${\mathbb{D}}^{m\times n}$ containing $0$,
such that ${\cal M}\neq {\cal M}_1$ and  $\varphi({\cal M})\subseteq {\cal M}'$,
where ${\cal M}'$ is a type one maximal set  in ${\mathbb{D}'}^{m'\times n'}$ containing $0$ and ${\cal M}'\neq T_0{\cal M}_1'$.
Using Lemma \ref{maximalset022}, there exists  $T_1\in GL_{m'}(\mathbb{D}')$ such that $T_0\mathcal{M}_1'=T_1\mathcal{M}_1'$
and $\mathcal{M}'=T_1\mathcal{M}_2'$. Thus, $\varphi(\mathcal{M}_1)\subseteq T_1\mathcal{M}_1'$ and $\varphi(\mathcal{M})\subseteq T_1\mathcal{M}_2'$.
Since
$$\mbox{$\mathcal{M}=\left\{\left(\begin{array}{c}
a_1x \\
a_2x \\
\vdots\\
a_m x\\
\end{array}\right): x\in \mathbb{D}^n
\right\}$ \ with \ $\left(
                    \begin{array}{c}
                      a_2 \\
                      \vdots \\
                      a_m \\
                    \end{array}
                  \right)
\neq 0$,}$$
 there is
$\scriptsize T_3^{-1}=\left(
\begin{array}{cc}
 1 & * \\
 0 & * \\
 \end{array}
 \right)\in GL_m(\mathbb{D})$ such that $\mathcal{M}=T_3^{-1}\mathcal{M}_2$ and $\mathcal{M}_1=T_3^{-1}\mathcal{M}_1$. Thus
 $$\mbox{$\varphi(T_3^{-1}\mathcal{M}_i)\subseteq T_1\mathcal{M}_i'$, \ $i=1,2$.}$$

Similarly, there exists a type two maximal set ${\cal N}$ in ${\mathbb{D}}^{m\times n}$ containing $0$,
 such that $\mathcal{N}\neq \mathcal{N}_1$,
 $\varphi(\mathcal{N}_1)\subseteq \mathcal{N}_1'T_2=T_1\mathcal{N}_1'T_2$ and $\varphi(\mathcal{N})\subseteq \mathcal{N}_2'T_2=T_1\mathcal{N}_2'T_2$,
where $T_2\in GL_{n'}(\mathbb{D}')$. Moreover, there is
$\scriptsize T_4^{-1}=\left(
\begin{array}{cc}
 1 & 0 \\
 * & * \\
 \end{array}
 \right)\in GL_n(\mathbb{D})$ such that $\mathcal{N}=\mathcal{N}_2T_4^{-1}=T_3^{-1}\mathcal{N}_2T_4^{-1}$ and
 $\mathcal{N}_1=\mathcal{N}_1T_4^{-1}=T_3^{-1}\mathcal{N}_1T_4^{-1}$.
It follows that $\varphi(T_3^{-1}\mathcal{N}_iT_4^{-1})\subseteq T_1\mathcal{N}_i'T_2$,  $i=1,2$.
Clearly, $T_3^{-1}\mathcal{M}_i=T_3^{-1}\mathcal{M}_iT_4^{-1}$ and  $T_1\mathcal{M}_i'=T_1\mathcal{M}_i'T_2$, $i=1,2$. Therefore, we obtain
$$\mbox{$\varphi(T_3^{-1}\mathcal{M}_iT_4^{-1})\subseteq T_1\mathcal{M}_i'T_2$, \ $i=1,2$.}$$

{\em Case 2}. \  $\varphi(\mathcal{M}_1)$ is contained in a maximal set of type two. Similarly,  there exist two invertible matrices  $T_1$ and $T_2$ over
$\mathbb{D}'$, and two invertible matrices  $\scriptsize T_3^{-1}=\left(
\begin{array}{cc}
 1 & * \\
 0 & * \\
 \end{array}
 \right)\in GL_m(\mathbb{D})$ and  $\scriptsize T_4^{-1}=\left(
\begin{array}{cc}
 1 & 0 \\
 * & * \\
 \end{array}
 \right)\in GL_n(\mathbb{D})$,  such that $\varphi(T_3^{-1}\mathcal{M}_iT_4^{-1})\subseteq T_1\mathcal{N}_i'T_2$  with
  $\varphi(T_3^{-1}\mathcal{N}_iT_4^{-1})\subseteq T_1\mathcal{M}_i'T_2$, $i=1,2$.
$\qed$

\begin{thm}\label{additive-non-degenerate1}  Suppose $\mathbb{D}, \mathbb{D}'$ are division rings and  $m',n'\geq 2$ are integers. Let
$\phi:  {\mathbb{D}}^{2\times 2}\rightarrow  {\mathbb{D}'}^{m'\times n'}$
be a non-degenerate additive graph homomorphism. Then either there exist  matrices  $P\in GL_{m'}(\mathbb{D}')$ and $Q\in GL_{n'}(\mathbb{D}')$,
and a nonzero ring homomorphism $\tau: \mathbb{D} \rightarrow\mathbb{D}'$,  such that
\begin{equation}\label{additiuytnv86h}
\phi(X)=P\left(\begin{array}{cc}
 X^\tau  & 0 \\
  0 & 0 \\
   \end{array}
  \right)Q, \ \, X\in {\mathbb{D}}^{2\times 2};
\end{equation}
or there exist matrices  $P\in GL_{m'}(\mathbb{D}')$ and $Q\in GL_{n'}(\mathbb{D}')$, and
a nonzero ring anti-homomorphism $\sigma: \mathbb{D}\rightarrow\mathbb{D}'$,  such that
\begin{equation}\label{additive-ere34nn}
\phi(X)=P\left(\begin{array}{cc}
{^tX^\sigma} & 0 \\
0 & 0 \\
\end{array}
\right)Q, \ \,  X\in {\mathbb{D}}^{2\times 2}.
\end{equation}
\end{thm}
\proof
By Lemma \ref{s4245549bcnfte}, we may assume with no loss of generality  either
$\phi(\mathcal{M}_i)\subseteq \mathcal{M}_i'$ with $\phi(\mathcal{N}_i)\subseteq \mathcal{N}_i'$ $(i=1,2)$,
or $\phi(\mathcal{M}_i)\subseteq \mathcal{N}_i'$  with  $\phi(\mathcal{N}_i)\subseteq \mathcal{M}_i'$ $(i=1,2)$.
We prove this result only for the first case; the second case is similar.
From now on we assume that
\begin{equation}\label{additive-fsdf5353}
\mbox{$\phi(\mathcal{M}_i)\subseteq \mathcal{M}_i'$ \ \ with \ \ $\phi(\mathcal{N}_i)\subseteq \mathcal{N}_i'$, \ $i=1,2$.}
\end{equation}
We show that $\phi$ is of the form (\ref{additiuytnv86h}) as follows. (For the second case, we can prove similarly $\phi$ is of the form (\ref{additive-ere34nn}).)

Let $E_{ij}=E_{ij}^{2\times 2}$ and $E_{ij}'=E_{ij}^{m'\times n'}$. By (\ref{additive-fsdf5353}) we have that $\phi(xE_{ij})=x^{\sigma_{ij}}E_{ij}'$ for all $x\in \mathbb{D}$, $i,j=1,2$,
where $\sigma_{ij}: \mathbb{D}\rightarrow \mathbb{D}'$ is an additive map with $0^{\sigma_{ij}}=0$. Since $\phi$ is additive, we get
$$\phi\left(
        \begin{array}{cc}
          x_{11} & x_{12}  \\
          x_{21}  & x_{22}  \\
        \end{array}
      \right)=\left(
                \begin{array}{cc}
                  \left(
        \begin{array}{cc}
          x_{11}^{\sigma_{11}} & x_{12}^{\sigma_{12}}  \\
          x_{21}^{\sigma_{21}}  & x_{22}^{\sigma_{22}}  \\
        \end{array}
      \right) & 0 \\
                  0 & 0 \\
                \end{array}
              \right), \ x_{ij}\in \mathbb{D}.$$
Without loss of generality, we assume that $1^{\sigma_{11}}=1^{\sigma_{22}}=1$.

Since $\scriptsize\left(
        \begin{array}{cc}
          x & x \\
        1& 1 \\
        \end{array}
      \right)\sim 0$ and $\scriptsize\left(
        \begin{array}{cc}
        1 & x \\
        1& x \\
        \end{array}
      \right)\sim 0$ for all $x\in \mathbb{D}$,
$\scriptsize\left(
        \begin{array}{cc}
          x^{\sigma_{11}} & x^{\sigma_{12}} \\
        1^{\sigma_{21}}& 1 \\
        \end{array}
      \right)\sim 0$ and $\scriptsize\left(
        \begin{array}{cc}
        1 & x^{\sigma_{12}} \\
        1^{\sigma_{21}}& x^{\sigma_{22}} \\
        \end{array}
      \right)\sim 0$ for all $x\in \mathbb{D}$. It follows that $x^{\sigma_{11}}=x^{\sigma_{12}}1^{\sigma_{21}}$ and $x^{\sigma_{22}}=1^{\sigma_{21}}x^{\sigma_{12}}$
for all $x\in \mathbb{D}$. Therefore, $1=1^{\sigma_{21}}1^{\sigma_{12}}=1^{\sigma_{12}}1^{\sigma_{21}}$ and hence $1^{\sigma_{12}}=(1^{\sigma_{21}})^{-1}$.
Replacing the map $\phi$ by the map
$$X\longmapsto {\rm diag}(1, (1^{\sigma_{21}})^{-1}, 1, \ldots, 1)\phi(X){\rm diag}(1, (1^{\sigma_{12}})^{-1}, 1, \ldots, 1),$$
we have that $1^{\sigma_{ij}}=1$ for $i,j=1,2$. Then $x^{\sigma_{11}}=x^{\sigma_{22}}=x^{\sigma_{12}}$
for all $x\in \mathbb{D}$. Write $\tau=\sigma_{11}$.
Since $\scriptsize\left(
        \begin{array}{cc}
          x & x \\
        x& x \\
        \end{array}
      \right)\sim 0$  for all $x\in \mathbb{D}^*$, $\scriptsize\left(
        \begin{array}{cc}
          x^\tau & x^\tau \\
        x^{\sigma_{21}} & x^\tau \\
        \end{array}
      \right)\sim 0$  for all $x\in \mathbb{D}^*$, which implies that $\sigma_{21}=\tau$. Since
$\scriptsize\left(
        \begin{array}{cc}
          xy & x \\
        y& 1 \\
        \end{array}
      \right)\sim 0$  for all $x, y\in \mathbb{D}$, $\scriptsize\left(
        \begin{array}{cc}
          (xy)^\tau & x^\tau \\
        y^\tau & 1 \\
        \end{array}
      \right)\sim 0$  for all $x, y\in \mathbb{D}$.
Consequently $(xy)^\tau=x^\tau y^\tau$ for all $x, y\in \mathbb{D}$. Thus  $\tau$ is a nonzero ring homomorphism,
and hence (\ref{additiuytnv86h}) holds.
$\qed$

\begin{cor}\label{addi-non-degeneratecor1}    If $\phi:  {\mathbb{D}}^{2\times 2}\rightarrow  {\mathbb{D}'}^{m'\times n'}$ $(m',n'\geq 2)$
is a non-degenerate additive graph homomorphism, then $\phi$ is a distance preserving map.
\end{cor}

\begin{thm}\label{additive-non-degenerate03}   Suppose $\mathbb{D}, \mathbb{D}'$ are division rings and $m,n,m',n'\geq 2$ are integers. Let
$\phi:  {\mathbb{D}}^{m\times n}\rightarrow  {\mathbb{D}'}^{m'\times n'}$
be a non-degenerate additive graph homomorphism. Then either there exist  matrices  $P\in {\mathbb{D}'}^{m'\times m}$ and $Q\in {\mathbb{D}'}^{n\times n'}$ with their ranks $\geq 2$,
and a nonzero ring homomorphism $\tau: \mathbb{D} \rightarrow\mathbb{D}'$,  such that
\begin{equation}\label{additive665cnc00}
 \phi(X)=P X^\tau Q, \ \, X\in {\mathbb{D}}^{m\times n};
\end{equation}
or there exist matrices  $P\in {\mathbb{D}'}^{m'\times n}$ and $Q\in {\mathbb{D}'}^{m\times n'}$ with their ranks $\geq 2$,
and a nonzero ring anti-homomorphism $\sigma: \mathbb{D}\rightarrow\mathbb{D}'$,  such that
\begin{equation}\label{additive-bb6600}
\phi(X)=P \,{^tX^\sigma}Q, \ \,  X\in {\mathbb{D}}^{m\times n}.
\end{equation}
\end{thm}
\proof
By Lemma \ref{s4245549bcnfte}, we may assume with no loss of generality  either
$\phi(\mathcal{M}_i)\subseteq \mathcal{M}_i'$ with $\phi(\mathcal{N}_i)\subseteq \mathcal{N}_i'$ $(i=1,2)$,
or $\phi(\mathcal{M}_i)\subseteq \mathcal{N}_i'$  with  $\phi(\mathcal{N}_i)\subseteq \mathcal{M}_i'$ $(i=1,2)$.
If the first case [resp. second case] happens, we can prove that $\phi$ is of the form (\ref{additive665cnc00}) [resp. (\ref{additive-bb6600})].
We prove this theorem only for  the first case; the second case is similar. From now on we assume that
\begin{equation}\label{additive-utyibvc5}
\mbox{$\phi(\mathcal{M}_i)\subseteq \mathcal{M}_i'$,  \  \ \ $\phi(\mathcal{N}_i)\subseteq \mathcal{N}_i'$, \ $i=1,2$.}
\end{equation}

By Theorem \ref{additive-non-degenerate1} and its proof, we may assume with no loss of generality that
\begin{equation}\label{addit756hfyh09vxv}
\phi\left(\begin{array}{cc}
 X & 0 \\
  0 & 0 \\
   \end{array}
  \right)=\left(\begin{array}{cc}
 X^\tau  & 0 \\
  0 & 0 \\
   \end{array}
  \right), \ \, X\in {\mathbb{D}}^{2\times 2},
\end{equation}
where $\tau: \mathbb{D} \rightarrow\mathbb{D}'$ is a nonzero ring homomorphism.

Write $E_{ij}=E_{ij}^{m\times n}$ and $E_{ij}'=E_{ij}^{m'\times n'}$. For  $i\in\{1,2\}$ and $j\in\{3, \ldots, n\}$, we let
$$\psi_{ij}\left(
          \begin{array}{cc}
            x_{1i} & x_{1j} \\
            x_{2i} & x_{2j} \\
          \end{array}
        \right)
=\phi\left(x_{1i}E_{1i}+ x_{1j}E_{1j}+x_{2i}E_{2i}+x_{2j}E_{2j}\right), \ \ x_{1i}, x_{1j}, x_{2i}, x_{2j}\in \mathbb{D}.$$
Then $\psi_{ij}:  {\mathbb{D}}^{2\times 2}\rightarrow  {\mathbb{D}'}^{m'\times n'}$ is an additive graph homomorphism.
For  $j\in\{3, \ldots, n\}$, there exists some $k\in \{1,2\}$ such that $\psi_{kj}$ is non-degenerate.
 Otherwise, both $\psi_{1j}$ and $\psi_{2j}$ are degenerate, and  Theorem \ref{additiverank} implies that $\psi_{1j}({\mathbb{D}}^{2\times 2})$
 and $\psi_{2j}({\mathbb{D}}^{2\times 2})$ are two adjacent sets.
   By Corollary \ref{Rectangular-1-13}, (\ref{addit756hfyh09vxv})  and  Lemma \ref{Rectangular-PID2-11}, it is easy to see that
 $\psi_{1j}({\mathbb{D}}^{2\times 2})\subseteq \mathcal{N}_1'$ and  $\psi_{2j}({\mathbb{D}}^{2\times 2})\subseteq \mathcal{N}_2'$. It follows that
 $\phi\left(x_{1j}E_{1j}+x_{2j}E_{2j}\right)\in \mathcal{N}_1'\cap \mathcal{N}_2'=\{0\}$ for all $x_{1j}, x_{2j}\in \mathbb{D}$, a contradiction.
Thus $\psi_{kj}$ is non-degenerate for some $k\in \{1,2\}$.

 By Theorem \ref{additive-non-degenerate1} and (\ref{addit756hfyh09vxv}), we can assume that
 \begin{equation}\label{gd53523fvxdmnv}
\psi_{kj}\left(
          \begin{array}{cc}
            x_{1k} & x_{1j} \\
            x_{2k} & x_{2j} \\
          \end{array}
        \right)
=\phi\left(x_{1k}E_{1k}+ x_{1j}E_{1j}+x_{2k}E_{2k}+x_{2j}E_{2j}\right)=P_k\left(
                \begin{array}{cc}
                  \left(
        \begin{array}{cc}
          x_{1k}^\mu & x_{1j}^\mu  \\
          x_{2k}^\mu  & x_{2j}^\mu \\
        \end{array}
      \right) & 0 \\
                  0 & 0 \\
                \end{array}
              \right)Q_j,
\end{equation}
for all  $x_{1k}, x_{1j}, x_{2k}, x_{2j}\in \mathbb{D}$, where  $P_k\in GL_{m'}(\mathbb{D}')$ is a diagonal matrix, $Q_j\in GL_{n'}(\mathbb{D}')$ and
$\mu: \mathbb{D} \rightarrow\mathbb{D}'$ is a nonzero ring homomorphism.
 Let $q_{ji}=(q_{ji}^{(1)}, \ldots, q_{ji}^{(n')})$ be the $i$-th row of  $Q_j$ and let $P_k={\rm diag}(p_{k1}, \ldots, p_{km'})$. Then (\ref{gd53523fvxdmnv}) and (\ref{addit756hfyh09vxv}) imply that
 \begin{equation}\label{hfg43xad13p}
\phi(xE_{1k})=\left(
                   \begin{array}{c}
                     p_{k1}x^\mu q_{j1} \\
                     0\\
                     0_{m'-2, n'} \\
                   \end{array}
                 \right)=x^\tau E_{1k}', \ \ \ \phi(xE_{2k})=\left(
                   \begin{array}{c}
                   0\\
                     p_{k2}x^\mu q_{j1} \\
                     0_{m'-2, n'} \\
                   \end{array}
                 \right)=x^\tau E_{2k}', \ \ x\in \mathbb{D};
 \end{equation}
 \begin{equation}\label{97dghnmqweq53}
\phi(xE_{1j})=\left(
                   \begin{array}{c}
                     p_{k1}x^\mu q_{j2} \\
                     0\\
                      0_{m'-2, n'} \\
                   \end{array}
                 \right), \ \ \  \phi(xE_{2j})=\left(
                   \begin{array}{c}
                   0\\
                   p_{k2}x^\mu q_{j2} \\
                   0_{m'-2, n'} \\
                   \end{array}
                 \right), \ \ x\in \mathbb{D}.
 \end{equation}
By (\ref{hfg43xad13p}), we have $x^\tau= p_{k1}x^\mu q_{j1}^{(k)}$ and $x^\tau= p_{k2}x^\mu q_{j1}^{(k)}$ for all $x\in \mathbb{D}$. Thus $p_{k1}=p_{k2}=(q_{j1}^{(k)})^{-1}$.

Let $e_i$  be the $i$-th row of $I_{m'}$, and let  $e_i'$  be the $i$-th row of $I_{n'}$.
Write $q_{j2}'=(q_{j1}^{(k)})^{-1}q_{j2}$.
Then (\ref{97dghnmqweq53}) can be written as
 \begin{equation}\label{te3hfk1v30}
\phi(xE_{1j})= \,^te_1  x^\tau q_{j2}', \ \ \  \phi(xE_{2j})= \,^te_2  x^\tau q_{j2}', \ \ x\in \mathbb{D}, \ j=3,\ldots, n.
 \end{equation}

 For  $s\in\{1,2\}$ and $i\in\{3, \ldots, n\}$, we define
$$\psi_{si}'\left(
          \begin{array}{cc}
            x_{s1} & x_{s2} \\
            x_{i1} & x_{i2} \\
          \end{array}
        \right)
=\phi\left(x_{s1}E_{s1}+ x_{s2}E_{s2}+x_{i1}E_{i1}+x_{i2}E_{i2}\right), \ \ x_{s1}, x_{s2}, x_{i1}, x_{i2}\in \mathbb{D}.$$
Then $\psi_{si}':  {\mathbb{D}}^{2\times 2}\rightarrow  {\mathbb{D}'}^{m'\times n'}$ is an additive graph homomorphism.
Similarly, there exists some $s\in \{1,2\}$ such that $\psi_{si}'$ is non-degenerate. Similar to the proof of (\ref{te3hfk1v30}),
there is $0\neq p_{i2}'\in {\mathbb{D}'}^{m'\times 1}$ such that
\begin{equation}\label{ry6iadmnv20b}
\phi(xE_{i1})= p_{i2}' x^\tau e_1', \ \ \  \phi(xE_{i2})= p_{i2}' x^\tau e_2', \ \ x\in \mathbb{D}, \ i=3,\ldots, n.
 \end{equation}

Recall that $\psi_{kj}$ is non-degenerate for some $k\in \{1,2\}$.
Put $i\in\{3, \ldots, m\}$,  $j\in\{3, \ldots, n\}$  and $r\in \{1,2\}$. We let
$$\theta_{ij}^{(r)}\left(
          \begin{array}{cc}
            x_{rk} & x_{rj} \\
            x_{ik} & x_{ij} \\
          \end{array}
        \right)
=\phi\left(x_{rk}E_{rk}+ x_{rj}E_{rj}+x_{ik}E_{ik}+x_{ij}E_{ij}\right), \ \ x_{rk}, x_{rj}, x_{ik}, x_{ij}\in \mathbb{D}.$$
Then $\theta_{ij}^{(r)}:  {\mathbb{D}}^{2\times 2}\rightarrow  {\mathbb{D}'}^{m'\times n'}$ is an additive graph homomorphism.
There exists some $r\in \{1,2\}$ such that $\theta_{ij}^{(r)}$ is non-degenerate. Otherwise, both $\theta_{ij}^{(1)}$ and $\theta_{ij}^{(2)}$ are degenerate.
By Theorem \ref{additiverank}, we have $\theta_{ij}^{(1)}({\mathbb{D}}^{2\times 2})\subseteq \mathcal{M}$, where $\mathcal{M}$ is a maximal set in ${\mathbb{D}'}^{m'\times n'}$ containing $0$.
From (\ref{addit756hfyh09vxv}) we get $\scriptsize\theta_{ij}^{(1)}\left(
          \begin{array}{cc}
            x_{1k} & 0 \\
            0 & 0 \\
          \end{array}
        \right)=x_{1k}^\tau E_{1k}'$ for all $x_{1k}\in \mathbb{D}$. Thus, Corollaries \ref{Rectangular-PID2-7} and \ref{Rectangular-1-13}  imply that
 $\mathcal{M}=\mathcal{M}_1'$ or $\mathcal{M}=\mathcal{N}_k'$.
Suppose  that $\mathcal{M}=\mathcal{N}_k'$. Then (\ref{additive-utyibvc5}) implies that
$\scriptsize\theta_{ij}^{(1)}\left(
          \begin{array}{cc}
            x_{1k} & x_{1j} \\
            0 & 0 \\
          \end{array}
        \right)\in \mathcal{N}_k'\cap \mathcal{M}_1'=\mathbb{D}'E_{1k}'$, and hence
 $\phi(E_{1j})=a_1E_{1k}'$ where $a_1\in {\mathbb{D}'}^*$. By (\ref{hfg43xad13p}) and (\ref{97dghnmqweq53}), it is easy to see that $q_{j1}$ and $q_{j2}$ are left
linearly dependent, a contradiction to $Q_j$ being invertible.
Therefore, we must have  $\mathcal{M}=\mathcal{M}_1'$ and hence $\theta_{ij}^{(1)}({\mathbb{D}}^{2\times 2})\subseteq \mathcal{M}_1'$.  Similarly,
 $\theta_{ij}^{(2)}({\mathbb{D}}^{2\times 2})\subseteq \mathcal{M}_2'$. It follows that
 $\phi\left(x_{ik}E_{ik}+x_{ij}E_{ij}\right)\in \mathcal{M}_1'\cap \mathcal{M}_2'=\{0\}$ for all $x_{ik}, x_{ij}\in \mathbb{D}$, a contradiction.
Thus, $\theta_{ij}^{(r)}$ is non-degenerate for some $r\in \{1,2\}$.

 By  Theorem \ref{additive-non-degenerate1} and (\ref{additive-utyibvc5}), we have
 \begin{equation}\label{ns13jgkl0szxz}
\theta_{ij}^{(r)}\left(
          \begin{array}{cc}
            x_{rk} & x_{rj} \\
            x_{ik} & x_{ij} \\
          \end{array}
        \right)
=\phi\left(x_{rk}E_{rk}+ x_{rj}E_{rj}+x_{ik}E_{ik}+x_{ij}E_{ij}\right)=S_i\left(
                \begin{array}{cc}
                  \left(
        \begin{array}{cc}
          x_{rk}^\delta & x_{rj}^\delta  \\
          x_{ik}^\delta  & x_{ij}^\delta \\
        \end{array}
      \right) & 0 \\
                  0 & 0 \\
                \end{array}
              \right)T_j,
\end{equation}
for all  $x_{rk}, x_{rj}, x_{ik}, x_{ij}\in \mathbb{D}$, where  $S_i\in GL_{m'}(\mathbb{D}')$, $T_j\in GL_{n'}(\mathbb{D}')$ and
$\delta: \mathbb{D} \rightarrow\mathbb{D}'$ is a nonzero ring homomorphism.
 Let $t_{ji}=(t_{ji}^{(1)}, \ldots, t_{ji}^{(n')})$ be the $i$-th row of  $T_j$,  and let $S_i=(s_{i1}, \ldots, s_{im'})$ where $s_{ij}$
be the $j$-th column of  $S_i$. Then (\ref{ns13jgkl0szxz}) and (\ref{addit756hfyh09vxv}) imply that
\begin{equation}\label{c357hfisag}
\phi(xE_{rk})=s_{i1}x^\delta t_{j1}=x^\tau E_{rk}'=\,^te_rx^{\tau} e_k',  \ \ x\in \mathbb{D};
 \end{equation}
$$
\phi(xE_{rj})=s_{i1}x^\delta t_{j2},  \ \ x\in \mathbb{D};
$$
$$
\phi(xE_{ik})=s_{i2}x^\delta t_{j1}, \ \ x\in \mathbb{D};
$$
$$
 \phi(xE_{ij})=s_{i2}x^\delta t_{j2}, \ \ x\in \mathbb{D}.
 $$

By (\ref{c357hfisag}), there are $a,b\in {\mathbb{D}'}^*$ such that $s_{i1}=\,^te_ra$ and $t_{j1}=be_k'$. Moreover, $x^\delta=a^{-1}x^\tau b^{-1}$ for all $x\in \mathbb{D}$.
 Since $1^\delta=1^\tau=1$, we get that  $b^{-1}=a$. Hence
$$
\mbox{ $s_{i1}=\,^te_ra$, \ $at_{j1}=e_k'$ \ and \ $x^\delta =a^{-1}x^\tau a$, \ $x\in \mathbb{D}$.}
$$
Therefore, for above  two fixed $k, r\in \{1,2\}$, we obtain that
\begin{equation}\label{bvc242fsfs529}
\phi(xE_{rj})=\,^te_r x^\tau a t_{j2},  \ \ x\in \mathbb{D}, j=3, \ldots, n;
 \end{equation}
 \begin{equation}\label{bcadar3gd075}
\phi(xE_{ik})=s_{i2}a^{-1}x^\tau at_{j1}=s_{i2}a^{-1}x^\tau e_k', \ \ x\in \mathbb{D},  \ i=3, \ldots, m;
 \end{equation}
\begin{equation}\label{o935fsfs3edf7}
 \phi(xE_{ij})=s_{i2} a^{-1}x^\tau a t_{j2}, \ \ x\in \mathbb{D}, \ i=3, \ldots, m, j=3, \ldots, n.
 \end{equation}
By (\ref{bvc242fsfs529}) and (\ref{te3hfk1v30}), we get $q_{j2}'=at_{j2}$. Thus (\ref{te3hfk1v30}) can be written as
 \begin{equation}\label{647dweqhazx64h}
\phi(xE_{1j})= \,^te_1  x^\tau at_{j2}, \ \ \  \phi(xE_{2j})= \,^te_2  x^\tau at_{j2}, \ \ x\in \mathbb{D}, \ j=3,\ldots, n.
 \end{equation}
By (\ref{ry6iadmnv20b}) and (\ref{bcadar3gd075}), we have $p_{i2}'=s_{i2}a^{-1}$.  Hence (\ref{ry6iadmnv20b}) can be written as
\begin{equation}\label{5485utcvsmi6t}
\phi(xE_{i1})= s_{i2}a^{-1} x^\tau e_1', \ \ \  \phi(xE_{i2})= s_{i2}a^{-1} x^\tau e_2', \ \ x\in \mathbb{D}, \ i=3,\ldots, n.
 \end{equation}

Let $\scriptsize P=\left( ^te_1, \,^te_2, s_{32}a^{-1},  s_{42}a^{-1}, \ldots,  s_{m2}a^{-1} \right)=\left(
              \begin{array}{cc}
                I_2 & \ast \\
                 0& \ast \\
              \end{array}
            \right)\in {\mathbb{D}'}^{m'\times m}$, and let
$$Q=\left(
      \begin{array}{c}
        e_1' \\
        e_2' \\
        at_{32} \\
        at_{42} \\
        \vdots \\
        at_{n2} \\
      \end{array}
    \right)=\left(
              \begin{array}{cc}
                I_2 & 0 \\
                \ast & \ast \\
              \end{array}
            \right)\in {\mathbb{D}'}^{n\times n'}.$$
Then ${\rm rank}(P)\geq 2$ and ${\rm rank}(Q)\geq 2$. By (\ref{addit756hfyh09vxv}) with (\ref{o935fsfs3edf7})-(\ref{5485utcvsmi6t}), we obtain that
$$\phi(xE_{ij})=Px^\tau E_{ij}'Q, \ \ x\in \mathbb{D}, \ i=1,\ldots, m, j=1,\ldots, n.$$
Since $\phi$ is additive, we get (\ref{additive665cnc00}).
$\qed$

By Theorem \ref{additive-non-degenerate03} and Theorem \ref{additiverank}, it is clear that Theorem \ref{additivehomomorphism0a} holds.

\section{Graph homomorphisms on $n\times n$ matrices}

\ \ \ \ \ \
In this section, we will expand \cite[Theorem 4.1]{optimal} (which is due to \v{S}emrl) to the cases of two division rings and $n=2$.

\begin{thm}\label{non-degenerate-b} \  Let $\mathbb{D}, \mathbb{D}'$ be division rings with $|\mathbb{D}|\geq 4$,
and let $m',n',n$ be integers with $m', n'\geq n\geq 2$. Suppose that
$\varphi:  {\mathbb{D}}^{n\times n}\rightarrow  {\mathbb{D}'}^{m'\times n'}$
is a graph homomorphism with $\varphi(0)=0$,  and there exists $A_0\in {\mathbb{D}}^{n\times n}$ such that ${\rm rank}(\varphi(A_0))=n$.
Then either there exist  matrices  $P\in GL_{m'}(\mathbb{D}')$ and $Q\in GL_{n'}(\mathbb{D}')$,
a nonzero ring homomorphism $\tau: \mathbb{D} \rightarrow\mathbb{D}'$, and a matrix $L\in {\mathbb{D}'}^{n\times n}$ with the property that
$I_n+X^\tau L\in GL_n(\mathbb{D}')$ for every $X\in {\mathbb{D}}^{n\times n}$, such that
\begin{equation}\label{3654665cnc00}
 \varphi(X)=P\left(\begin{array}{cc}
 (I_n+X^\tau L)^{-1}X^\tau  & 0 \\
  0 & 0 \\
   \end{array}
  \right)Q, \ \, X\in {\mathbb{D}}^{n\times n};
\end{equation}
or there exist matrices  $P\in GL_{m'}(\mathbb{D}')$ and $Q\in GL_{n'}(\mathbb{D}')$,
a nonzero ring anti-homomorphism $\sigma: \mathbb{D}\rightarrow\mathbb{D}'$, and a matrix $L\in {\mathbb{D}'}^{n\times n}$
with the property that $I_n+L\,^tX^\sigma\in GL_n(\mathbb{D}')$ for every $X\in {\mathbb{D}}^{n\times n}$, such that
\begin{equation}\label{CX3mmbb6600}
\varphi(X)=P\left(\begin{array}{cc}
{^tX^\sigma}(I_n+L\,^tX^\sigma)^{-1} & 0 \\
0 & 0 \\
\end{array}
\right)Q, \ \,  X\in {\mathbb{D}}^{n\times n};
\end{equation}
or $\varphi\left(\mathbb{D}_{\leq 1}^{n\times n}\right)$  is  an adjacent set.
In particular, if $\tau$ $[$resp. $\sigma$$]$ is surjective, then  $L=0$ and  $\tau$ is a ring isomorphism
$[$resp. $\sigma$ is a ring anti-isomorphism$]$.
\end{thm}

 By Theorem \ref{non-degenerate-b}, it is easy to prove the following corollary.

\begin{cor}\label{non-degenerate-cc02} \  Let $\mathbb{D}, \mathbb{D}'$ be division rings with $|\mathbb{D}|\geq 4$,
and let $m',n',n$ be integers with $m', n'\geq n\geq 2$. Suppose that
$\varphi:  {\mathbb{D}}^{n\times n}\rightarrow  {\mathbb{D}'}^{m'\times n'}$
is a graph homomorphism,  and there exist $A_0, B_0\in {\mathbb{D}}^{n\times n}$ such that ${\rm rank}(\varphi(B_0)-\varphi(A_0))=n$.
Then either there exist  matrices  $P\in GL_{m'}(\mathbb{D}')$ and $Q\in GL_{n'}(\mathbb{D}')$,
a nonzero ring homomorphism $\tau: \mathbb{D} \rightarrow\mathbb{D}'$, and a matrix $L\in {\mathbb{D}'}^{n\times n}$ with the property that
$I_n+X^\tau L\in GL_n(\mathbb{D}')$ for every $X\in {\mathbb{D}}^{n\times n}$, such that
\begin{equation}\label{3654665cnc0c2}
 \varphi(X)=P\left(\begin{array}{cc}
 (I_n+X^\tau L)^{-1}(X^\tau-A_0^\tau)  & 0 \\
  0 & 0 \\
   \end{array}
  \right)Q+\varphi(A_0), \ \, X\in{\mathbb{D}}^{n\times n};
\end{equation}
or there exist matrices  $P\in GL_{m'}(\mathbb{D}')$ and $Q\in GL_{n'}(\mathbb{D}')$,
a nonzero ring anti-homomorphism $\sigma: \mathbb{D}\rightarrow\mathbb{D}'$, and a matrix $L\in {\mathbb{D}'}^{n\times n}$
with the property that $I_n+L\,^tX^\sigma\in GL_n(\mathbb{D}')$ for every $X\in {\mathbb{D}}^{n\times n}$, such that
\begin{equation}\label{CX3mmbb660c2}
\varphi(X)=P\left(\begin{array}{cc}
({^tX^\sigma}-{^tA_0^\sigma})(I_n+L\,^tX^\sigma)^{-1} & 0 \\
0 & 0 \\
\end{array}
\right)Q+\varphi(A_0), \ \,  X\in {\mathbb{D}}^{n\times n};
\end{equation}
or both $\varphi\left(\mathbb{B}_{A_0}\right)$ and $\varphi\left(\mathbb{B}_{B_0}\right)$ are  adjacent sets.
In particular, if $\tau$ $[$resp. $\sigma$$]$ is surjective, then  $L=0$ and  $\tau$ is a ring isomorphism
$[$resp. $\sigma$ is a ring anti-isomorphism$]$.
\end{cor}

Note that maps (\ref{3654665cnc00})-(\ref{CX3mmbb660c2}) are distance preserving maps (cf. \cite[Examples 2.5 and 2.6]{Huangli-I}. We have:

\begin{cor}\label{stronglydegenerate2} \ Let $\mathbb{D}, \mathbb{D}'$ be division rings with $|\mathbb{D}|\geq 4$,
and let $m',n',n$ be integers with $m', n'\geq n\geq 2$.  Suppose that $\varphi:  {\mathbb{D}}^{n\times n}\rightarrow  {\mathbb{D}'}^{m'\times n'}$
is a graph homomorphism,  and there exist $A_0, B_0\in {\mathbb{D}}^{n\times n}$ such that ${\rm rank}(\varphi(B_0)-\varphi(A_0))=n$.
Then either $\varphi$ is a distance preserving map, or  both $\varphi\left(\mathbb{B}_{A_0}\right)$ and
$\varphi\left(\mathbb{B}_{B_0}\right)$ are  adjacent sets.
\end{cor}

In order to prove  Theorem \ref{non-degenerate-b} we need the following  knowledge and lemmas.

It is well-known that $\mathbb{D}^{m\times n}$ ($m,n\geq 2$) is a partially ordered set (poset) with partial order defined by
$A\leq B$ if ${\rm rank}(B-A)={\rm rank}(B)-{\rm rank}(A)$. This partial order is also called {\em minus order} (or {\em minus  partial order}) of $D^{m\times n}$
(cf. \cite{Mitra}). In particular, we write $A< B$ $\Leftrightarrow$ $A\leq B$ with $A\neq B$.

For any $X\in {\mathbb{D}}^{m\times n}$, a g-inverse of $X$ will be denoted by $X^-$ and is understood as a  matrix (over $\mathbb{D}$)
for which $XX^-X=X$.
Two matrices $A,B\in {\mathbb{D}}^{m\times n}$ are said to be {\em equivalent}, denoted by $A\cong B$, if $B$ may be obtained from
$A$ by a finite sequence of elementary row and column operations.

\begin{lem}\label{minuspar1}{\rm (cf. \cite[Section 3.3]{Mitra})} \ Let $\mathbb{D}$ be a division ring and let $A,B\in \mathbb{D}^{m\times n}$.
Then the following results are equivalent:
\begin{itemize}
\item[\rm(a)] \  $A\leq B$;
\item[\rm(b)] \  there are g-inverses $G_1$ and $G_2$ of $A$ such that $AG_1=BG_1$ and $G_2A=G_2B$;
\item[\rm(c)] \ there are $P\in GL_m(\mathbb{D})$ and $Q\in GL_n(\mathbb{D})$ such that
$B=P{\rm diag}(I_{r+s}, 0)Q$ and
$A=P{\rm diag}(I_{r}, 0)Q$;
\item[\rm(d)] \ $PAQ\leq PBQ$, where $P\in GL_m(\mathbb{D})$ and $Q\in GL_n(\mathbb{D})$.
\end{itemize}
\end{lem}

\begin{lem}\label{minuspar1-3} {\rm (see \cite[Theorem 3.6.5]{Mitra})} \
Let $A, B\in \mathbb{D}^{n\times n}$ with $B^2=B$. Then  $A\leq B$ if and only if $A=A^2=AB=BA$.
\end{lem}

\begin{lem}\label{midsad23par}  \ Let  $A, B, C\in {\mathbb{D}}^{m\times n}$  and
$d(A,B)=d(B,C)-d(A,C)$. Suppose that $\varphi:  {\mathbb{D}}^{m\times n} \rightarrow {\mathbb{D}'}^{m'\times n'}$ is a graph homomorphism and
$d(B,C)=d(\varphi(B), \varphi(C))$. Then $d(\varphi(A),\varphi(B))=d(\varphi(B),\varphi(C))-d(\varphi(A),\varphi(C))$.
\end{lem}
\proof  We have
$d(\varphi(A), \varphi(B))\leq d(A,B)=d(B,C)-d(A,C)=d(\varphi(B),\varphi(C))-d(A,C)\leq d(\varphi(B),\varphi(C))-d(\varphi(A),\varphi(C)).$
 Hence
 $$d(\varphi(A), \varphi(B))\leq d(\varphi(B),\varphi(C))-d(\varphi(A),\varphi(C)).$$
On the other hand, it is clear that $d(\varphi(A), \varphi(B))\geq d(\varphi(B),\varphi(C))-d(\varphi(A),\varphi(C))$.
Thus $d(\varphi(A), \varphi(B))=d(\varphi(B),\varphi(C))-d(\varphi(A), \varphi(C))$.
$\qed$

\begin{cor}\label{mYR25HLLsad23} \ Let
$\varphi:  {\mathbb{D}}^{m\times n} \rightarrow {\mathbb{D}'}^{m'\times n'}$ be a graph homomorphism with $\varphi(0)=0$, and let
 $A,B\in{\mathbb{D}}^{m\times n}$.
  If $A\leq B$ and ${\rm rank}(\varphi(B))={\rm rank}(B)$, then $\varphi(A)\leq \varphi(B)$ and  ${\rm rank}(\varphi(A))={\rm rank}(A)$.
\end{cor}
\proof
Assume that $A\leq B$ and ${\rm rank}(\varphi(B))={\rm rank}(B)$. Then $d(A,B)=d(B,0)-d(A,0)$ and $d(B,0)=d(\varphi(B), \varphi(0))$.
By Lemma \ref{midsad23par},  $d(\varphi(A),\varphi(B))=d(\varphi(B),\varphi(0))-d(\varphi(A),\varphi(0))$,  hence $\varphi(A)\leq \varphi(B)$.
It follows that ${\rm rank}(\varphi(A))=d(\varphi(A),\varphi(0))=d(\varphi(B),\varphi(0))-d(\varphi(A),\varphi(B))\geq {\rm rank}(B)-d(A,B)={\rm rank}(A)$.
Since ${\rm rank}(\varphi(A))\leq {\rm rank}(A)$, we obtain  ${\rm rank}(\varphi(A))={\rm rank}(A)$.
$\qed$

\begin{lem}\label{Recastatre43} \ Let $\mathbb{D}, \mathbb{D}'$ be division rings,  and let $n,m', n'$ be integers with $m',n'\geq n\geq 2$.
Suppose that $\varphi:  {\mathbb{D}}^{n\times n}\rightarrow  {\mathbb{D}'}^{m'\times n'}$
is a graph homomorphism with $\varphi(0)=0$,  and there exists $A_0\in {\mathbb{D}}^{n\times n}$ such that $d(\varphi(A_0), 0)=n$.
 Assume that ${\cal M}$ is a maximal set  in ${\mathbb{D}}^{n\times n}$ containing $0$ and  $\varphi(\mathcal{M})\subseteq {\cal M}'$
where ${\cal M}'$  is a maximal set in ${\mathbb{D}'}^{m'\times n'}$.
 Then $\varphi({\cal M})$ is not contained in any line in $AG({\cal M}')$, and
 ${\cal M}'$ is the unique  maximal set  containing $\varphi(\mathcal{M})$ in ${\mathbb{D}'}^{m'\times n'}$.
 \end{lem}
\proof
Without loss of generality, we assume that  ${\cal M}$ is a maximal set of type two,  and  ${\cal M}'$ is a maximal set of type one.
 By Lemma \ref{Rectangular-PID2-4}, there are invertible matrices $P_1, P_2$ such that  $\mathcal{M}={\cal N}_1P_1$
and $\mathcal{M}'= P_2{\cal M}'_1$. Replacing  $\varphi$ by the map $X\mapsto P_2^{-1}\varphi(XP_1)$,
we have ${\cal M}={\cal N}_1$, ${\cal M}'={\cal M}_1'$ and  $\varphi({\cal N}_1)\subseteq {\cal M}'_{1}$.

 Suppose $\varphi({\cal N}_1)$ is contained in a line $\ell$ in $AG({\cal M}'_1)$. We  show  a contradiction as follows.
 By the parametric equation  of a line, there exists an invertible matrix $Q$ such that $\ell=(\mathbb{D}'E_{11}+B)Q$, where $E_{11}=E_{11}^{m'\times n'}$
and $B\in{\cal M}'_{1}$. Since  $0\in \ell$, we can assume that  $B=0$ and $\ell=\mathbb{D}'E_{11}Q$.
 Replacing $\varphi$ by the map $X\mapsto \varphi(X)Q^{-1}$,
 we have $\varphi(0)=0$ and $\varphi({\cal N}_1)\subseteq \ell=\mathbb{D}'E_{11}$.
By the conditions, there is an $A_0\in {\mathbb{D}}^{n\times n}$ such that  $d(\varphi(A_0), 0)=n=d(A_0,0)$.
Let $A_0=(B_1, B_2)$ where $B_1\in \mathbb{D}_2^{n\times 2}$ and $B_2\in \mathbb{D}_{n-2}^{n\times (n-2)}$.
When $n=2$, $B_2$ is  absent.
Put $A_1=(B_1,  0)\in\mathbb{D}_2^{n\times n}$.
Then  $A_1\leq A_0$.  Using Corollary \ref{mYR25HLLsad23}, we get  ${\rm rank}(\varphi(A_1))=2$.
Clearly, there are two distinct points $Y_1, Y_2\in {\cal N}_1$ such that $A_1\sim Y_i$, $i=1,2$. Since $\varphi(Y_i)\subseteq \ell$,
$\varphi(A_1)$ are adjacent to  two distinct points in $\mathbb{D}'E_{11}$. By Lemma \ref{Matrix-PID5-4bb},
$\varphi(A_1)\in {\cal M}_1'$ or  $\varphi(A_1)\in {\cal N}_1'$, and
hence ${\rm rank}(\varphi(A_1))\leq 1$, a contradiction.

Thus $\varphi({\cal N}_1)$ is not contained in any line of $AG({\cal M}'_1)$. By Corollary \ref{Rectangular-1-13} and Lemma \ref{Rectangular-PID2-11},
it is easy to see that  ${\cal M}'$ is the unique  maximal set containing $\varphi(\mathcal{M})$.
$\qed$

\begin{lem}\label{strong-degenerate-1}{\rm (cf. \cite{optimal})} \  Let $\mathbb{D}, \mathbb{D}'$ be division rings,
and let $m',n',n$ be integers with $m', n'\geq n\geq 2$. Suppose that
$\varphi:  {\mathbb{D}}^{n\times n}\rightarrow  {\mathbb{D}'}^{m'\times n'}$
is a graph homomorphism with $\varphi(0)=0$ and  ${\rm rank}(\varphi(I_n))=n$.
Then there exist invertible matrices $P\in GL_{m'}(\mathbb{D}')$ and $Q\in GL_{n'}(\mathbb{D}')$ such that
\begin{equation}\label{vcx2354fgfd96s}
\varphi\left({\rm diag}(I_r, 0)\right)=P{\rm diag}\left(I_r, 0_{m'-r,n'-r}\right)Q, \ \ r=1,\ldots,n.
 \end{equation}
 \end{lem}
\proof
 Since ${\rm diag}(I_{n-1},0)<I_n$ and  ${\rm rank}(\varphi(I_n))=n$, from Corollary \ref{mYR25HLLsad23} we get that
 $$\mbox{${\rm rank}\left(\varphi\left({\rm diag}(I_{n-1},0)\right)\right)=n-1$ \ and \ $\varphi\left({\rm diag}(I_{n-1},0)\right)< \varphi(I_n)$.}$$
Using Lemma \ref{minuspar1}(c), there exist $Q_1\in GL_{m'}(\mathbb{D}')$ and $Q_2\in GL_{n'}(\mathbb{D}')$ such that
$\varphi(I_n)=Q_1{\rm diag}(I_{n}, 0)Q_2$ and $\varphi\left({\rm diag}(I_{n-1},0)\right)=Q_1{\rm diag}(I_{n-1}, 0)Q_2$.
 Replacing $\varphi$ by the map $X\mapsto Q_1^{-1}\varphi(X)Q_2^{-1}$, we have that $\varphi(I_n)={\rm diag}\left(I_{n}, 0\right)$ and
$\varphi\left({\rm diag}(I_{n-1},0)\right)={\rm diag}\left(I_{n-1}, 0_{m'-n+1,n'-n+1}\right)$.

Now, let $n\geq 3$. Since ${\rm diag}(I_{n-2},0)<{\rm diag}(I_{n-1},0)$,
 Corollary \ref{mYR25HLLsad23} implies that ${\rm rank}\left(\varphi\left({\rm diag}(I_{n-2},0\right)\right)=n-2$ and
\begin{equation}\label{vczczwq5369lj}
\varphi\left({\rm diag}(I_{n-2},0)\right)< {\rm diag}(I_{n-1},0_{m'-n+1,n'-n+1}).
\end{equation}
Write $\scriptsize\varphi\left({\rm diag}(I_{n-2},0)\right)=\left(
                                                   \begin{array}{cc}
                                                     A_{11} & A_{12} \\
                                                    A_{21} &   A_{22} \\
                                                   \end{array}
                                                 \right)$ where $A_{11}\in \mathbb{{D}'}^{(n-1)\times (n-1)}$.
By (\ref{vczczwq5369lj}), Lemma \ref{minuspar1-3} and multiplications of matrices,
it is easy to see that $A_{12}=0$, $A_{21}=0$ and $A_{11}$ is an idempotent matrix.
 (Note that when $m'\neq n'$ and  using Lemma \ref{minuspar1-3}, we can get two suitable square matrices by adding some zero elements
on $\varphi\left({\rm diag}(I_{n-2},0)\right)$ and ${\rm diag}(I_{n-1},0_{m'-n+1,n'-n+1})$.)
Since $\varphi\left({\rm diag}(I_{n-2},0)\right)\sim{\rm diag}(I_{n-1},0_{m'-n+1,n'-n+1})$, we have $A_{22}=0$. Hence
$\varphi\left({\rm diag}(I_{n-2},0)\right)={\rm diag}\left(A_{11}, 0_{m'-n+1,n'-n+1}\right)$, where $A_{11}$ is an idempotent matrix of rank $n-2$.
There is $T\in GL_{n-1}(\mathbb{D}')$ such that $A_{11}=T^{-1}{\rm diag}(I_{n-2},0)T$. Let $Q_3={\rm diag}(T^{-1}, I_{m'-n+1})$ and
$Q_4={\rm diag}(T, I_{n'-n+1})$.  Replacing $\varphi$ by the map $X\mapsto Q_3^{-1}\varphi(X)Q_4^{-1}$, we get that
$\varphi(I_n)={\rm diag}(I_{n}, 0_{m'-n,n'-n})$, $\varphi\left({\rm diag}(I_{n-1},0)\right)={\rm diag}\left(I_{n-1}, 0_{m'-n+1,n'-n+1}\right)$ and
$\varphi\left({\rm diag}(I_{n-2},0)\right)={\rm diag}\left(I_{n-2}, 0_{m'-n+2,n'-n+2}\right)$.

Similarly, after by some transformations of the form  $X\mapsto P'\varphi(X)Q'$, we can get
$\varphi\left({\rm diag}(I_{r},0)\right)={\rm diag}\left(I_{r}, 0_{m'-r,n'-r}\right)$, $r=1,\ldots, n$. Thus (\ref{vcx2354fgfd96s}) holds.
$\qed$

 A graph homomorphism $\varphi:  {\mathbb{D}}^{m\times n}\rightarrow  {\mathbb{D}'}^{m'\times n'}$
with $\varphi(0)=0$ is called {\em to satisfy Condition (I)}, if for any two  maximal sets ${\cal M}$ and ${\cal N}$ of different types
 in ${\mathbb{D}}^{m\times n}$ with $0\in{\cal M}\cap {\cal N}$, there are two maximal sets ${\cal M}'$ and ${\cal N}'$ of different types
in ${\mathbb{D}'}^{m'\times n'}$, such that $\varphi({\cal M})\subseteq {\cal M}'$ and $\varphi({\cal N})\subseteq {\cal N}'$.

\begin{lem}\label{slynon-degenerate-2}  Let $\mathbb{D}, \mathbb{D}'$ be division rings,
and let $m',n',n$ be integers with $m', n'\geq n\geq 2$. Suppose that
$\varphi:  {\mathbb{D}}^{n\times n}\rightarrow  {\mathbb{D}'}^{m'\times n'}$
is a graph homomorphism with $\varphi(0)=0$,  and there exists $A_0\in {\mathbb{D}}^{n\times n}$ such that ${\rm rank}(\varphi(A_0))=n$.
Assume further that $\varphi$ satisfies the Condition (I). Then
\begin{equation}\label{d2345fgds53n}
{\rm dim}(\varphi(\mathcal{M}_1))={\rm dim}(\varphi(\mathcal{N}_1))=n.
\end{equation}
\end{lem}
\proof By $\varphi(0)=0$ and (\ref{uy87mbmm}), $A_0$ is invertible.
Replacing  $\varphi$ by the map  $X\mapsto \varphi(A_0^{-1}X)$, we can assume that $A_0=I_n$.
By Lemma \ref{strong-degenerate-1}, there exist invertible matrices $P\in GL_{m'}(\mathbb{D}')$ and $Q\in GL_{n'}(\mathbb{D}')$ such that
$\varphi\left({\rm diag}(I_r, 0)\right)=P{\rm diag}(I_r, 0)Q$, $r=1,\ldots,n$.
Replacing  $\varphi$ by the map $X\mapsto P^{-1}\varphi(X)Q^{-1}$, we get that
\begin{equation}\label{vc32863wrdf66}
\varphi\left({\rm diag}(I_{r},0)\right)={\rm diag}\left(I_{r}, 0_{m'-r,n'-r}\right), \ \ r=1,\ldots, n.
\end{equation}

Write $E_{ij}=E_{ij}^{n\times n}$ and $E_{ij}'=E_{ij}^{m'\times n'}$.
 Note that $\varphi(E_{11})=E_{11}'$.  Since  $\varphi$ satisfies the Condition (I), from  Corollary \ref{Rectangular-PID2-7}
 we have either $\varphi(\mathcal{M}_1)\subseteq \mathcal{M}_1'$ with $\varphi(\mathcal{N}_1)\subseteq \mathcal{N}_1'$, or
 $\varphi(\mathcal{M}_1)\subseteq \mathcal{N}_1'$ with $\varphi(\mathcal{N}_1)\subseteq \mathcal{M}_1'$. We prove (\ref{d2345fgds53n}) only for the first case;
 the second case is similar. From now on we assume that
\begin{equation}\label{CXVWEEQ42fsd443}
\varphi(\mathcal{M}_1)\subseteq \mathcal{M}_1' \ \ {\rm and} \ \ \varphi(\mathcal{N}_1)\subseteq \mathcal{N}_1'.
\end{equation}

For any  $x\in \mathbb{D}^*$, by $\varphi(\mathcal{M}_1)\subseteq \mathcal{M}_1'$, we can assume that $\scriptsize\varphi(E_{11}+xE_{12})=\left(
              \begin{array}{ccc}
                1^* & x^{\sigma_1} & 0^* \\
                0 & 0& 0 \\
              \end{array}
            \right)$. Since  $E_{11}+xE_{12}\sim E_{11}+E_{22}$,
it follows from (\ref{vc32863wrdf66})  that
$\scriptsize{\rm rank}\left(
              \begin{array}{ccc}
                 1^*-1& x^{\sigma_1} & 0^* \\
                0&-1& 0 \\
              \end{array}\right)=1$. Thus $1^*=1$ and $0^*=0$.
Therefore,
\begin{equation}\label{v2164yr890bjo32}
\varphi(E_{11}+xE_{12})=E_{11}'+x^{\sigma_1}E_{12}', \ x\in \mathbb{D},
\end{equation}
where $\sigma_1: \mathbb{D}\rightarrow \mathbb{D}'$ is an injective map with $0^{\sigma_1}=0$.
Similarly, we have
\begin{equation}\label{164yryrrui5}
\varphi(E_{11}+xE_{21})=E_{11}'+x^{\mu_1}E_{21}', \ x\in \mathbb{D},
\end{equation}
where $\mu_1: \mathbb{D}\rightarrow \mathbb{D}'$ is an injective map with $0^{\mu_1}=0$.

For $2\leq k\leq n$ and $x_i, y_i\in \mathbb{D}$ ($i=2, \ldots, k$), we have  $E_{11}+x_2E_{12}+\cdots+x_kE_{1k}<{\rm diag}\left(I_k, 0\right)$ and
$E_{11}+y_2E_{21}+\cdots+y_kE_{k1}<{\rm diag}\left(I_k, 0\right)$.
Thus Corollary \ref{mYR25HLLsad23} and (\ref{vc32863wrdf66}) imply that  $\varphi(E_{11}+x_2E_{12}+\cdots+x_kE_{1k})<{\rm diag}\left(I_k, 0_{m'-k,n'-k}\right)$
and  $\varphi(E_{11}+y_2E_{21}+\cdots+y_kE_{k1})<{\rm diag}\left(I_k, 0_{m'-k,n'-k}\right)$.
Using (\ref{CXVWEEQ42fsd443}), it is easy to see that
\begin{equation}\label{dbcnbrer078rew}
\varphi(E_{11}+x_2E_{12}+\cdots+x_kE_{1k})=E_{11}'+x_2^*E_{12}'+\cdots+x_k^*E_{1k}', \ \ 2\leq k\leq n;
\end{equation}
\begin{equation}\label{dbcnb53iyo0789}
\varphi(E_{11}+y_2E_{21}+\cdots+y_kE_{k1})=E_{11}'+y_2^*E_{21}'+\cdots+y_k^*E_{k1}', \ \ 2\leq k\leq n,
\end{equation}
where  $x_j^*$, $y_j^*\in \mathbb{D}'$, $j=2,\ldots, k$. Moreover, by $\varphi(E_{11})=E_{11}'$ and the adjacency,
$(x_2, \ldots, x_k)\neq 0$ [resp. $^t(y_2, \ldots, y_k)\neq 0$] if and only if
$(x_2^*, \ldots, x_k^*)\neq 0$ [resp. $^t(y_2^*, \ldots, y_k^*)\neq 0$].

For any $1\leq r\leq n-2$ and  $2\leq k\leq n-r$,
we let $\scriptsize A_{k}=E_{11}^{k\times k}+x_2E_{12}^{k\times k}+\cdots+x_kE_{1k}^{k\times k}=\left(
      \begin{array}{cc}
        1 & \alpha_k \\
        0 & 0 \\
      \end{array}
    \right)$,
where $x_i\in \mathbb{D}$ and $\alpha_k:=(x_2,\ldots,x_k)\neq 0$. Then $A_k$  is a $k\times k$ idempotent matrices of rank one.
By  Corollary \ref{mYR25HLLsad23} and (\ref{vc32863wrdf66}), we have that
$$
\left(
          \begin{array}{cc}
            I_r & \\
             & 0 \\
          \end{array}
        \right)<\varphi\left(
          \begin{array}{ccc}
            I_r & &  \\
             & A_{k}& \\
             &&0\\
          \end{array}
        \right)<\left(
          \begin{array}{cc}
            I_{r+k} &  \\
             & 0 \\
          \end{array}
        \right),
$$
and ${\rm rank}\left(\varphi\left({\rm diag}\left(I_r, A_k, 0\right)\right)\right)=r+1$.
Write $\scriptsize\varphi\left({\rm diag}\left(I_r, A_k, 0\right)\right)=\left(
\begin{array}{cc}
 A_{11} & A_{12} \\
  A_{21} &   A_{22} \\
   \end{array}
   \right)$ where $A_{11}\in \mathbb{{D}'}^{(r+k)\times (r+k)}$.
By (\ref{vc32863wrdf66}), $\varphi\left({\rm diag}\left(I_r, A_k, 0\right)\right)$ is adjacent with both ${\rm diag}\left(I_{r+1}, 0\right)$
and ${\rm diag}\left(I_{r}, 0\right)$. It follows from Lemma \ref{Matrix-PID5-4bb} that $A_{22}=0$.
By $\varphi\left({\rm diag}\left(I_r, A_k, 0\right)\right)< {\rm diag}\left(I_{r+k}, 0\right)$ and
 Lemma \ref{minuspar1-3}, we  get that $A_{12}=0$,  $A_{21}=0$ and  $A_{11}$ is an idempotent matrix of rank $r+1$.
(Note that when $m'\neq n'$ and  using Lemma \ref{minuspar1-3}, we can get two suitable square matrices by adding some zero elements
on $\varphi\left({\rm diag}\left(I_r, A_k, 0\right)\right)$ and ${\rm diag}(I_{r+k},0)$.)
Thus $\varphi\left({\rm diag}\left(I_r, A_k, 0\right)\right)={\rm diag}\left(A_{11}, 0\right)$.

Write $\scriptsize A_{11}=\left(
 \begin{array}{cc}
  B_{11} & B_{12} \\
   B_{21} &  B_{22} \\
    \end{array}
    \right)$ where $B_{11}\in \mathbb{{D}'}^{r\times r}$ and $B_{22}\in \mathbb{{D}'}^{k\times k}$.
 Using ${\rm diag}\left(I_r,  0_{k}\right)<A_{11}$ and Lemma \ref{minuspar1-3},
we have similarly that $B_{12}=0$ and  $B_{21}=0$. By $B_{22}\neq 0$ and ${\rm diag}\left(I_r,  0_{k}\right)\sim A_{11}$, it is clear that $B_{11}=I_r$.
Therefore, we obtain
$$
\mbox{$\varphi\left({\rm diag}\left(I_r, A_k, 0\right)\right)={\rm diag}\left(I_r, B_{22}, 0\right)$, }
$$
where $B_{22}\in \mathbb{{D}'}^{k\times k}$ is an idempotent matrix of rank one.

 Since  ${\rm diag}\left(I_r, A_k, 0\right)\sim {\rm diag}\left(I_{r+1}, 0\right)$, ${\rm diag}\left(I_r, B_{22}, 0\right)\sim {\rm diag}\left(I_{r+1}, 0\right)$,
 and hence $B_{22}\sim E_{11}^{k\times k}$.
 By  Lemma \ref{Matrix-PID5-4bb}, we obtain either $\scriptsize B_{22}=\left(
  \begin{array}{cc}
   1 & \alpha_k^* \\
    0 & 0 \\
    \end{array}
     \right)$  or
 $\scriptsize B_{22}=\left(
 \begin{array}{cc}
 1 & 0 \\
 ^t\alpha_k^* & 0 \\
 \end{array}
 \right)$, where $0\neq \alpha_k^*\in {\mathbb{D}'}^{k-1}$.
Since $E_{11}+E_{r+1,1}\leq {\rm diag}\left(I_r, A_k, 0\right)$,  Corollary \ref{mYR25HLLsad23} implies that
$\small \varphi(E_{11}+E_{r+1,1})\leq \varphi\left({\rm diag}(I_r, A_k, 0)\right)$.
From (\ref{dbcnb53iyo0789}) we must have $\scriptsize B_{22}=\left(
                                                                \begin{array}{cc}
                                                                  1 & \alpha_k^* \\
                                                                  0 & 0 \\
                                                                \end{array}
                                                              \right)$.
Therefore, for any $0\neq \alpha_k\in \mathbb{D}^{k-1}$, there is $0\neq \alpha_k^*\in {\mathbb{D}'}^{k-1}$ such that
\begin{equation}\label{vcx24tyi890hfgj}
\varphi\left(
\begin{array}{ccc}
I_r &  &\\
 & \left(
      \begin{array}{cc}
        1 & \alpha_k\\
        0 & 0 \\
      \end{array}
    \right)&\\
&&0\\
\end{array}
\right)=
\left(\begin{array}{ccc}
I_r & & \\
 & \left(
      \begin{array}{cc}
        1 & \alpha_k^* \\
        0 & 0 \\
      \end{array}
    \right)& \\
&&0\\
\end{array}
\right), \ 1\leq r\leq n-2, \ 2\leq k\leq n-r.
\end{equation}

Similarly,  for any $0\neq \beta_k\in \,^{k-1}\mathbb{D}$, there is $0\neq \beta_k^*\in \,^{k-1}{\mathbb{D}'}$ such that
\begin{equation}\label{v24tyi890hfg}
\varphi\left(
\begin{array}{ccc}
I_r &  &\\
 & \left(
      \begin{array}{cc}
        1 & 0\\
        \beta_k & 0 \\
      \end{array}
    \right)&\\
&&0\\
\end{array}
\right)=
\left(\begin{array}{ccc}
I_r & & \\
 & \left(
      \begin{array}{cc}
        1 & 0 \\
        \beta_k^* & 0 \\
      \end{array}
    \right)& \\
&&0\\
\end{array}
\right), \ 1\leq r\leq n-2, \ 2\leq k\leq n-r.
\end{equation}

For any $2\leq k\leq n$, by  (\ref{dbcnbrer078rew}) we  write
\begin{equation}\label{b656tgfdgd8iads5}
\varphi\left(E_{11}+\sum_{j=2}^kE_{1j}\right)=E_{11}'+\sum_{j=2}^ka_{1j}^{(k)}E_{1j}',
\end{equation}
where $a_{1j}^{(k)}\in \mathbb{D}'$ and $(a_{12}^{(k)}, \ldots, a_{1k}^{(k)})\neq (0,\ldots,0)$.
From (\ref{v2164yr890bjo32}) we have $a_{12}^{(2)}\neq 0$.

We  prove $a_{1k}^{(k)}\neq 0$ ($3\leq k\leq n$) as follows. For any $3\leq k\leq n$ and $2\leq r< k$,
by (\ref{vcx24tyi890hfgj}), we  let
$$\varphi\left(E_{11}+\cdots +E_{rr}+\sum_{j=r+1}^kE_{rj}\right)=E_{11}'+\cdots +E_{rr}'+\sum_{j=r+1}^ka_{rj}^{(k)}E_{rj}',$$
where $a_{rj}^{(k)}\in \mathbb{D}'$ and $(a_{r,r+1}^{(k)},\ldots, a_{rk}^{(k)})\neq (0,\ldots, 0)$.

By (\ref{vc32863wrdf66}) and the adjacency, when $r=k-1$ we have $a_{k-1,k}^{(k)}\neq 0$.
We show $a_{rk}^{(k)}\neq 0$, $r=2, \ldots, k-2$.

Since
$\varphi\left(E_{11}+\cdots +E_{k-1,k-1}+E_{k-1,k}\right)\sim \varphi\left(E_{11}+\cdots +E_{k-2,k-2}+E_{k-2,k-1}+E_{k-2,k}\right),$
we get
$$E_{11}'+\cdots +E_{k-1,k-1}'+a_{k-1,k}^{(k)}E_{k-1,k}'\sim E_{11}'+\cdots +E_{k-2,k-2}'+a_{k-2,k-1}^{(k)}E_{k-2,k-1}'+a_{k-2,k}^{(k)}E_{k-2,k}',$$
which implies that $\small{\rm rank}\left(
                                \begin{array}{cc}
                                  a_{k-2,k-1}^{(k)} & a_{k-2,k}^{(k)} \\
                                  1 & a_{k-1,k}^{(k)} \\
                                \end{array}
                              \right)=1$.  By $(a_{k-2,k-1}^{(k)}, a_{k-2,k}^{(k)})\neq (0,0)$ and $a_{k-1,k}^{(k)}\neq 0$,
we have  $a_{k-2,k}^{(k)}\neq 0.$

Since
$\varphi\left(E_{11}+\cdots +E_{k-2,k-2}+E_{k-2,k-1}+E_{k-2,k}\right)\sim \varphi\left(E_{11}+\cdots +E_{k-3,k-3}+\sum_{j=k-2}^kE_{k-3,j}\right),$
one has
$$E_{11}'+\cdots +E_{k-2,k-2}'+a_{k-2,k-1}^{(k)}E_{k-2,k-1}'+a_{k-2,k}^{(k)}E_{k-2,k}'\sim E_{11}'+\cdots +E_{k-3,k-3}'+\sum_{j=k-2}^ka_{k-3,j}^{(k)}E_{k-3,j}'.$$
Similarly, we  have $a_{k-3,k}^{(k)}\neq 0.$
In the similar way, we can prove that $a_{k-4,k}^{(k)}\neq 0$, $\ldots$, $a_{2k}^{(k)}\neq 0$.

Since
$\varphi\left(E_{11}+\sum_{j=2}^kE_{1j}\right)\sim \varphi\left(E_{11}+E_{22}+\sum_{j=3}^kE_{2j}\right),$
$E_{11}'+\sum_{j=2}^ka_{1j}^{(k)}E_{1j}'\sim E_{11}'+E_{22}'+\sum_{j=3}^ka_{2j}^{(k)}E_{2j}'.$
It follows that
$${\rm rank}\left(
                                \begin{array}{cccc}
                                  a_{12}^{(k)} & a_{13}^{(k)}&\cdots & a_{1k}^{(k)}\\
                                  1 & a_{23}^{(k)} &\cdots & a_{2k}^{(k)}\\
                                \end{array}
                              \right)=1.$$
  By  $(a_{12}^{(k)},\ldots,  a_{1k}^{(k)})\neq (0,\ldots,0)$ and $a_{2k}^{(k)}\neq 0$, we must have $a_{1k}^{(k)}\neq 0$.
Therefore, we have proved that
$$a_{1k}^{(k)}\neq 0, \ \ k=2,\ldots,n.$$

Let  $\gamma_1=E_{11}$, $\gamma_k=E_{11}+\sum_{j=2}^kE_{1j}$, $\gamma_1^*=E_{11}'$ and
$\gamma_k^*=E_{11}'+\sum_{j=2}^ka_{1j}^{(k)}E_{1j}'$, $k=2,\ldots,n$. Then
$\gamma_i\in \mathcal{M}_1$, $\gamma_i^*\in \mathcal{M}_1'$, $\gamma_1, \ldots, \gamma_n$ are left linearly independent over $\mathbb{D}$,
and $\gamma_1^*, \ldots, \gamma_n^*$ are left linearly independent over $\mathbb{D}'$.
By (\ref{b656tgfdgd8iads5}) we have $\varphi(\gamma_k)=\gamma_k^*$, $k=1,\ldots,n$, and hence  ${\rm dim}(\varphi(\mathcal{M}_1))\geq n$.
By Lemma \ref{degenerate-2}, we get  ${\rm dim}(\varphi(\mathcal{M}_1))\leq n$.
Hence ${\rm dim}(\varphi(\mathcal{M}_1))=n$.

Similarly, we can prove ${\rm dim}(\varphi(\mathcal{N}_1))=n$. Hence ${\rm dim}(\varphi(\mathcal{M}_1))={\rm dim}(\varphi(\mathcal{N}_1))=n.$
$\qed$

\begin{lem}\label{slynsf32524asd} \ Let $\mathbb{D}, \mathbb{D}'$ be division rings,
and let $m',n',n$ be integers with $m', n'\geq n\geq 2$. Suppose that
$\varphi:  {\mathbb{D}}^{n\times n}\rightarrow  {\mathbb{D}'}^{m'\times n'}$
is a graph homomorphism with $\varphi(0)=0$,  and there exists $A_0\in {\mathbb{D}}^{n\times n}$ such that ${\rm rank}(\varphi(A_0))=n$.
 Assume further that $\varphi$ does not satisfies the Condition (I). Then  there exist two
maximal sets ${\cal M}'$ and ${\cal N}'$ of different types in ${\mathbb{D}'}^{m'\times n'}$ such that
 $\varphi\left(\mathbb{D}_{\leq 1}^{n\times n}\right)\subseteq {\cal M}'\cup {\cal N}'$. Moreover, if $|\mathbb{D}|\geq 4$,
 then $\varphi\left(\mathbb{D}_{\leq 1}^{n\times n}\right)$  is  an adjacent set.
\end{lem}
\proof  {\bf Step 1.} \
Since  $\varphi$ does not satisfies the Condition (I),   there are two maximal sets  ${\cal M}$ and ${\cal R}$ of different types
in ${\mathbb{D}}^{m\times n}$ with $0\in{\cal M}\cap {\cal R}$, such that  $\varphi({\cal M})\subseteq {\cal M}'$ and $\varphi({\cal R})\subseteq {\cal R}'$,
where ${\cal M}'$ and ${\cal R}'$ are two maximal sets of the same type containing $0$ in ${\mathbb{D}'}^{m'\times n'}$.
Since $|{\cal M}\cap {\cal R}|\geq 2$, $|{\cal M}'\cap {\cal R}'|\geq 2$.  It follows from Lemma \ref{Rectangular-1-13} that ${\cal M}'={\cal R}'$. Thus
\begin{equation}\label{gf32fdsmvzqjj}
\mbox{$\varphi({\cal M})\subseteq {\cal M}'$ \ and \ $\varphi({\cal R})\subseteq {\cal M}'$.}
\end{equation}
Without loss of generality, we assume that both $\mathcal{M}$ and $\mathcal{M}'$ are of type one. Hence $\mathcal{R}$ is of type two.
By Lemma \ref{non-degeneratelemma00a}(b), every non-degenerate graph homomorphism satisfies the Condition (I). Thus $\varphi$ is degenerate.

 We prove that there exist a type two maximal set ${\cal N}'$ such that
 $\varphi\left(\mathbb{D}_{\leq 1}^{n\times n}\right)\subseteq {\cal M}'\cup {\cal N}'$ as follows.
Suppose  $\varphi\left(\mathbb{D}_{\leq 1}^{n\times n}\right)$ is  an adjacent set. By (\ref{gf32fdsmvzqjj}) and Lemma \ref{Recastatre43}, it is easy to see that
 $\varphi\left(\mathbb{D}_{\leq 1}^{n\times n}\right)\subseteq {\cal M}'$.
Now, suppose  $\varphi\left(\mathbb{D}_{\leq 1}^{n\times n}\right)$  is not an adjacent set.
 Then there exists a type one maximal set ${\cal N}$  containing $0$ in ${\mathbb{D}}^{n\times n}$, such that ${\cal M}\neq {\cal N}$,
$\varphi({\cal N})\subseteq {\cal N}'$ and ${\cal N}'\neq {\cal M}'$, where ${\cal N}'$ is a  maximal set containing $0$ in ${\mathbb{D}'}^{m'\times n'}$.
Otherwise, for any type one maximal set ${\cal N}$  containing $0$ in ${\mathbb{D}}^{n\times n}$,
we have $\varphi({\cal N})\subseteq {\cal M}'$. Since every matrix in $\mathbb{D}_{\leq 1}^{n\times n}$ is contained in some type one maximal set containing $0$,
we get $\varphi\left(\mathbb{D}_{\leq 1}^{n\times n}\right)\subseteq {\cal M}'$, a contradiction to assumption.
By Corollary \ref{Rectangular-1-13}, we have $|{\cal N}\cap {\cal R}|\geq 2$, hence
$|{\cal N}'\cap {\cal M}'|\geq 2$. It follows that ${\cal N}'$ is of type two.
Let ${\cal S}$ be any type two  maximal set containing $0$ in ${\mathbb{D}}^{n\times n}$, and let $\varphi({\cal S})\subseteq {\cal S}'$
where ${\cal S}'$ is a maximal set containing $0$ in ${\mathbb{D}'}^{m'\times n'}$.
By Corollary \ref{Rectangular-1-13}, we have
 $|{\cal S}\cap {\cal M}|\geq 2$ and $|{\cal S}\cap {\cal N}|\geq 2$,    hence
$|{\cal S}'\cap {\cal M}'|\geq 2$ and $|{\cal S}'\cap {\cal N}'|\geq 2$. Applying Corollary \ref{Rectangular-1-13} again, we get either
${\cal S}'={\cal M}'$ or ${\cal S}'={\cal N}'$, and hence ${\cal S}'\subseteq {\cal M}'\cup {\cal N}'$.
Consequently,  $\varphi({\cal S})\subseteq {\cal M}'\cup {\cal N}'$ for any type two maximal set ${\cal S}$ containing $0$ in ${\mathbb{D}}^{n\times n}$.
Since every matrix in $\mathbb{D}_{\leq 1}^{n\times n}$ is contained in some type two maximal set containing $0$,
we get $\varphi\left(\mathbb{D}_{\leq 1}^{n\times n}\right)\subseteq {\cal M}'\cup {\cal N}'$.
Therefore,  we always have
 \begin{equation}\label{fs25nvm58gjg}
 \varphi\left(\mathbb{D}_{\leq 1}^{n\times n}\right)\subseteq {\cal M}'\cup {\cal N}'.
 \end{equation}

{\bf Step 2.} \
 From now on we assume $|\mathbb{D}|\geq 4$. In this step, we will prove that $\varphi\left(\mathbb{D}_{\leq 1}^{n\times n}\right)$ is  an adjacent set.

By Lemma \ref{maximalset022}(i), there are two invertible matrices $P_1$ and $Q_1$ over $\mathbb{D}$, such that
$\mathcal{M}=P_1\mathcal{M}_1Q_1$ and $\mathcal{R}=P_1\mathcal{N}_1Q_1$. Also, there is an invertible matrix $P_2$ over $\mathbb{D}'$, such that
$\mathcal{M}'=P_2\mathcal{M}_1'$. Replacing  $\varphi$ by the map $X\mapsto P_2^{-1}\varphi(P_1XQ_1)$, (\ref{gf32fdsmvzqjj}) becomes
\begin{equation}\label{fgd3ncnit9}
\mbox{$\varphi({\cal M}_1)\subseteq {\cal M}_1'$ \ \  and \ \ $\varphi({\cal N}_1)\subseteq {\cal M}_1'$.}
\end{equation}

 Suppose that $\varphi\left(\mathbb{D}_{\leq 1}^{n\times n}\right)$  is not an adjacent set. We show a contradiction as follows.

 There exists a type one maximal set ${\cal R}_1$   containing $0$ in ${\mathbb{D}}^{n\times n}$ with ${\cal R}_1\neq {\cal M}_1$, such that
$\varphi({\cal R}_1)\nsubseteq {\cal M}_1'$. Otherwise, for any type one maximal set ${\cal R}_1$  containing $0$ in ${\mathbb{D}}^{n\times n}$,
we have $\varphi({\cal R}_1)\subseteq {\cal M}_1'$. Since every matrix in $\mathbb{D}_{\leq 1}^{n\times n}$ is contained in some type one maximal set containing $0$,
we get $\varphi\left(\mathbb{D}_{\leq 1}^{n\times n}\right)\subseteq {\cal M}_1'$, a contradiction.
By (\ref{fs25nvm58gjg}), it is easy to see that
$$\varphi({\cal R}_1)\subseteq{\cal N}'.$$
Using Lemma \ref{maximalset022}(ii),  there is $P_3\in GL_{m}(\mathbb{D})$ such that ${\cal R}_1=P_3{\cal M}_2$ and ${\cal M}_1=P_3{\cal M}_1$.
On the other hand, there is $Q_2\in GL_{n'}(\mathbb{D}')$ such that ${\cal N}'={\cal N}'_1Q_2$. Modifying the map $\varphi$ by the map
$X\mapsto \varphi(P_3X)Q_2^{-1}$.
We get  $\varphi({\cal M}_2)\subseteq {\cal N}_1'$ and (\ref{fgd3ncnit9}) holds. Therefore, we may assume that
\begin{equation}\label{VCXWR5HF34sds}
\mbox{$\varphi({\cal M}_2)\subseteq {\cal N}_1'$, \ $\varphi({\cal M}_1)\subseteq {\cal M}_1'$ \ \  and \ \ $\varphi({\cal N}_1)\subseteq {\cal M}_1'$.}
\end{equation}

By the conditions of this lemma, there exists $A_0\in GL_n(\mathbb{D})$ such that ${\rm rank}(\varphi(A_0))=n$.
Write
$\scriptsize A_0=\left(
                        \begin{array}{c}
                          \alpha_1 \\
                          \vdots \\
                          \alpha_n \\
                        \end{array}
                      \right)=\left(\beta_1,\ldots,\beta_{n}\right)$  where  $\alpha_{i}\in {\mathbb{D}}^{n}$ and $\beta_{i}\in \,^{n}{\mathbb{D}}$;
$\scriptsize\varphi(A_0)=\left(
                 \begin{array}{cc}
                   A_{11} & A_{12}\\
                   A_{21} & A_{22}\\
                 \end{array}
               \right)$ where  $A_{11}\in {\mathbb{D}'}^{1\times 1}$ and $A_{22}\in {\mathbb{D}'}^{(m'-1)\times (n'-1)}$.
Let $\scriptsize B=\left(
         \begin{array}{c}
           \alpha_1 \\
           0_{n-1,n} \\
         \end{array}
       \right)$,  $\small C=\left(\beta_1, 0, \ldots, 0\right)$ and $\scriptsize J_n=\left(
   \begin{array}{ccc}
     0 & 0 & 1 \\
     0 & I_{n-2} & 0 \\
     1 & 0 & 0 \\
   \end{array}
 \right)$.
Then $B <A_0$ and $C<A_0$. By  Corollary \ref{mYR25HLLsad23},  we have $\varphi(B)<\varphi(A_0)$ and $\varphi(C)<\varphi(A_0)$.
Since $\varphi(B)\in \mathcal{M}_1'$ and $\varphi(C)\in \mathcal{N}_1'$, it is clear that
\begin{equation}\label{f53fsfsxbc235}
{\rm rank}\left(A_{21}, A_{22}\right)={\rm rank}\left(
\begin{array}{c}
A_{12} \\
A_{22}\\
\end{array}
\right)=n-1.
\end{equation}

By (\ref{f53fsfsxbc235}) and ${\rm rank}(\varphi(A_0))=n$, we have  either $A_{22}\cong {\rm diag}(I_{n-1}, 0)$
or $A_{22}\cong {\rm diag}(I_{n-2}, 0)$.  By appropriate  elementary row
and column operations on matrices, we can obtain either ${\rm diag}(I_{n}, 0)$  (if $A_{22}\cong {\rm diag}(I_{n-1}, 0)$) or
${\rm diag}(J_n, 0)$ (if $A_{22}\cong {\rm diag}(I_{n-2}, 0)$) from $\varphi(A_0)$.  Since our matrix elementary
  operations  do not change ${\cal N}_1'$ and  ${\cal M}_1'$,  (\ref{VCXWR5HF34sds})  still holds.
Thus, without loss of generality we may assume  either $\varphi(A_0)={\rm diag}(I_{n}, 0)$ or $\varphi(A_0)={\rm diag}(J_n, 0)$.

{\em Case 1.} \  $\varphi(A_0)={\rm diag}(I_{n}, 0)$.
Since $\scriptsize\left(
\begin{array}{c}
0\\
\alpha_2+\lambda\alpha_1 \\
 0_{n-2, n} \\
  \end{array}
  \right)< A_0$ for all $\lambda\in \mathbb{D}^*$,  Corollary \ref{mYR25HLLsad23} implies that
$\scriptsize \varphi\left(
\begin{array}{c}
0\\
\alpha_2+\lambda\alpha_1 \\
 0_{n-2,\, n} \\
  \end{array}\right)<\varphi(A_0)$, $\lambda\in \mathbb{D}^*$. Thus $\varphi({\cal M}_2)\subseteq {\cal N}_1'$ implies that
$\scriptsize\varphi\left(
\begin{array}{c}
0\\
\alpha_2+\lambda\alpha_1 \\
 0_{n-2,\, n} \\
  \end{array}\right)=\left(
                       \begin{array}{cc}
                         1&0 \\
                          Z_\lambda &0_{m'-1,  n'-1} \\
                       \end{array}
                     \right)$, where $Z_\lambda\in \,^{m'-1}{\mathbb{D}'}$, $\lambda\in \mathbb{D}^*$.
Since $\scriptsize\varphi\left(
\begin{array}{c}
\lambda^{-1}\alpha_2+\alpha_1 \\
0\\
 0_{n-2, n} \\
  \end{array}\right)<\varphi(A_0)$ for all $\lambda\in \mathbb{D}^*$,  it follows from $\varphi({\cal M}_1)\subseteq{\cal M}_1'$ that
 $\scriptsize\varphi\left(
\begin{array}{c}
\lambda^{-1}\alpha_2+\alpha_1 \\
0\\
 0_{n-2, n} \\
  \end{array}\right)=\left(
                       \begin{array}{cc}
                         1 & T_\lambda\\
                         0 & 0 \\
                       \end{array}
                     \right)$,
where $T_\lambda\in {\mathbb{D}'}^{n'-1}$, $\lambda\in \mathbb{D}^*$. By $|\mathbb{D}^*|\geq 3$ and
the adjacency, there is $\lambda_0\in\mathbb{D}^*$
such that $Z_{\lambda_0}\neq 0$ and $T_{\lambda_0}\neq 0$.
By $\scriptsize\varphi\left(
                \begin{array}{c}
                  \lambda^{-1}_0\alpha_2+\alpha_1 \\
                  0 \\
                  0\\
                \end{array}
              \right)\sim  \varphi\left(
\begin{array}{c}
0\\
\alpha_2+\lambda_0\alpha_1 \\
 0_{n-2, n} \\
  \end{array}\right)$, we obtain $\scriptsize\left(
                       \begin{array}{cc}
                         1 & T_{\lambda_0}\\
                         0 & 0 \\
                       \end{array}
                     \right)\sim \left(
                       \begin{array}{cc}
                         1&0 \\
                         Z_{\lambda_0} &0_{m'-1,  n'-2} \\
                       \end{array}
                     \right)$, a contradiction.

{\em Case 2.} \  $\varphi(A_0)={\rm diag}(J_n, 0)$.
By $\scriptsize\left(
         \begin{array}{c}
           \alpha_1+\alpha_2 \\
           0\\
           0_{n-2, n} \\
         \end{array}
       \right)< A_0$ and Corollary \ref{mYR25HLLsad23},
$\scriptsize \varphi\left(
         \begin{array}{c}
           \alpha_1+\alpha_2 \\
           0\\
           0_{n-2, n} \\
         \end{array}
       \right)<\varphi(A_0)$. It follows from $\varphi({\cal M}_1)\subseteq {\cal M}_1'$ that
$\scriptsize\varphi\left(
         \begin{array}{c}
           \alpha_1+\alpha_2 \\
           0\\
           0_{n-2, n} \\
         \end{array}
       \right)=\left(
       \begin{array}{ccc}
        \alpha &1 &0\\
         0&0 &0\\
          0&0 &0_{m'-n,n'-n}\\
         \end{array}
        \right)$ where $\alpha\in {\mathbb{D}'}^{n-1}$.
Since $\scriptsize\left(
\begin{array}{c}
0\\
\alpha_1+\alpha_2 \\
 0_{n-2,\, n} \\
  \end{array}
  \right)< A_0$ and Corollary \ref{mYR25HLLsad23},
$\scriptsize\varphi\left(
\begin{array}{c}
0\\
\alpha_1+\alpha_2 \\
 0_{n-2, n} \\
  \end{array}\right)<\varphi(A_0)$. By $\varphi({\cal M}_2)\subseteq {\cal N}_1'$, we get
$\scriptsize\varphi\left(
\begin{array}{c}
0\\
\alpha_1+\alpha_2 \\
 0_{n-2, n} \\
  \end{array}\right)=\left(
                       \begin{array}{ccc}
                         \beta & 0 &0\\
                         1 & 0&0 \\
                         0&0&0_{m'-n,n'-n}\\
                       \end{array}
                     \right)$  where  $\beta\in \,^{n-1}\mathbb{D}'$.
Since $\small\varphi\left(
                \begin{array}{c}
                \alpha_1+\alpha_2 \\
                  0 \\
                  0_{n-2,n}\\
                \end{array}
              \right)\sim  \varphi\left(
\begin{array}{c}
0\\
\alpha_1+\alpha_2 \\
 0_{n-2,n} \\
  \end{array}\right)$,
 $\small\left(
       \begin{array}{ccc}
        \alpha &1 &0\\
         0&0 &0\\
          0&0 &0_{m'-n,n'-n}\\
         \end{array}
        \right)\sim \left(
                       \begin{array}{ccc}
                         \beta & 0 &0\\
                         1 & 0&0 \\
                         0&0&0_{m'-n,n'-n}\\
                       \end{array}
                     \right),$
  a contradiction.

Thus both Case 1 and Case 2  cannot occur. It follows that $\varphi\left(\mathbb{D}_{\leq 1}^{n\times n}\right)$  is  an adjacent set.
$\qed$

\begin{lem}\label{slycxzwe46e4dfz} \ Let $\mathbb{D}, \mathbb{D}'$ be division rings with $|\mathbb{D}|\geq 4$,
and let $m',n',n$ be integers with $m', n'\geq n\geq 2$. Suppose that
$\varphi:  {\mathbb{D}}^{n\times n}\rightarrow  {\mathbb{D}'}^{m'\times n'}$
is a graph homomorphism with $\varphi(0)=0$,  and there exists $A_0\in {\mathbb{D}}^{n\times n}$ such that ${\rm rank}(\varphi(A_0))=n$.
Assume further  that $\varphi$ satisfies the Condition (I). Then $\varphi$ is non-degenerate.
\end{lem}
\proof \ We prove that $\varphi$ is non-degenerate by contradiction.
Suppose that $\varphi$ is degenerate. Then there exists a matrix $A\in\mathbb{D}^{n\times n}_{\leq 1}$ and there are two
maximal sets $\mathcal{R}$ and $\mathcal{S}$ of different types in ${\mathbb{D}'}^{m'\times n'}$ such that
\begin{equation}\label{ghf54sbxbnpa12}
\mbox{$\varphi\left(\mathbb{B}_A\right)\subseteq \mathcal{R}\cup \mathcal{S}$ \ with \ $\varphi(A)\in\mathcal{R}\cap\mathcal{S}$.}
\end{equation}
Without loss of generality, we assume that $\mathcal{R}$ is of type one and  $\mathcal{S}$ is of type two.
 By ${\rm rank}(\varphi(A_0))=n$ and (\ref{uy87mbmm}), we have ${\rm rank}(A_0)=n$. Since $d(\varphi(A), 0)\leq 1$, we have either $d(\varphi(A_0),\varphi(A))=n$ or  $d(\varphi(A_0),\varphi(A))=n-1$.
 Note either $d(A_0,A)=n$ or $d(A_0,A)=n-1$. Let $E_{ij}=E_{ij}^{n\times n}$ and $E_{ij}'=E_{ij}^{m'\times n'}$.
We distinguish the following cases to show a contradiction.

{\bf Case 1.} \ $d(A_0,A)=n$ and $d(\varphi(A_0),\varphi(A))=n-1$. Then  ${\rm rank}(\varphi(A))={\rm rank}(A)=1$ and $\varphi(A)<\varphi(A_0)$.
By Lemma \ref{minuspar1}(c), there are invertible matrices $P_1,Q_1$ over $\mathbb{D}'$ such that
$\varphi(A)=P_1E_{11}'Q_1$  and  $\varphi(A_0)=P_1{\rm diag}(I_n,0)Q_1$.
Also, there are invertible matrices $P_2,Q_2$ over $\mathbb{D}$ such that $P_2AQ_2=E_{11}$. Let
$\scriptsize P_2A_0Q_2=\left(
                 \begin{array}{cc}
                   A_{11} & A_{12}\\
                   A_{21} & A_{22}\\
                 \end{array}
               \right)$  where $A_{11}\in {\mathbb{D}}^{1\times 1}$ and $A_{22}\in {\mathbb{D}}^{(n-1)\times (n-1)}$,
Put
$\scriptsize J_n=\left(
   \begin{array}{ccc}
     0 & 0 & 1 \\
     0 & I_{n-2} & 0 \\
     1 & 0 & 0 \\
   \end{array}
 \right)$.

 By ${\rm rank}(A_0)={\rm rank}(A_0-A)=n$, it is clear that $A_{22}\cong I_{n-1}$
or $A_{22}\cong {\rm diag}(I_{n-2}, 0)$.
Using appropriate  elementary row and column operations on matrices, we can obtain either ${\rm diag}(a, I_{n-1})$ where $a\notin \{0,1\}$ (if  $A_{22}\cong I_{n-1}$)
or $J_n$ (if $A_{22}\cong {\rm diag}(I_{n-2}, 0)$) from  $P_2A_0Q_2$. Moreover, $P_2AQ_2=E_{11}$ is  unchanged.
Thus, there are invertible matrices $P_3,Q_3$ over $\mathbb{D}$, such that $P_3AQ_3=E_{11}$  and either  $P_3A_0Q_3={\rm diag}(a, I_{n-1})$
or $P_3A_0Q_3=J_n$.
Replacing $\varphi$ by the map $X\mapsto P_1^{-1}\varphi(P_3^{-1}XQ_3^{-1})Q_1^{-1}$, we have that
$A=E_{11}$, and either  $A_0={\rm diag}(a, I_{n-1})$  or $A_0=J_n$. Moreover,
$\varphi(E_{11})=E_{11}'$  and  $\varphi({\rm diag}(a, I_{n-1}))={\rm diag}(I_n,0)$.
On the other hand, (\ref{ghf54sbxbnpa12}) becomes
$\varphi\left(\mathbb{D}_{\leq 1}^{n\times n}+E_{11}\right)\subseteq \mathcal{R}\cup \mathcal{S}$  with  $\varphi(E_{11})=E_{11}'\in\mathcal{R}\cap\mathcal{S}$.

Recall that $\varphi$ satisfies the Condition (I). By $\varphi(E_{11})=E_{11}'$, $\varphi(0)=0$ and Corollary \ref{Rectangular-PID2-7},  we may assume with no loss of generality that
\begin{equation}\label{VCTR24tfdgh67}
\varphi(\mathcal{M}_1)\subseteq \mathcal{M}_1', \ \ \varphi(\mathcal{N}_1)\subseteq \mathcal{N}_1'.
\end{equation}
Since $\mathcal{M}_1, \mathcal{N}_1\subseteq \mathbb{D}_{\leq 1}^{n\times n}+E_{11}$,
we get that $\varphi(\mathcal{M}_1)\subseteq {\cal M}'_1\cap (\mathcal{R}\cup \mathcal{S})$
and $\varphi(\mathcal{N}_1)\subseteq {\cal N}'_1\cap (\mathcal{R}\cup \mathcal{S})$.
Suppose that ${\cal M}'_1\neq\mathcal{R}$. Then ${\cal M}'_1\cap (\mathcal{R}\cup \mathcal{S})={\cal M}'_1\cap \mathcal{S}$, thus
$\ell:={\cal M}'_1\cap \mathcal{S}$ is a line in $AG({\cal M}'_1)$ and $\varphi({\cal M}'_1)\subseteq \ell$, a contradiction to Lemma \ref{Recastatre43}.
Therefore, we must have  ${\cal M}'_1=\mathcal{R}$. Similarly,  ${\cal N}'_1=\mathcal{S}$. It follows that
\begin{equation}\label{VC2EHFH46fghbv}
\mbox{$\varphi\left(\mathbb{D}_{\leq 1}^{n\times n}+E_{11}\right)\subseteq \mathcal{M}_1'\cup \mathcal{N}_1'$.}
\end{equation}

{\em Subcase 1.1.} \ $A_0={\rm diag}(a, I_{n-1})$. By (\ref{VCTR24tfdgh67}), we have $\varphi(aE_{11})=a^*E_{11}'$ where $a^*\in {\mathbb{D}'}^*$.
Since $aE_{11}<A_0$, Corollary \ref{mYR25HLLsad23} implies that $\varphi(aE_{11})<\varphi(A_0)={\rm diag}(I_n,0)$. Consequently, $a^*=1$ and hence
$\varphi(aE_{11})=E_{11}'$, a contradiction to $\varphi(aE_{11})\sim \varphi(E_{11})$.

{\em Subcase 1.2.} \ $A_0=J_n$. By (\ref{VC2EHFH46fghbv}), we get $\varphi\left(E_{nn}+E_{11}\right)\in\mathcal{M}_1'\cup \mathcal{N}_1'$.
Without loss of generality, we assume that $\varphi\left(E_{nn}+E_{11}\right)\in\mathcal{M}_1'$. Since $d(E_{nn}+E_{11}, J_n)=n-1$,
 $d\left(\varphi(E_{nn}+E_{11}), {\rm diag}(I_n,0)\right)\leq n-1$. It follows that
 $\scriptsize\varphi(E_{nn}+E_{11})=\left(
                            \begin{array}{cc}
                              1 & \alpha \\
                              0 & 0 \\
                            \end{array}
                          \right)$, where $0\neq \alpha\in {\mathbb{D}'}^{n'-1}$ because $\varphi(E_{nn}+E_{11})\sim \varphi(E_{11})=E_{11}'$.
Since $d(E_{n1}+E_{11}, J_n)=n-1$,
 $d\left(\varphi(E_{n1}+E_{11}), {\rm diag}(I_n,0)\right)\leq n-1$. Thus $ \varphi(\mathcal{N}_1)\subseteq \mathcal{N}_1'$ implies that
 $\scriptsize \varphi(E_{n1}+E_{11})=\left(
                            \begin{array}{cc}
                              1 & 0\\
                              \beta & 0 \\
                            \end{array}
                          \right)$, where $0\neq \beta\in \,^{m'-1}{\mathbb{D}'}$ because $\varphi(E_{n1}+E_{11})\sim \varphi(E_{11})$.
Since $\varphi(E_{nn}+E_{11})\sim \varphi(E_{n1}+E_{11})$, one gets
$\scriptsize\left(
\begin{array}{cc}
1 & \alpha \\
 0 & 0 \\
 \end{array}
\right)\sim \left(
                            \begin{array}{cc}
                              1 & 0\\
                              \beta & 0 \\
                            \end{array}
                          \right)$, a contradiction.

By Subcases 1.1-1.2, we always have a contradiction in Case 1.

{\bf Case 2.} \ $d(A_0,A)=n$ and $d(\varphi(A_0),\varphi(A))=n$. Let $\psi(X)=\varphi(X+A)-\varphi(A)$, $X\in{\mathbb{D}}^{n\times n}$. Then
$\psi:  {\mathbb{D}}^{n\times n}\rightarrow  {\mathbb{D}'}^{m'\times n'}$
is a graph homomorphism with $\psi(0)=0$  and  ${\rm rank}(\psi(A_0-A))=n$. Moreover, from (\ref{ghf54sbxbnpa12}) we have
\begin{equation}\label{vxcqe25mbfy67}
\mbox{ $\psi\left(\mathbb{D}_{\leq 1}^{n\times n}\right)\subseteq \mathcal{R}'\cup \mathcal{S}'$
\ \  with \ \ $0\in\mathcal{R}'\cap\mathcal{S}'$,}
 \end{equation}
 where $\mathcal{R}'=\mathcal{R}-\varphi(A)$ is a maximal set of type one, and $\mathcal{S}'=\mathcal{S}-\varphi(A)$ is a maximal set of type two.

{\em Subcase 2.1.} \  $\psi$ satisfies the Condition (I). Then, by Lemma \ref{slynon-degenerate-2} we have
\begin{equation}\label{645dgdh9js53n}
{\rm dim}(\psi(\mathcal{M}_1))={\rm dim}(\psi(\mathcal{N}_1))=n.
\end{equation}
By Lemma \ref{Rectangular-PID2-4}, without loss of generality,  we may assume either $\varphi(\mathcal{M}_1)\subseteq \mathcal{M}_1'$ with $\varphi(\mathcal{N}_1)\subseteq \mathcal{N}_1'$, or
 $\varphi(\mathcal{M}_1)\subseteq \mathcal{N}_1'$ with $\varphi(\mathcal{N}_1)\subseteq \mathcal{M}_1'$. We prove this result only for the first case;
 the second case is similar. From now on we assume that
\begin{equation}\label{54fvs72fsd443}
\psi(\mathcal{M}_1)\subseteq \mathcal{M}_1' \ \ {\rm and} \ \ \psi(\mathcal{N}_1)\subseteq \mathcal{N}_1'.
\end{equation}

By (\ref{vxcqe25mbfy67}) and (\ref{54fvs72fsd443}), one has $\psi(\mathcal{M}_1)\subseteq \mathcal{M}_1'\cap(\mathcal{R}'\cup\mathcal{S}')$
and $\psi(\mathcal{N}_1)\subseteq \mathcal{N}_1'\cap(\mathcal{R}'\cup\mathcal{S}')$. Clearly, $\mathcal{M}_1'\cap\mathcal{R}'=\{0\}=\mathcal{N}_1'\cap\mathcal{S}'$.
We assert $\mathcal{M}_1'=\mathcal{R}'$.  Otherwise,  we have that $\mathcal{M}_1'\neq\mathcal{R}'$ and
 $\psi(\mathcal{M}_1)\subseteq \mathcal{M}_1'\cap(\mathcal{R}'\cup\mathcal{S}')=\mathcal{M}_1'\cap\mathcal{S}'$.
   Hence Lemma \ref{Rectangular-PID2-11} implies that $\ell:=\mathcal{M}_1'\cap\mathcal{S}'$ is a line in $AG({\cal M}'_1)$,  a contradiction to Lemma \ref{Recastatre43}.
 Similarly,  $\mathcal{N}_1'=\mathcal{S}'$. Therefore,
\begin{equation}\label{5649kgaafsd43}
\psi\left(\mathbb{D}_{\leq 1}^{n\times n}\right)\subseteq \mathcal{M}'_1\cup \mathcal{N}_1'.
\end{equation}

By $\psi$ satisfying the Condition (I) and $\psi(\mathcal{M}_1)\subseteq \mathcal{M}_1'$, $\psi$ maps any type two maximal set containing $0$
into a type two maximal set.  Thus (\ref{5649kgaafsd43}) and Lemma \ref{Recastatre43} imply that $\psi$ maps any type two maximal set containing $0$
into $\mathcal{N}_1'$. By (\ref{645dgdh9js53n}), let $\{B_1, \ldots,B_n\}$ be a maximal left linear independent subset in $\psi(\mathcal{M}_1)$,
 and let $\psi(C_i)=B_i$ where $C_i\in \mathcal{M}_1$, $i=1, \ldots, n$. Since $C_i\sim 0$ and $\psi$ satisfies the Condition (I),
 there is a type two maximal set $\mathcal{R}_i$ containing $C_i$ and $0$ in $\mathbb{D}^{n\times n}$,  such that
 $\psi(\mathcal{R}_i)\subseteq \mathcal{N}_1'$,  $i=1, \ldots, n$. Thus, $\{B_1, \ldots,B_n\}\subset (\mathcal{M}_1'\cap\mathcal{N}_1')=\mathbb{D}'E_{11}'$,
 a contradiction.

{\em Subcase 2.2.} \  $\psi$ does not satisfies the Condition (I).  By Lemma \ref{slynsf32524asd},
$\psi\left(\mathbb{D}_{\leq 1}^{n\times n}\right)$  is  an adjacent set. Thus $\varphi\left(\mathbb{B}_A\right)$  is  an adjacent set.
Then, there is a  maximal set $\mathcal{M}'$ containing $\varphi(A)$ in ${\mathbb{D}'}^{m'\times n'}$, such that
\begin{equation}\label{6gvnao785wew}
\varphi\left(\mathbb{B}_A\right)\subseteq \mathcal{M}'.
\end{equation}
Without loss of generality, we assume that $\mathcal{M}'$ is of type one.

Since  ${\rm rank}(A)\leq 1$,
there is a type one  maximal set $\mathcal{M}_A$ [resp. type two maximal set $\mathcal{S}_A$] containing $A$ and $0$ in $\mathbb{D}^{n\times n}$.
Since $\varphi$ satisfies the Condition (I), there are two maximal sets ${\cal M}_A'$ and ${\cal S}_A'$ of different types
in ${\mathbb{D}'}^{m'\times n'}$, such that $\varphi({\cal M}_A)\subseteq {\cal M}'_A$ and $\varphi({\cal S}_A)\subseteq {\cal S}'_A$.
Without loss of generality, we assume that $\mathcal{M}'_A$ is of type one and $\mathcal{S}'_A$ is of type two.
Since   $\mathcal{S}_A\subseteq \mathbb{D}_{\leq 1}^{n\times n}+A$,
it follows from (\ref{6gvnao785wew}) that $\varphi(\mathcal{S}_A)\subseteq ({\cal S}'_A\cap {\cal M}')$.
By Lemma \ref{Rectangular-PID2-11}, $\ell:={\cal S}'_A\cap {\cal M}'$ is a line in $AG({\cal S}'_A)$,  a contradiction to Lemma \ref{Recastatre43}.

Therefore, there is always a contradiction in  Case 2.

{\bf Case 3.} \ $d(A_0,A)=n-1$. Then  $A<A_0$. Since  ${\rm rank}(\varphi(A_0))={\rm rank}(A_0)=n$ and
 Corollary \ref{mYR25HLLsad23},  we get $\varphi(A)<\varphi(A_0)$.
By Lemma \ref{minuspar1}(c), there are invertible matrices $P_1,Q_1$ over $\mathbb{D}$ such that $A=P_1E_{11}Q_1$ and $A_0=P_1I_nQ_1$,
and there are invertible matrices $P_2,Q_2$ over $\mathbb{D}'$ such that $\varphi(A)=P_2E_{11}'Q_2$ and $\varphi(A_0)=P_2{\rm diag}(I_n,0)Q_2$.
Replacing  $\varphi$ by the map $X\mapsto P_2^{-1}\varphi(P_1XQ_1)Q_2^{-1}$,
we have that $A=E_{11}$,  $A_0=I_n$,  $\varphi(E_{11})=E_{11}'$  and  $\varphi(I_n)={\rm diag}(I_n,0)$.

Since $\varphi$ satisfies the Condition (I) and $\varphi(E_{11})=E_{11}'$,
 Corollary \ref{Rectangular-PID2-7} implies either $\varphi(\mathcal{M}_1)\subseteq \mathcal{M}_1'$ with $\varphi(\mathcal{N}_1)\subseteq \mathcal{N}_1'$, or
 $\varphi(\mathcal{M}_1)\subseteq \mathcal{N}_1'$ with $\varphi(\mathcal{N}_1)\subseteq \mathcal{M}_1'$. We  prove this result only for the first case;
the second case is similar. From now on we assume that
$$
\varphi(\mathcal{M}_1)\subseteq \mathcal{M}_1' \ \ {\rm and} \ \ \varphi(\mathcal{N}_1)\subseteq \mathcal{N}_1'.
$$
Recalling (\ref{ghf54sbxbnpa12}),  we have $\varphi\left(\mathbb{D}_{\leq 1}^{n\times n}+E_{11}\right)\subseteq \mathcal{R}\cup \mathcal{S}$  with
$\varphi(E_{11})=E_{11}'\in\mathcal{R}\cap\mathcal{S}$.
 Since $\mathcal{M}_1, \mathcal{N}_1\subset (\mathbb{D}_{\leq 1}^{n\times n}+E_{11})$, we obtain that
 $\varphi(\mathcal{M}_1)\subseteq \mathcal{M}_1'\cap(\mathcal{R}\cup \mathcal{S})=(\mathcal{M}_1'\cap\mathcal{R})\cup(\mathcal{M}_1'\cap\mathcal{S})$ and
$\varphi(\mathcal{N}_1)\subseteq \mathcal{N}_1'\cap(\mathcal{R}\cup \mathcal{S})=(\mathcal{N}_1'\cap\mathcal{R})\cup(\mathcal{N}_1'\cap\mathcal{S})$.

We show $\mathcal{M}_1'=\mathcal{R}$ by contradiction.
Suppose that $\mathcal{M}_1'\neq\mathcal{R}$. Then $\mathcal{M}_1'\cap\mathcal{R}=\{E_{11}'\}$ and
 hence $\varphi(\mathcal{M}_1)\subseteq (\mathcal{M}_1'\cap \mathcal{S})$.
By Lemma \ref{Rectangular-PID2-11}, $\ell:=\mathcal{M}_1'\cap \mathcal{S}$ is a line in $AG({\cal M}'_1)$,  a contradiction to Lemma \ref{Recastatre43}.
Hence $\mathcal{M}_1'=\mathcal{R}$. Similarly,  $\mathcal{N}_1'=\mathcal{S}$. Thus
$$
\mbox{$\varphi\left(\mathbb{D}_{\leq 1}^{n\times n}+E_{11}\right)\subseteq \mathcal{M}_1'\cup \mathcal{N}_1'$.}
$$
Since $E_{11}+E_{22}\in (\mathbb{D}_{\leq 1}^{n\times n}+E_{11})$, we get $\varphi(E_{11}+E_{22})\in \mathcal{M}_1'$ or $\varphi(E_{11}+E_{22})\in \mathcal{N}_1'$.
Without loss of generality, we assume that $\varphi(E_{11}+E_{22})\in \mathcal{M}_1'$. Recall $\varphi(I_n)={\rm diag}(I_n,0)$.
We have $d(E_{11}+E_{22}, I_n)=n-2$ and $d\left(\varphi(E_{11}+E_{22}), {\rm diag}(I_n,0)\right)\geq n-1$,  a contradiction.

Combining Cases 1-3, we always have a contradiction. Hence $\varphi$ is non-degenerate.
$\qed$

Now, we prove Theorem \ref{non-degenerate-b} as follows.

\noindent{\bf Proof of Theorem \ref{non-degenerate-b}.} \
Suppose  the graph homomorphism $\varphi$ satisfies the Condition (I). Then, by Lemmas \ref{slynon-degenerate-2} and \ref{slycxzwe46e4dfz},
$\varphi$ is non-degenerate  and ${\rm dim}(\varphi(\mathcal{M}_1))={\rm dim}(\varphi(\mathcal{N}_1))=n$.
By  \cite[Theorem 3.1]{Huangli-I}, $\varphi$ is of the form either (\ref{3654665cnc00}) or (\ref{CX3mmbb6600}).
Now, we assume that the homomorphism $\varphi$ does not satisfies the Condition (I). Then Lemma \ref{slynsf32524asd} implies that
$\varphi\left(\mathbb{D}_{\leq 1}^{n\times n}\right)$  is  an adjacent set.
$\qed$

\section{Degenerate graph homomorphisms}

\ \ \ \ \ \ In this section, we  discuss the degenerate graph homomorphisms.
For a degenerate graph homomorphism $\varphi:{\mathbb{D}}^{m\times n}\rightarrow  {\mathbb{D}'}^{m'\times n'}$,
there are no explicit algebraic formulas of $\varphi$. Thus, we discuss the ranges and some properties
on degenerate  graph homomorphisms.

Let $\varphi:{\mathbb{D}}^{m\times n}\rightarrow  {\mathbb{D}'}^{m'\times n'}$ be a graph homomorphism. We called that
{\em $\varphi$ maps  distinct maximal sets  of the same type $[$resp. different types$]$  into  distinct maximal sets of the same type $[$resp. different types$]$},
if for any two distinct maximal sets $\mathcal{M}$ and $\mathcal{N}$ of the same type [resp. different types] in $\mathbb{D}^{m\times n}$, there are
two distinct maximal sets $\mathcal{M}'$ and $\mathcal{N}'$  of the same type [resp. different types] in ${\mathbb{D}'}^{m'\times n'}$,
such that $\varphi(\mathcal{M})\subseteq \mathcal{M}'$ and $\varphi(\mathcal{N})\subseteq \mathcal{N}'$.
We called that {\em $\varphi$ maps  maximal sets  of the same type $[$resp. different types$]$ into maximal sets of the same type $[$resp. different types$]$},
if for any two distinct maximal sets $\mathcal{M}$ and $\mathcal{N}$ of the same type $[$resp. different types$]$ in $\mathbb{D}^{m\times n}$, there are
two  maximal sets $\mathcal{M}'$ and $\mathcal{N}'$  of the same type [resp. different types] in ${\mathbb{D}'}^{m'\times n'}$, such that
$\varphi(\mathcal{M})\subseteq \mathcal{M}'$ and $\varphi(\mathcal{N})\subseteq \mathcal{N}'$
(where $\mathcal{N}'$ and $\mathcal{N}'$  may be equal if $\mathcal{M}$ and $\mathcal{N}$ are of the same type).

 Our main results in this section are the following two theorems.

\begin{thm}\label{distinguishing-non01} \ Let $\mathbb{D}, \mathbb{D}'$ be division rings  with $|\mathbb{D}|\geq 4$,
 and let $m,n,m', n'\geq 2$ be integers with $m',n'\geq {\rm min}\{m,n\}$. Suppose  $\varphi:{\mathbb{D}}^{m\times n}\rightarrow  {\mathbb{D}'}^{m'\times n'}$
is a degenerate graph homomorphism with $\varphi(0)=0$,  and there exists $A_0\in {\mathbb{D}}^{m\times n}$ such that ${\rm rank}(\varphi(A_0))={\rm min}\{m,n\}$.
Then  $\varphi\left(\mathbb{D}_{\leq 1}^{m\times n} \right)$ and $\varphi(\mathbb{B}_{A_0})$ are two adjacent sets.
\end{thm}

\begin{thm}\label{degenerate7654}  Let $\mathbb{D}, \mathbb{D}'$ be division rings with $|\mathbb{D}|\geq 4$, and let $m,n,m', n'\geq 2$ be integers
with ${\rm min}\{m,n\}=2$.
Assume that $\varphi: \mathbb{D}^{m\times n}\rightarrow  {\mathbb{D}'}^{m'\times n'}$
is a degenerate graph homomorphism. Then there are two fixed maximal sets $\mathcal{M}$ and $\mathcal{N}$ of different types containing $0$ in ${\mathbb{D}'}^{m'\times n'}$,
 such that
\begin{equation}\label{432fgmv8772m}
\varphi(\mathbb{D}^{m\times n})\subseteq (\mathcal{M}+R)\cup(\mathcal{N}+R),
\end{equation}
where $R\in {\mathbb{D}'}^{m'\times n'}$ is fixed.
\end{thm}

To prove Theorems \ref{distinguishing-non01} and \ref{degenerate7654}, we need  the following lemmas.

\begin{lem}\label{stronglydegenerate3} Let $\mathbb{D}, \mathbb{D}'$ be division rings  with $|\mathbb{D}|\geq 4$,
 and let $m,n,m', n'\geq 2$ be integers with $n, m', n'\geq m$. Suppose that
$\varphi:{\mathbb{D}}^{m\times n}\rightarrow  {\mathbb{D}'}^{m'\times n'}$
is a graph homomorphism with $\varphi(0)=0$,  and there exists $A_0\in {\mathbb{D}}^{m\times n}$ such that ${\rm rank}(\varphi(A_0))=m$.
Assume further that $\varphi\left(\mathbb{D}_{\leq 1}^{m\times n} \right)$ is not an adjacent set. Then:
\begin{itemize}
\item[{\rm (i)}] $\varphi$ maps distinct maximal sets of type one containing $0$ into distinct maximal sets of the same type, and
$\varphi$ maps  maximal sets of different types  [resp. the same type]  containing $0$ into  maximal sets of different types  [resp. the same type].
Moreover, if ${\cal M}$ is a maximal set containing $0$ in ${\mathbb{D}}^{m\times n}$
and $\varphi(\mathcal{M})\subseteq \mathcal{M}'$ where ${\cal M}'$ is a maximal set in ${\mathbb{D}'}^{m'\times n'}$,
then $\varphi({\cal M})$ is not contained in any line in $AG({\cal M}')$;

\item[{\rm (ii)}] if $A\in\mathbb{D}_1^{m\times n}$ and $d(A_0,A)=m-1$, then $\varphi\left(\mathbb{B}_A\right)$
 is not contained in a union of two maximal sets of different types containing $\varphi(A)$.
\end{itemize}
 \end{lem}
\proof
 By Theorem \ref{non-degenerate-b} or Corollary \ref{stronglydegenerate2}, without loss of generality, we  assume that $n>m$. There is $Q_1\in GL_n(\mathbb{D})$  such that
$A_0=(A_1, 0)Q_1$ where $A_1\in GL_m(\mathbb{D})$. Let $\psi(X)=\varphi(XQ_1)$  for all $X\in {\mathbb{D}}^{m\times n}$.
Then $\psi$ is a graph homomorphism with $\psi(0)=0$, and ${\rm rank}(\psi(A_1, 0))=m={\rm rank}(A_1,0)$.
Since  $\varphi\left(\mathbb{D}_{\leq 1}^{m\times n} \right)$ is not an adjacent set,
$\psi\left(\mathbb{D}_{\leq 1}^{m\times n} \right)$ is not an adjacent set.
Definite the map $g:{\mathbb{D}}^{m\times m}\rightarrow  {\mathbb{D}'}^{m'\times n'}$ by
$$g(Y)=\psi((Y, 0))=\varphi((Y, 0)Q_1), \ \ Y\in {\mathbb{D}}^{m\times m}.$$
 Then $g$ is a graph homomorphism with $g(0)=0$,
${\rm rank}(g(A_1))=m={\rm rank}(A_1)$. By Theorem \ref{non-degenerate-b} or Corollary \ref{stronglydegenerate2},  either
$g\left(\mathbb{D}_{\leq 1}^{m\times m} \right)$ is an adjacent set, or $g$ is a distance preserving map.

We assert that $g\left(\mathbb{D}_{\leq 1}^{m\times m} \right)$ is not  an adjacent set. Otherwise,  $g(\mathbb{D}_{\leq 1}^{m\times m})\subseteq \mathcal{R}'$
where $\mathcal{R}'$ is a maximal set  containing $0$ in $\mathbb{D'}^{m'\times n'}$. We show  a contradiction as follows.

By Lemma \ref{Rectangular-PID2-4}, $\mathcal{R}'=P_2\mathcal{M}'_1Q_2$ or $\mathcal{R}'=P_2\mathcal{N}'_1Q_2$
where $P_2\in GL_{m'}(\mathbb{D}')$ and $Q_2\in GL_{n'}(\mathbb{D}')$. Replacing the map $g$ by the map $X\mapsto P_2^{-1}g(X)Q_2^{-1}$,
we have  either $\mathcal{R}'=\mathcal{M}'_1$ or $\mathcal{R}'=\mathcal{N}'_1$.
Without loss of generality, we may assume that $\mathcal{R}'=\mathcal{M}_1'$. Then
\begin{equation}\label{VC2HFJ9780AS44}
g(\mathbb{D}_{\leq 1}^{m\times m})\subseteq \mathcal{M}_1'.
\end{equation}
Thus, we obtain
\begin{equation}\label{cxvn24kdd}
\psi(\mathcal{N}_j)\subseteq \mathcal{M}_1', \ \ j=1,\ldots, m.
\end{equation}
Set $\left(\mathbb{D}_{\leq 1}^{m\times m}, 0 \right)= \left\{ (X,0): X\in \mathbb{D}_{\leq 1}^{m\times m} \right\}\subset \mathbb{D}_{\leq 1}^{m\times n}$.
By (\ref{VC2HFJ9780AS44}), we have $\psi\left(\mathbb{D}_{\leq 1}^{m\times m}, 0 \right)\subseteq \mathcal{M}_1'$.
Let $\mathcal{S}$ be any  maximal set of type one containing $0$ in $\mathbb{D}^{m\times n}$, and let $\psi(\mathcal{S})\subseteq {\cal S}'$
where ${\cal S}'$ is a  maximal set in ${\mathbb{D}'}^{m'\times n'}$.
 Then $|\mathcal{S}\cap\mathcal{N}_1|\geq 2$. Thus (\ref{cxvn24kdd}) implies that  $|\mathcal{S}'\cap\mathcal{M}_1'|\geq 2$. Using  Corollary \ref{Rectangular-1-13},
 either ${\cal S}'={\cal M}'_1$ or ${\cal S}'$ is of type two.
Put $(\mathcal{S}^{(m)},0)=\mathcal{S}\cap\left(\mathbb{D}_{\leq 1}^{m\times m}, 0 \right)$. Then $\mathcal{S}^{(m)}$ is a  maximal set of type one in $\mathbb{D}^{m\times m}$.
From  (\ref{VC2HFJ9780AS44})
we get $g(\mathcal{S}^{(m)})\subseteq \mathcal{M}_1'$. Since $g(\mathcal{S}^{(m)})\subseteq \mathcal{S}'$,
$g(\mathcal{S}^{(m)})\subseteq \mathcal{M}_1'\cap {\cal S}'$. If ${\cal S}'$ is of type two, then
$\ell:=\mathcal{M}_1'\cap {\cal S}'$ is a line in $AG(\mathcal{M}_1')$, which is a contradiction because Lemma \ref{Recastatre43}. Thus, ${\cal S}'={\cal M}'_1$.
Then, $\psi$ maps every   maximal set of type one containing $0$ in $\mathbb{D}^{m\times n}$ into $\mathcal{M}_1'$.
Since every matrix of rank one is contained in a maximal set of type one containing $0$, we obtain
$\psi(\mathbb{D}_{\leq 1}^{m\times n})\subseteq \mathcal{M}_1'$. Since $\psi\left(\mathbb{D}_{\leq 1}^{m\times n} \right)$ is not an adjacent set,
we get a contradiction.

Therefore, $g$  is a distance preserving map, and hence $g$ is non-degenerate.

(i). \ By  Lemma \ref{non-degeneratelemma00a}(b),
$g$ maps two distinct maximal sets of different types  [resp. the same type] containing $0$ into two distinct maximal sets of different types  [resp. the same type].
Since $g(Y)=\varphi((Y, 0)Q_1)$ ($Y\in {\mathbb{D}}^{m\times m}$), it is easy to see that $\varphi$ maps two distinct maximal sets of type one
containing $0$ into two distinct maximal sets of the same type.

Let ${\cal M}$ be any type one maximal set containing $0$ in ${\mathbb{D}}^{m\times n}$,  and let $\varphi(\mathcal{M})\subseteq \mathcal{M}'$
where ${\cal M}'$ is a maximal set containing $0$ in ${\mathbb{D}'}^{m'\times n'}$.
 Set $(\mathcal{M}^{(m)},0)=\mathcal{M}\cap\left(\mathbb{D}_{\leq 1}^{m\times m}, 0 \right)$.
Then $\mathcal{M}^{(m)}$ is a  maximal set of type one containing $0$ in $\mathbb{D}^{m\times m}$.
 Recall that $g(Y)=\varphi((Y, 0)Q_1)$ ($Y\in {\mathbb{D}}^{m\times m}$).
We get
$(\mathcal{M}^{(m)},0)Q_1\subseteq \mathcal{M}Q_1=\mathcal{M}$. Thus,
$$\varphi((\mathcal{M}^{(m)},0)Q_1)=g(\mathcal{M}^{(m)})\subseteq\varphi(\mathcal{M})\subseteq \mathcal{M}'.$$
By Lemma \ref{non-degeneratelemma00a}(c), $g(\mathcal{M}^{(m)})$ is not contained in any line in $AG({\cal M}')$.
It follows that  $\varphi(\mathcal{M})$ is not contained in any line in $AG({\cal M}')$.

 Thus,  if ${\cal M}$ is a type one maximal set containing $0$ in ${\mathbb{D}}^{m\times n}$
and $\varphi(\mathcal{M})\subseteq \mathcal{M}'$ where ${\cal M}'$ is a maximal set in ${\mathbb{D}'}^{m'\times n'}$,
then $\varphi({\cal M})$ is not contained in any line in $AG({\cal M}')$.

Next, we discuss maximal sets of type two.
Let ${\cal S}_1$, ${\cal S}_2$ be any two distinct maximal sets of type two containing $0$ in ${\mathbb{D}}^{m\times n}$, and let $\varphi(\mathcal{S}_i)\subseteq \mathcal{S}_i'$
where ${\cal S}_i'$ is a maximal set in ${\mathbb{D}'}^{m'\times n'}$, $i=1,2$. Assume that $\mathcal{R}_1, \mathcal{R}_2$ are any two
distinct type one maximal sets containing $0$ in ${\mathbb{D}}^{m\times n}$. By the above results, there are two distinct  maximal sets $\mathcal{R}_1', \mathcal{R}_2'$
containing $0$ in ${\mathbb{D}'}^{m'\times n'}$, such that $\mathcal{R}_1'$ and $\mathcal{R}_2'$ are of the same type and
$\varphi(\mathcal{R}_i)\subseteq \mathcal{R}_i'$, $i=1,2$. By Corollary \ref{Rectangular-1-13}, we have $|\mathcal{S}_i'\cap\mathcal{R}_j'|\geq 2$, $i,j=1,2$.
Thus $\mathcal{R}_1'\neq \mathcal{R}_2'$ implies that $\mathcal{S}_i'$ and $\mathcal{R}_1'$ are of different types $i=1,2$.
Consequently  ${\cal S}_1'$ and ${\cal S}_2'$  are of the same type.

Therefore, we have proved that
$\varphi$ maps two distinct maximal sets of different types [resp. the same type] containing $0$ into two  maximal sets of different types [resp. the same type].

On the other hand, by Lemma \ref{Rectangular-PID2-11},
$\ell_1:=\mathcal{S}_1\cap\mathcal{R}_1$  and $\ell_2:=\mathcal{S}_1\cap\mathcal{R}_2$ are two distinct lines in $AG(\mathcal{S}_1)$. Also,
 $\ell_1':=\mathcal{S}_1'\cap\mathcal{R}_1'$  and $\ell_2':=\mathcal{S}_2'\cap\mathcal{R}_2'$ are two distinct lines in $AG(\mathcal{S}_1')$.
 Since $\varphi(\ell_i)\subseteq \ell_i'$, $i=1,2$, it follows that $\varphi({\cal S}_1)$ is not contained in any line in $AG({\cal S}_1')$ because
 two different lines have at most a common point.
 Thus,  if ${\cal S}_1$ is a type two maximal set containing $0$ in ${\mathbb{D}}^{m\times n}$
and $\varphi(\mathcal{S}_1)\subseteq \mathcal{S}_1'$ where ${\cal S}_1'$ is a maximal set in ${\mathbb{D}'}^{m'\times n'}$,
then $\varphi({\cal S}_1)$ is not contained in any line in $AG({\cal S}_1')$.

Therefore,  if ${\cal M}$ is a maximal set containing $0$ in ${\mathbb{D}}^{m\times n}$ and $\varphi(\mathcal{M})\subseteq \mathcal{M}'$ where ${\cal M}'$ is a maximal set in ${\mathbb{D}'}^{m'\times n'}$,
we have prove that $\varphi({\cal M})$ is not contained in any line in $AG({\cal M}')$. Thus, the (i) of this lemma is proved.

(ii). \  Suppose $A\in\mathbb{D}_1^{m\times n}$ and $d(A_0,A)=m-1$. Then $A<A_0$.
By Lemma \ref{minuspar1}(c), there are $P_3\in GL_m(\mathbb{D})$ and $Q_3\in GL_n(\mathbb{D})$ such that
$A=P_3E_{11}Q_3$ and $A_0=P_3(I_{m}, 0)Q_3$. Replacing $\varphi$ by the map $X\mapsto \varphi(P_3XQ_3)$,
we can assume that $A_0=(I_m, 0)$ and $A=E_{11}$.

Definite the map $f:{\mathbb{D}}^{m\times m}\rightarrow  {\mathbb{D}'}^{m'\times n'}$ by
$f(Y)=\varphi(Y, 0)$, $Y\in {\mathbb{D}}^{m\times m}$.
 Then $f$ is a graph homomorphism such that $f(0)=0$ and ${\rm rank}(f(I_m))=m$. By Theorem \ref{non-degenerate-b} or Corollary \ref{stronglydegenerate2},  either
$f\left(\mathbb{D}_{\leq 1}^{m\times m} \right)$ is an adjacent set or $f$ is a distance preserving map.
Similar to the proof on the  homomorphism $g$, we can prove that $f$ is a distance preserving map. Since every distance preserving map carries
 distinct maximal sets into  distinct maximal sets, $f\left(\mathbb{D}_{\leq 1}^{m\times m}+E_{11}^{m\times m}\right)$
 is not contained in a union of two maximal sets of different types containing $f(E_{11}^{m\times m})$.
Since $f\left(\mathbb{D}_{\leq 1}^{m\times m}+E_{11}^{m\times m}\right)=
\varphi\left(\mathbb{D}_{\leq 1}^{m\times m}+E_{11}^{m\times m}, \, 0\right)\subseteq \varphi\left(\mathbb{B}_{E_{11}}\right),$
 $\varphi\left(\mathbb{B}_{E_{11}}\right)$
 is not contained in a union of two maximal sets of different types containing $\varphi(E_{11})$. Thus the (ii) of this lemma is proved.
$\qed$

By the symmetry of rows and columns of a matrix, we have similarly the following lemma.

\begin{lem}\label{s546hma97e3} Let $\mathbb{D}, \mathbb{D}'$ be division rings  with $|\mathbb{D}|\geq 4$,
 and let $m,n,m', n'\geq 2$ be integers with  $m, m', n'\geq n$. Suppose that
$\varphi:{\mathbb{D}}^{m\times n}\rightarrow  {\mathbb{D}'}^{m'\times n'}$
is a graph homomorphism with $\varphi(0)=0$,  and there exists $A_0\in {\mathbb{D}}^{m\times n}$ such that ${\rm rank}(\varphi(A_0))=n$.
Assume further that $\varphi\left(\mathbb{D}_{\leq 1}^{m\times n} \right)$ is not an adjacent set. Then:
\begin{itemize}

\item[{\rm (i)}]  $\varphi$ maps distinct type two maximal sets  containing $0$ into distinct maximal sets of the same type, and
$\varphi$ maps  maximal sets of different types   [resp. the same type] containing $0$ into  maximal sets of different types  [resp. the same type].
Moreover, if ${\cal M}$ is a maximal set containing $0$ in ${\mathbb{D}}^{m\times n}$ and $\varphi(\mathcal{M})\subseteq \mathcal{M}'$ where
${\cal M}'$ is a maximal set in ${\mathbb{D}'}^{m'\times n'}$, then $\varphi({\cal M})$ is not contained in any line in $AG({\cal M}')$;

\item[{\rm (ii)}] if $A\in\mathbb{D}_1^{m\times n}$ and $d(A_0,A)=n-1$, then $\varphi\left(\mathbb{B}_A\right)$
 is not contained in a union of two maximal sets of different types containing $\varphi(A)$.
\end{itemize}
 \end{lem}

Now, we prove Theorem \ref{distinguishing-non01} as follows.

\noindent{\bf Proof of Theorem \ref{distinguishing-non01}}. \
We prove this theorem only for the case of $m={\rm min}\{m,n\}$; the case of $n={\rm min}\{m,n\}$ is similar by using Lemma \ref{s546hma97e3}.
Since $\varphi$ is degenerate, there exists a matrix $A\in\mathbb{D}^{m\times n}_{\leq 1}$ and there are two
maximal sets $\mathcal{M}$ and $\mathcal{N}$ of different types in ${\mathbb{D}'}^{m'\times n'}$, such that
\begin{equation}\label{gd4nfdkkl8qw}
\mbox{$\varphi\left(\mathbb{B}_A\right)\subseteq \mathcal{M}\cup \mathcal{N}$ \
 with \ $\varphi(A)\in\mathcal{M}\cap\mathcal{N}$.}
\end{equation}
Without loss of generality, we assume that $\mathcal{M}$ is of type one  and $\mathcal{N}$ is of type two.
Since ${\rm rank}(\varphi(A_0))=m$ and ${\rm rank}(\varphi(A))\leq 1$,  we have either  $d(\varphi(A_0),\varphi(A))=m$ or $m-1$.

{\bf Step 1.} \ In this step, we will prove that $\varphi\left(\mathbb{D}_{\leq 1}^{m\times n} \right)$ is an adjacent set by contradiction.
By Corollary \ref{stronglydegenerate2}, without loss of generality, we  assume that $n>m$.
Suppose that  $\varphi\left(\mathbb{D}_{\leq 1}^{m\times n} \right)$ is not an adjacent set.
We distinguish the following two cases to show  a contradiction.

{\bf Case 1.1}. \ $d(\varphi(A_0),\varphi(A))=m$.
Let $\psi(X)=\varphi(X+A)-\varphi(A)$, $X\in{\mathbb{D}}^{m\times n}$. Then $\psi:  {\mathbb{D}}^{m\times n}\rightarrow  {\mathbb{D}'}^{m'\times n'}$
is a graph homomorphism with $\psi(0)=0$  and  ${\rm rank}(\psi(A_0-A))=m$. We show that $\psi\left(\mathbb{D}_{\leq 1}^{m\times n} \right)$ is  an adjacent set.
Otherwise, if $\psi\left(\mathbb{D}_{\leq 1}^{m\times n} \right)$ is not an adjacent set, then by Lemma \ref{stronglydegenerate3}(i),
$\psi\left(\mathbb{D}_{\leq 1}^{m\times n} \right)$ is not contained in a union of two maximal sets of different types containing $0$,
which implies that $\varphi\left(\mathbb{B}_A\right)$ is not contained in a union of two maximal sets of different types containing $\varphi(A)$,
a contradiction to (\ref{gd4nfdkkl8qw}). Therefore,  $\psi\left(\mathbb{D}_{\leq 1}^{m\times n}\right)$ is an adjacent set.
It follows that  $\varphi\left(\mathbb{B}_A\right)$ is an adjacent set. Thus, there exists a maximal set $\mathcal{S}$ such that
$$\varphi\left(\mathbb{B}_A\right)\subseteq \mathcal{S}.$$

Since $\varphi\left(\mathbb{D}_{\leq 1}^{m\times n} \right)$ is not an adjacent set,  $A$ is of rank one.
There is a type one  maximal set $\mathcal{M}_A$ [resp. type two maximal set $\mathcal{S}_A$] containing $A$ and $0$ in $\mathbb{D}^{m\times n}$.
 By Lemma \ref{stronglydegenerate3}(i), we have  $\varphi({\cal M}_A)\subseteq {\cal M}'_A$ and $\varphi({\cal S}_A)\subseteq {\cal S}'_A$, where
 ${\cal M}_A'$ and ${\cal S}_A'$ are two maximal sets  of different types in ${\mathbb{D}'}^{m'\times n'}$ containing $\varphi(A)$ and $0$.
Since   $\mathcal{M}_A, \mathcal{S}_A\subseteq \mathbb{B}_A$, we get  $\varphi(\mathcal{M}_A)\subseteq ({\cal M}'_A\cap {\cal S})$ and $\varphi(\mathcal{S}_A)\subseteq ({\cal S}'_A\cap {\cal S})$.
Without loss of generality, we assume that $\mathcal{M}'_A$ and  ${\cal S}$ are of the same type. Then Corollary \ref{Rectangular-1-13}
implies that ${\cal M}'_A\cap {\cal S}=\{\varphi(A)\}$ and $\varphi(\mathcal{M}_A)\subseteq \{\varphi(A)\}$, a contradiction.

{\bf Case 1.2}. \ $d(\varphi(A_0),\varphi(A))=m-1$.
 Then, we have ${\rm rank}(\varphi(A))={\rm rank}(A)=1$ and $\varphi(A)<\varphi(A_0)$. Write $E_{ij}'=E_{ij}^{m'\times n'}$ and $E_{ij}=E_{ij}^{m\times n}$.
By Lemma \ref{minuspar1}(c), there are invertible matrices $P_1,Q_1$ over $\mathbb{D}'$ such that
$$\mbox{$\varphi(A)=P_1E_{11}'Q_1$ \ and \ $\varphi(A_0)=P_1{\rm diag}(I_m,0)Q_1$.}$$
Also, there are invertible matrices $P_2,Q_2$ over $\mathbb{D}$ such that $P_2AQ_2=E_{11}$.
Let
$\scriptsize P_2A_0Q_2=\left(
                 \begin{array}{cc}
                   A_{11} & A_{12}\\
                   A_{21} & A_{22}\\
                 \end{array}
               \right)$ where $A_{11}\in {\mathbb{D}}^{1\times 1}$ and $A_{22}\in {\mathbb{D}}^{(m-1)\times (n-1)}$.
 Clearly, $d(A_0,A)=m$ or $m-1$. By Lemma \ref{stronglydegenerate3}(ii) and (\ref{gd4nfdkkl8qw}), we have $d(A_0,A)=m$.
 Thus, either $A_{22}\cong {\rm diag}(I_{m-1}, 0)$ or $A_{22}\cong {\rm diag}(I_{m-2}, 0)$. Put
 $\scriptsize J_m=\left(
   \begin{array}{ccc}
     0 & 0 & 1 \\
     0 & I_{m-2} & 0 \\
     1 & 0 & 0 \\
   \end{array}
 \right)$.

When $A_{22}\cong {\rm diag}(I_{m-1}, 0)$, by appropriate elementary row and column operations on matrices, we can obtain
either  $\left({\rm diag}(a, I_{m-1}), 0\right)$ (where  $a\in \mathbb{D}^*$ with $a\neq 1$)  or $\scriptsize\left(
                                               \begin{array}{cccc}
                                                 0_1 & 0 & 1&0 \\
                                                 0 & I_{m-1} & 0&0 \\
                                               \end{array}
                                             \right)$
from  $P_2A_0Q_2$.  Moreover,  $P_2AQ_2=E_{11}$ is unchanged.

When $A_{22}\cong {\rm diag}(I_{m-2}, 0)$, by appropriate elementary row and column operations on matrices, we can get
$(J_m, 0)$ from $P_2A_0Q_2$. Moreover, $P_2AQ_2=E_{11}$ is unchanged.

Thus, there are invertible matrices $P_3,Q_3$ over $\mathbb{D}$, such that $P_3AQ_3=E_{11}$ and either
$$\mbox{$P_3A_0Q_3= \left({\rm diag}(a, I_{m-1}),0\right)$, or $(J_m,0)$,
  or $\left(
                                               \begin{array}{cccc}
                                                 0_1 & 0 & 1&0 \\
                                                 0 & I_{m-1} & 0&0 \\
                                               \end{array}
                                             \right)$.}$$
Replacing  $\varphi$ by the map $X\mapsto P_1^{-1}\varphi(P_3^{-1}XQ_3^{-1})Q_1^{-1}$, we have  $A=E_{11}$, and either
$A_0=\left({\rm diag}(a, I_{m-1}),0\right)$,   or $(J_m,0)$, or $\scriptsize\left(
                                               \begin{array}{cccc}
                                                 0_1 & 0 & 1&0 \\
                                                 0 & I_{m-1} & 0&0 \\
                                               \end{array}
                                             \right)$.
Moreover, we have
\begin{equation}\label{BVCEWRW3UTITOA967H}
\mbox{$\varphi(E_{11})=E_{11}'$ \ and \ $\varphi(A_0)={\rm diag}(I_m,0)$.}
\end{equation}
On the other hand, (\ref{gd4nfdkkl8qw}) becomes
$$\mbox{$\varphi\left(\mathbb{B}_{E_{11}}\right)\subseteq \mathcal{M}\cup \mathcal{N}$ \ with \ $\varphi(E_{11})\in\mathcal{M}\cap\mathcal{N}$.}$$

In $\mathbb{D}^{m\times n}$ [resp.  ${\mathbb{D}'}^{m'\times n'}$], by Corollary \ref{Rectangular-PID2-7},  there are only two maximal sets
of different types containing $E_{11}$ and $0$ [resp. $E_{11}'$ and $0$], they are $\mathcal{M}_1$ and $\mathcal{N}_1$ [resp.
$\mathcal{M}_1'$ and $\mathcal{N}_1'$]. Thus, by Lemma \ref{stronglydegenerate3}(i) and $\varphi(E_{11})=E_{11}'$,
we have either $\varphi(\mathcal{M}_1)\subseteq \mathcal{M}_1'$  with $\varphi(\mathcal{N}_1)\subseteq \mathcal{N}_1'$, or
$\varphi(\mathcal{M}_1)\subseteq \mathcal{N}_1'$  with $\varphi(\mathcal{N}_1)\subseteq \mathcal{M}_1'$.  Without loss of generality, we  assume that
\begin{equation}\label{VCTR24tfdgh6700}
\varphi(\mathcal{M}_1)\subseteq \mathcal{M}_1', \ \ \varphi(\mathcal{N}_1)\subseteq \mathcal{N}_1'.
\end{equation}
Since $\mathcal{M}_1, \mathcal{N}_1\subseteq (\mathbb{D}_{\leq 1}^{m\times n}+E_{11})$,
we get that $\varphi(\mathcal{M}_1)\subseteq {\cal M}'_1\cap (\mathcal{M}\cup \mathcal{N})$
and $\varphi(\mathcal{N}_1)\subseteq {\cal N}'_1\cap (\mathcal{M}\cup \mathcal{N})$.
Suppose that ${\cal M}'_1\neq\mathcal{M}$. Then ${\cal M}'_1\cap\mathcal{M}=\{E_{11}'\}$ and
${\cal M}'_1\cap (\mathcal{M}\cup \mathcal{N})={\cal M}'_1\cap \mathcal{N}$, thus
$\ell:={\cal M}'_1\cap \mathcal{N}$ is a line in $AG({\cal M}'_1)$ and $\varphi({\cal M}'_1)\subseteq \ell$.
By  Lemma \ref{stronglydegenerate3}(i), this is a contradiction.
Therefore, we must have  ${\cal M}'_1=\mathcal{M}$. Similarly,  ${\cal N}'_1=\mathcal{N}$. It follows that
\begin{equation}\label{VC2EHFH46fghbv00}
\mbox{$\varphi\left(\mathbb{B}_{E_{11}}\right)\subseteq \mathcal{M}_1'\cup \mathcal{N}_1'$.}
\end{equation}

In order to give a contradiction, we distinguish the following subcases.

{\em Subcase 1.2.1}. \ $A_0=\left({\rm diag}(a, I_{m-1}), 0\right)$ where $a\in \mathbb{D}^*$ with $a\neq 1$. By (\ref{VCTR24tfdgh6700}),
we have $\varphi(aE_{11})=a^*E_{11}'$ where $a^*\in {\mathbb{D}'}^*$.
Since $aE_{11}<A_0$, Corollary \ref{mYR25HLLsad23} and (\ref{BVCEWRW3UTITOA967H}) imply that $\varphi(aE_{11})<\varphi(A_0)={\rm diag}(I_m,0)$. Thus, $a^*=1$ and hence
$\varphi(aE_{11})=E_{11}'$. Since $\varphi(aE_{11})\sim \varphi(E_{11})=E_{11}'$, we have a contradiction.

{\em Subcase 1.2.2}. \ $A_0=(J_m,0)$. By (\ref{VC2EHFH46fghbv00}), we get $\varphi\left(E_{mm}+E_{11}\right)\in\mathcal{M}_1'\cup \mathcal{N}_1'$.
Without loss of generality, we assume that $\varphi\left(E_{mm}+E_{11}\right)\in\mathcal{M}_1'$. Since $d\left(E_{mm}+E_{11}, (J_m,0)\right)=m-1$,
it follows from (\ref{BVCEWRW3UTITOA967H}) that $d\left(\varphi(E_{mm}+E_{11}), {\rm diag}(I_m,0)\right)\leq m-1$. Thus
 $\scriptsize \varphi(E_{mm}+E_{11})=\left(
                            \begin{array}{cc}
                              1 & \alpha \\
                              0 & 0 \\
                            \end{array}
                          \right)$ where $0\neq \alpha\in {\mathbb{D}'}^{n'-1}$ because $\varphi(E_{mm}+E_{11})\sim \varphi(E_{11})=E_{11}'$.
Since $d\left(E_{m1}+E_{11}, (J_m,0)\right)=m-1$, from (\ref{BVCEWRW3UTITOA967H}) we have
 $d\left(\varphi(E_{m1}+E_{11}), {\rm diag}(I_m,0)\right)\leq m-1$. Thus $ \varphi(\mathcal{N}_1)\subseteq \mathcal{N}_1'$ implies that
 $\scriptsize \varphi(E_{m1}+E_{11})=\left(
                            \begin{array}{cc}
                              1 & 0\\
                              \beta & 0 \\
                            \end{array}
                          \right)$ where $0\neq \beta\in \,^{m'-1}{\mathbb{D}'}$ because $\varphi(E_{m1}+E_{11})\sim E_{11}'$.
By $\varphi(E_{mm}+E_{11})\sim \varphi(E_{m1}+E_{11})$, one gets
$\scriptsize\left(
\begin{array}{cc}
1 & \alpha \\
 0 & 0 \\
 \end{array}
\right)\sim \left(
                            \begin{array}{cc}
                              1 & 0\\
                              \beta & 0 \\
                            \end{array}
                          \right)$, a contradiction.

{\em Subcase 1.2.3}. \  $\scriptsize A_0=\left(
\begin{array}{cccc}
0_1 & 0 & 1&0 \\
 0 & I_{m-1} & 0&0 \\
\end{array}
\right)$. Let $B=E_{11}+E_{mm}$.
Since $B\in\mathbb{B}_{E_{11}}$,  from (\ref{VC2EHFH46fghbv00}) we have $\varphi(B)\in \mathcal{M}_1'$ or $\varphi(B)\in \mathcal{N}_1'$.

 Assume that $\varphi(B)\in \mathcal{M}_1'$. Since $d(B,A_0)=m-1$, by (\ref{BVCEWRW3UTITOA967H}) we get $d(\varphi(B), {\rm diag}(I_m,0))\leq m-1$.
On the other hand, we have $\varphi(B)\sim \varphi(E_{11})=E_{11}'$, which implies that
$\scriptsize\varphi(B)=\left(
                            \begin{array}{cc}
                              1 & \alpha\\
                              0 & 0 \\
                            \end{array}
                          \right)$ where $0\neq \alpha\in \,{\mathbb{D}'}^{n'-1}$.
Let $C_1=E_{11}+E_{m1}$.
Since $C_1\in \mathcal{N}_1$ and $\varphi(C_1)\sim \varphi(B)$, it follows from $\varphi(\mathcal{N}_1)\subseteq \mathcal{N}_1'$  that $\varphi(C_1)=c_1E_{11}'$,  where $c_1\in {\mathbb{D}'}^*$.
Let $\mathcal{M}_{C}$ be the type one  maximal set containing $C_1$ and $0$ in $\mathbb{D}^{m\times n}$. Then  $\mathcal{M}_{C}\neq \mathcal{M}_{1}$.
Suppose  $\varphi(\mathcal{M}_{C})\subseteq \mathcal{M}_{C}'$ where $\mathcal{M}_{C}'$ is a maximal set containing $\varphi(C_1)$ and $0$
in ${\mathbb{D}'}^{m'\times n'}$. By  Lemma \ref{stronglydegenerate3}(i) and (\ref{VCTR24tfdgh6700}),  $\mathcal{M}_{C}'$ is of type one and
$\mathcal{M}_{C}'\neq \mathcal{M}_{1}'$. Thus, $\{0, \varphi(C_1)\}\subseteq \mathcal{M}_{C}'\cap\mathcal{M}_{1}'$, which is a contradiction to Corollary \ref{Rectangular-1-13}.

 Assume that $\varphi(B)\in \mathcal{N}_1'$. Similarly, we have
$\scriptsize\varphi(B)=\left(
                            \begin{array}{cc}
                              1 & 0\\
                              \beta & 0 \\
                            \end{array}
                          \right)$ where $0\neq \beta\in \,^{m'-1}{\mathbb{D}'}$. Let $B_1=E_{11}+E_{m1}+E_{mm}-E_{m,m+1}$.
Since $B_1\in\mathbb{B}_{E_{11}}$,  by (\ref{VC2EHFH46fghbv00}) we get $\varphi(B_1)\in \mathcal{M}_1'$ or $\varphi(B_1)\in \mathcal{N}_1'$.
By $d(B_1,A_0)=m-1$ and (\ref{BVCEWRW3UTITOA967H}), we have  $d(\varphi(B_1), {\rm diag}(I_m,0))\leq m-1$.
Thus $\scriptsize\varphi(B_1)=\left(
                            \begin{array}{cc}
                              1 & \alpha_1\\
                              0 & 0 \\
                            \end{array}
                          \right)$ or $\scriptsize\varphi(B_1)=\left(
                            \begin{array}{cc}
                              1 & 0\\
                              \beta_1 & 0 \\
                            \end{array}
                          \right)$, where $0\neq \alpha_1\in \,{\mathbb{D}'}^{n'-1}$ and $0\neq \beta_1\in \,^{m'-1}{\mathbb{D}'}$.
Since $\varphi(B_1)\sim \varphi(B)$, we must have $\scriptsize\varphi(B_1)=\left(
                            \begin{array}{cc}
                              1 & 0\\
                              \beta_1 & 0 \\
                            \end{array}
\right)$. Let $B_2=-E_{1m}+E_{1,m+1}\in \mathcal{M}_1$. Then $B_2\sim B_1$ and $B_2\sim E_{11}$,  hence (\ref{BVCEWRW3UTITOA967H}) and (\ref{VCTR24tfdgh6700})
imply that $\varphi(B_2)=b_2E_{11}'$, where $b_2\in {\mathbb{D}'}^*$ with $b_2\neq 1$. However, we have $d(B_2,A_0)=m-1$ and $d(\varphi(B_2), \varphi(A_0))=
d(b_2E_{11}', {\rm diag}(I_m,0)))= m$, a contradiction.

Combining Case 1.1 with Case 1.2, there is always a contradiction.
Therefore, $\varphi\left(\mathbb{D}_{\leq 1}^{m\times n} \right)$ must be  an adjacent set.

{\bf Step 2.} \ In this step, we will prove that $\varphi(\mathbb{B}_{A_0})$ is an adjacent set. By Step 1,
 $\varphi\left(\mathbb{D}_{\leq 1}^{m\times n} \right)$ is an adjacent set. Thus, there is a  maximal set $\mathcal{M}'$ containing $0$ in ${\mathbb{D}'}^{m'\times n'}$,
such that $\varphi\left(\mathbb{D}_{\leq 1}^{m\times n} \right)\subseteq \mathcal{M}'$.
There is an invertible matrix $Q$ over $\mathbb{D}$ such that $A_0=(A_1,0)Q$ where $A_1\in GL_m(\mathbb{D})$.

Let $\varphi'(X)=\varphi(XQ)$
for all $X\in \mathbb{D}^{m\times n}$. Then $\varphi'$ is a graph homomorphism, ${\rm rank}(\varphi'(A_1,0))=m$ and
\begin{equation}\label{cxef554eer000i}
\varphi'\left(\mathbb{D}_{\leq 1}^{m\times n}\right)\subseteq \mathcal{M}'.
\end{equation}
Let $f(Y)=\varphi'(Y+A_1,0)-\varphi'(A_1,0)$, $Y\in\mathbb{D}^{m\times m}$. Then $f$ is a graph homomorphism  such that $f(0)=0$ and
${\rm rank}(f(-A_1))=m$.
By (\ref{cxef554eer000i}), it is clear that   $f$ does not preserve  distance  $2$. It follows from Corollary \ref{stronglydegenerate2} that
$f\left(\mathbb{D}_{\leq 1}^{m\times m} \right)$ is an adjacent set.
Hence there is a  maximal set $\mathcal{N}'$ containing $0$ in ${\mathbb{D}'}^{m'\times n'}$, such that
$$f\left(\mathbb{D}_{\leq 1}^{m\times m}\right)\subseteq \mathcal{N}'.$$

Let $h(X)=\varphi'(X+(A_1,0))-\varphi'(A_1,0)$, $X\in\mathbb{D}^{m\times n}$. Then $h$ is a graph homomorphism  with $h(0)=0$.
Moreover, $h(Y,0)=f(Y)$ for all  $Y\in\mathbb{D}^{m\times m}$.
We assert $h\left(\mathbb{D}_{\leq 1}^{m\times n}\right)\subseteq \mathcal{N}'$. Otherwise, there exists a type one maximal set $\mathcal{R}$ containing $0$
in ${\mathbb{D}}^{m\times n}$, such that
$h(\mathcal{R})\subseteq \mathcal{S}$ where $\mathcal{S}$ is a maximal set  containing $0$ in ${\mathbb{D}'}^{m'\times n'}$ and $\mathcal{S}\neq \mathcal{N}'$.
Let $\left(\mathcal{R}^{(m)}, 0\right)=\mathcal{R}\cap \left(\mathbb{D}_{\leq 1}^{m\times m}, 0\right)$.
Then $\mathcal{R}^{(m)}$ is a  type one maximal set  containing $0$ in $\mathbb{D}^{m\times m}$, and hence
$f(\mathcal{R}^{(m)})\subseteq \mathcal{N}'\cap \mathcal{S}$. By Corollary \ref{Rectangular-1-13}, $\mathcal{S}$ and $\mathcal{N}'$ must be of different type. Thus
$\ell:=\mathcal{N}'\cap {\cal S}$ is a line in $AG(\mathcal{N}')$ by Lemma \ref{Rectangular-PID2-11}, a contradiction to Lemma \ref{Recastatre43}.
Therefore, we obtain  $h\left(\mathbb{D}_{\leq 1}^{m\times n}\right)\subseteq \mathcal{N}'$, and hence
$\varphi'\left(\mathbb{B}_{(A_1,0)}\right)\subseteq \mathcal{N}'+\varphi'(A_1,0)$. Consequently,
$\varphi(\mathbb{B}_{A_0})\subseteq \mathcal{N}'+\varphi(A_0)$ and hence  $\varphi(\mathbb{B}_{A_0})$ is an adjacent set.
$\qed$

Next, we will prove Theorem \ref{degenerate7654}. We need the following lemma.

\begin{lem}\label{degenerate675235utigh}  Let $\mathbb{D}, \mathbb{D}'$ be division rings with $|\mathbb{D}|\geq 4$, and let $m,n,m', n'\geq 2$ be integers
with ${\rm min}\{m,n\}=2$. Suppose  $\varphi: \mathbb{D}^{m\times n}\rightarrow  {\mathbb{D}'}^{m'\times n'}$
is a degenerate graph homomorphism with $\varphi(0)=0$, and
there exists $A\in {\mathbb{D}}_2^{m\times n}$ such that ${\rm rank}(\varphi(A))=2$.
 Then there are two fixed maximal sets $\mathcal{M}$ and $\mathcal{N}$ of different types in ${\mathbb{D}'}^{m'\times n'}$,
 such that $0\in \mathcal{M}\cap\mathcal{N}$ and
\begin{equation}\label{vcadgjr772m}
\mbox{$\varphi(\mathbb{D}^{m\times n})\subseteq \mathcal{M}\cup(\mathcal{N}+\varphi(A))=(\mathcal{M}+R)\cup(\mathcal{N}+R)$,}
\end{equation}
where  $R\in\mathcal{M}$ is fixed and ${\rm rank}(R)=1$.
Moreover, there is $\mathcal{M}'\in \{\mathcal{M}, \mathcal{N}\}$ such that
$\varphi\left(\mathbb{D}_{\leq 1}^{m\times n}\right)\subseteq \mathcal{M}'$ and $\varphi(\mathbb{N}_A)\subseteq \mathcal{M}'\cap(\mathcal{N}'+R)$
where $\{\mathcal{M}', \mathcal{N}'\}=\{\mathcal{M}, \mathcal{N}\}$.
\end{lem}
\proof \
 We prove this lemma only for the case $m={\rm min}\{m,n\}=2$; the case  $n={\rm min}\{m,n\}=2$ is similar.
From now on we assume that $m={\rm min}\{m,n\}=2$.

{\bf Step 1}. \ By Theorem \ref{distinguishing-non01}, $\varphi\left(\mathbb{D}_{\leq 1}^{2\times n}\right)$ and  $\varphi(\mathbb{B}_A)$ are two adjacent sets.
Then, there are  maximal sets $\mathcal{M}$ and  $\mathcal{N}$ containing $0$ in ${\mathbb{D}'}^{m'\times n'}$,
such that  $\varphi\left(\mathbb{D}_{\leq 1}^{2\times n}\right)\subseteq \mathcal{M}$ and $\varphi(\mathbb{B}_A)\subseteq \mathcal{N}+\varphi(A)$.
By Lemma \ref{Rectangular-PID2-4}, $\mathcal{M}=P_1\mathcal{M}'_1Q_1$ or $\mathcal{M}=P_1\mathcal{N}'_1Q_1$
where $P_1\in GL_{m'}(\mathbb{D}')$ and $Q_1\in GL_{n'}(\mathbb{D}')$. Replacing  $\varphi$ by the map $X\mapsto P_1^{-1}\varphi(X)Q_1^{-1}$,
we have  either $\mathcal{M}=\mathcal{M}'_1$ or $\mathcal{M}=\mathcal{N}'_1$.
We prove this lemma only for the case of $\mathcal{M}=\mathcal{M}'_1$;   the case of $\mathcal{M}=\mathcal{N}'_1$ is similar.
Now, we assume that $\mathcal{M}=\mathcal{M}'_1$.

Write
$\scriptsize A=
\left(
  \begin{array}{c}
    \alpha_1 \\
    \alpha_2\\
  \end{array}
\right)$ where  $\alpha_1, \alpha_2\in \mathbb{D}^n$ are left linearly independent. Set
 $\scriptsize C_\lambda=
\left(
  \begin{array}{c}
    -\alpha_1+\lambda\alpha_2 \\
    0\\
  \end{array}
\right)$ where $\lambda\in \mathbb{D}$.
Since $C_\lambda+A\in \mathbb{B}_A\cap\mathbb{D}_{\leq 1}^{2\times n}$ for all $\lambda\in \mathbb{D}$, $|(\mathcal{N}+\varphi(A))\cap \mathcal{M}_1'|\geq 2$.
Since ${\rm rank}(\varphi(A))=2$, $\mathcal{N}+\varphi(A)\neq \mathcal{M}_1'$. Thus
Corollary \ref{Rectangular-1-13} implies that $\mathcal{N}$ is of type two.
By Lemma \ref{Rectangular-PID2-4}, ${\cal N}={\cal N}_1'Q_2$ where $Q_2\in GL_{n'}(\mathbb{D}')$.
Replacing  $\varphi$ by the map $X\mapsto \varphi(X)Q_2^{-1}$ ($X\in \mathbb{D}^{2\times n}$), we may assume with no loss of generality that
 $\mathcal{N}=\mathcal{N}_1'$. Thus
\begin{equation}\label{cx53gdgg3wevvv6}
\varphi\left(\mathbb{D}_{\leq 1}^{2\times n}\right)\subseteq \mathcal{M}_1', \ \ \varphi(\mathbb{B}_A)\subseteq \mathcal{N}_1'+\varphi(A).
\end{equation}

We prove  $\varphi(\mathbb{N}_A)\subseteq \mathcal{M}_1'$ as follows.

Write
 $\scriptsize\varphi(A)=
\left(
  \begin{array}{cc}
    a_{11} & A_{12}\\
    A_{21} & A_{22} \\
  \end{array}
\right)$ where $ a_{11}\in \mathbb{D}'$.
Since $C_\lambda+A\in \mathbb{N}_A\cap \mathbb{D}_{\leq 1}^{2\times n}$ for all $\lambda\in \mathbb{D}$,
$\varphi(C_\lambda+A)\in (\mathcal{N}_1'+\varphi(A))\cap \mathcal{M}_1'$ for all $\lambda\in \mathbb{D}$.
Hence $\varphi(A)$ is of the form
\begin{equation}\label{bv43ffgdb}
\varphi(A)=
\left(
  \begin{array}{cc}
    a_{11} & A_{12}\\
    A_{21} & 0 \\
  \end{array}
\right), \ \mbox{ where $a_{11}\in \mathbb{D}'$, $A_{21}\neq 0$ and $ A_{12}\neq 0$.}
\end{equation}
Otherwise,  $A_{22}\neq 0$ implies that $(\mathcal{N}_1'+\varphi(A))\cap \mathcal{M}_1'=\emptyset$, a contradiction.

Now, we affirm that ${\rm rank}(\varphi(Z))=1$ for every  $Z\in \mathbb{N}_A$.
Suppose that $Z\in\mathbb{N}_A$ and ${\rm rank}(\varphi(Z))=2$. Then ${\rm rank}(Z)=2$. We show a contradiction as follows.
By Theorem \ref{distinguishing-non01}, $\varphi(\mathbb{B}_Z)$ is an adjacent set. Thus, there exists a maximal set
$\mathcal{S}$   containing $0$ in  ${\mathbb{D}'}^{m'\times n'}$ such that
$$\varphi(\mathbb{B}_Z)\subseteq \mathcal{S}+\varphi(Z).$$
Write $\scriptsize Z=
\left(
  \begin{array}{c}
    \delta_1 \\
    \delta_2\\
  \end{array}
\right)$, where $\delta_1, \delta_2\in \mathbb{D}^n$ are left linearly independent.
Let
 $\scriptsize T_\lambda=
\left(
  \begin{array}{c}
    -\delta_1+\lambda\delta_2 \\
    0\\
  \end{array}
\right)$ where $\lambda\in \mathbb{D}$.
Since $T_\lambda+Z\in \mathbb{B}_Z\cap\mathbb{D}_{\leq 1}^{2\times n}$ for all $\lambda\in \mathbb{D}$, $|(\mathcal{S}+\varphi(Z))\cap \mathcal{M}_1'|\geq 2$.
By  ${\rm rank}(\varphi(Z))=2$, we have $\mathcal{S}+\varphi(Z)\neq\mathcal{M}_1'$. Hence
 Corollary \ref{Rectangular-1-13} implies that $\mathcal{S}+\varphi(Z)$ is a maximal set of type two.
Applying Corollary \ref{Rectangular-1-13} again, we get $\left|(\mathcal{S}+\varphi(Z))\cap(\mathcal{N}_1'+\varphi(A))\right|\leq 1$.
Let
$\scriptsize Z-A=
\left(\begin{array}{c}
    \gamma_1\\
    \gamma_2\\
  \end{array}
\right)\in \mathbb{D}_1^{2\times n}$, $\scriptsize Z_1=
\left(\begin{array}{c}
    0\\
    \gamma_1-\gamma_2\\
  \end{array}
\right)$,
$\scriptsize Z_2=
\left(\begin{array}{c}
    \gamma_2-\gamma_1\\
    0\\
  \end{array}
\right)$,
$\scriptsize Y_1=
\left(\begin{array}{c}
    \gamma_1\\
    \gamma_1\\
  \end{array}
\right)$ and
$\scriptsize Y_2=
\left(\begin{array}{c}
    \gamma_2\\
    \gamma_2\\
  \end{array}
\right)$. Then $Z_i+Z=Y_i+A$, $i=1,2$. Moreover, $Z_1+Z\sim Z_2+Z$. Clearly, $Z_i+Z\in\mathbb{B}_Z$ and $Y_i+A\in\mathbb{B}_A$, $i=1,2$.
It follows that $|\mathbb{B}_Z\cap \mathbb{B}_A|\geq 2$ and hence $\left|(\mathcal{S}+\varphi(Z))\cap(\mathcal{N}_1'+\varphi(A))\right|\geq 2$,
a contradiction.
Therefore,  ${\rm rank}(\varphi(Z))\neq 2$ for all $Z\in \mathbb{N}_A$. By ${\rm rank}(\varphi(A))=2$ and $\varphi(A)\sim \varphi(Z)$, one gets
${\rm rank}(\varphi(Z))=1$  for every $Z\in \mathbb{N}_A$.

Note that $\mathbb{N}_A\subset \mathbb{B}_A$.  By (\ref{cx53gdgg3wevvv6}) and (\ref{bv43ffgdb}), we obtain
\begin{equation}\label{c543fuadlj970}
\varphi(\mathbb{N}_A)\subseteq\left\{
\left(
  \begin{array}{cc}
    y & A_{12}\\
    0 & 0 \\
  \end{array}
\right): y\in \mathbb{D}'\right\}\subset\mathcal{M}_1'.
\end{equation}

{\bf Step 2}. \ Let $B\in {\mathbb{D}}_2^{2\times n}$ with $B\neq A$ and ${\rm rank}(\varphi(B))=2$. By Theorem \ref{distinguishing-non01},
$\varphi(\mathbb{B}_B)$ is an adjacent set. Thus,
there is a maximal set $\mathcal{N}''$ containing $0$ in ${\mathbb{D}'}^{m'\times n'}$ such that
\begin{equation}\label{V4NM543DDS}
\varphi(\mathbb{B}_B)\subseteq \mathcal{N}''+\varphi(B).
\end{equation}
 Write
$\scriptsize B=
\left(
  \begin{array}{c}
    \beta_1 \\
    \beta_2\\
  \end{array}
\right)$ where  $\beta_1, \beta_2\in \mathbb{D}^n$ are left linearly independent.
 Set
 $\scriptsize C_\lambda=
\left(
  \begin{array}{c}
    -\beta_1+\lambda\beta_2 \\
    0\\
  \end{array}
\right)$ where $\lambda\in \mathbb{D}$.
Since $C_\lambda+B\in\mathbb{B}_B\cap\mathbb{D}_{\leq 1}^{2\times n}$ for all $\lambda\in \mathbb{D}$, $|(\mathcal{N}''+\varphi(B))\cap \mathcal{M}_1'|\geq 2$.
Since ${\rm rank}(\varphi(B))=2$, $\mathcal{N}''+\varphi(B)\neq \mathcal{M}_1'$.
It follows from Corollary \ref{Rectangular-1-13}  that $\mathcal{N}''$ is of type two.
We prove $\mathcal{N}''=\mathcal{N}_1'$ by contradiction as follows.

Suppose that $\mathcal{N}''\neq\mathcal{N}_1'$. Write $\mathcal{N}''=\left\{(yq_1,yq_2,\ldots, yq_n): y\in \,^m\mathbb{D}\right\}$ where $q_j\in \mathbb{D}'$ and
$(q_2,\ldots,q_n)\neq 0$. Then there is $\scriptsize Q_2=\left(
                                                \begin{array}{cc}
                                                  1 & 0 \\
                                                  * & *\\
                                                \end{array}
                                              \right)\in GL_{n'}(\mathbb{D}')$ such that
$\mathcal{N}''=\mathcal{N}_2'Q_2$ and $\mathcal{N}'_1=\mathcal{N}_1'Q_2$.
Modifying the map $\varphi$ by the map $X\mapsto \varphi(X)Q_2^{-1}$ ($X\in \mathbb{D}^{2\times n}$). We may assume with no loss of
generality  that $\mathcal{N}''=\mathcal{N}_2'$
and (\ref{cx53gdgg3wevvv6})-(\ref{c543fuadlj970}) hold. Thus
$$\varphi(\mathbb{B}_B)\subseteq \mathcal{N}'_2+\varphi(B).$$
Write
 $\scriptsize\varphi(B)=
\left(
  \begin{array}{ccc}
    b_{11}&b_{12} & B_{13}\\
    B_{21}&B_{22} & B_{23} \\
  \end{array}
\right)$ where $b_{11}, b_{12}\in \mathbb{D}'$.
Since $C_\lambda+B\in \mathbb{B}_B\cap \mathbb{D}_{\leq 1}^{2\times n}$ for all $\lambda\in \mathbb{D}$,
$\varphi(C_\lambda+B)\in (\mathcal{N}_2'+\varphi(B))\cap \mathcal{M}_1'$ for all $\lambda\in \mathbb{D}$.
Hence $\varphi(B)$ is of the form
\begin{equation}\label{bv43CCFFNZZ}
\varphi(B)=
\left(
  \begin{array}{ccc}
  b_{11}&b_{12} & B_{13}\\
    0&B_{22} & 0 \\
  \end{array}
\right),
\end{equation}
where  $B_{22}\neq 0$ and $ (b_{11}, B_{13})\neq (0, 0)$.
Otherwise, we have $(B_{21}, B_{23})\neq (0,0)$, which implies that $(\mathcal{N}_2'+\varphi(B))\cap \mathcal{M}_1'=\emptyset$, a contradiction.
By $\varphi(B)\notin \mathcal{M}_1'$ and (\ref{c543fuadlj970}), we have $B\notin \mathbb{N}_A$ and hence $d(A,B)=2$.
 Write
 $\scriptsize B-A=
\left(
  \begin{array}{c}
    \nu_1\\
    \nu_2 \\
  \end{array}
\right)$, where $\nu_1,\nu_2\in \mathbb{D}^n$ are left linearly independent.
 Then
 $\scriptsize B\sim
\left(
  \begin{array}{c}
    \nu_1+\lambda\nu_2\\
     0\\
  \end{array}
\right)+A\sim A$ for all $\lambda\in\mathbb{D}$. Let $A_{12}=(a_{12}, A_{13})$ where $a_{12}\in \mathbb{D}'$. By (\ref{c543fuadlj970}), one has that
\begin{equation}\label{vcx344nvutiy}
\mbox{$\varphi(B)=
\left(
  \begin{array}{ccc}
  b_{11}&b_{12} & B_{13}\\
    0&B_{22} & 0 \\
  \end{array}
\right)\sim
\varphi\left(\left(
  \begin{array}{c}
    \nu_1+\lambda\nu_2\\
     0\\
  \end{array}
\right)+A\right)=:\left(
  \begin{array}{ccc}
    \lambda^\sigma &a_{12}&A_{13}\\
     0 &0&0\\
  \end{array}
\right)$, \ for all $\lambda\in\mathbb{D}$,}
\end{equation}
 where $\sigma$ is an injection from $\mathbb{D}$ to $\mathbb{D}'$.
Hence
 $$\mbox{${\rm rank}\left(
               \begin{array}{ccc}
                 b_{11}-\lambda^\sigma &b_{12}-a_{12} & B_{13}-A_{13} \\
                 0& B_{22}&0 \\
               \end{array}
             \right)=1$, \ for all $\lambda\in\mathbb{D}$,}$$
which is a contradiction because $B_{22}\neq 0$. Therefore,  we must have $\mathcal{N}''=\mathcal{N}_1'$.

 Similar to the proof of (\ref{vcx344nvutiy}), we have that $d(A,B)=2$ and
\begin{equation}\label{vcxn44vutiy}
\mbox{$\varphi(B)=:
\left(
  \begin{array}{cc}
  b_{11}& B_{12}'\\
    B_{21}'&B_{22}'\\
  \end{array}
\right)\sim
\varphi\left(\left(
  \begin{array}{c}
    \nu_1+\lambda\nu_2\\
     0\\
  \end{array}
\right)+A\right)=:\left(
  \begin{array}{cc}
    \lambda^\mu &A_{12}\\
     0 &0\\
  \end{array}
\right)$ \ for all $\lambda\in\mathbb{D}$,        }
\end{equation}
where $\mu$ is an injection from $\mathbb{D}$ to $\mathbb{D}'$. Therefore,
$$\mbox{${\rm rank}\left(
               \begin{array}{cc}
                 b_{11}-\lambda^\mu & B_{12}'-A_{12} \\
                 B_{21}'& B_{22}' \\
               \end{array}
             \right)=1$, \ for all $\lambda\in\mathbb{D}$.}$$
Using Lemma \ref{Matrix-PID5-4bb} and ${\rm rank}(\varphi(B))=2$, it is easy to see that $B_{22}'=0$ and $B_{12}'=A_{12}$.
Consequently, $\mathcal{N}_1'+\varphi(B)=\mathcal{N}_1'+\varphi(A)$, and (\ref{V4NM543DDS}) implies that
\begin{equation}\label{VC4GHFJII88}
\mbox{$\varphi(\mathbb{N}_B)\subseteq \varphi(\mathbb{B}_B)\subseteq \mathcal{N}_1'+\varphi(A)$, \
if $B\in {\mathbb{D}}_2^{2\times n}$ with $B\neq A$ and ${\rm rank}(\varphi(B))=2$.}
\end{equation}

{\bf Step 3}. \
Let $T\in {\mathbb{D}}_2^{2\times n}$ with ${\rm rank}(\varphi(T))=1$.  We show $\varphi(T)\in \mathcal{M}_1'$ as follows.

Put
$\scriptsize\varphi(T)=
\left(
  \begin{array}{cc}
    t_{11} & T_{12}\\
    T_{21} & T_{22} \\
  \end{array}
\right)$ where $t_{11}\in \mathbb{D}'$.
If $T\sim A$, then $T\in\mathbb{N}_A$ and hence
$\varphi(T)\in \mathcal{M}_1'$ by (\ref{c543fuadlj970}).
 Now we assume  $d(T, A)=2$. Write
 $\scriptsize T-A=
\left(
  \begin{array}{c}
    \delta_1\\
    \delta_2 \\
  \end{array}
\right)$ where $\delta_1, \delta_2\in \mathbb{D}^n$ are left linearly independent.
 Then
 $\scriptsize T\sim
\left(
  \begin{array}{c}
    \delta_1+\lambda \delta_2\\
     0\\
  \end{array}
\right)+A\sim A$ for all $\lambda\in\mathbb{D}$. By (\ref{c543fuadlj970}), we can let
 $\scriptsize
\varphi\left(\left(
  \begin{array}{c}
    \delta_1+\lambda\delta_2\\
     0\\
  \end{array}
\right)+A\right)=\left(
  \begin{array}{cc}
    \lambda^\tau & A_{12}\\
     0 &0\\
  \end{array}
\right)$ for all $\lambda\in\mathbb{D}$, where $\tau$ is an injection from $\mathbb{D}$ to $\mathbb{D}'$.
 Since $\scriptsize\varphi(T)\sim
 \left(
  \begin{array}{cc}
    \lambda^\tau & A_{12}\\
     0&0\\
  \end{array}
\right)$ for all $\lambda\in\mathbb{D}$, we get that
 $$\mbox{${\rm rank}\left(
               \begin{array}{cc}
                 t_{11}-\lambda^\tau & T_{12}-A_{12} \\
                 T_{21}& T_{22} \\
               \end{array}
             \right)={\rm rank}\left(
               \begin{array}{cc}
                 t_{11} & T_{12} \\
                 T_{21}& T_{22} \\
               \end{array}
             \right)=1$, \ for all $\lambda\in\mathbb{D}$.}$$
 By Lemma \ref{Matrix-PID5-4bb} and $A_{12}\neq 0$, we obtain $(T_{21},  T_{22})=0$. Then
 $\scriptsize\varphi(T)=
\left(
  \begin{array}{cc}
    t_{11} & T_{12}\\
  0 & 0 \\
  \end{array}
\right)\in\mathcal{M}_1'$.
Hence
 \begin{equation}\label{kjh5464dfsf}
 \mbox{$\varphi(T)\in\mathcal{M}_1'$, \ if $T\in {\mathbb{D}}_2^{2\times n}$ and ${\rm rank}(\varphi(T))=1$.}
 \end{equation}

By (\ref{cx53gdgg3wevvv6}), (\ref{VC4GHFJII88}) and (\ref{kjh5464dfsf}),  we have that $\varphi(X)\in\mathcal{N}_1'+\varphi(A)$ whenever ${\rm rank}(\varphi(X))= 2$, and
 $\varphi(X)\in\mathcal{M}_1'$ whenever ${\rm rank}(\varphi(X))\neq 2$. Thus
$$\varphi(\mathbb{D}^{2\times n})\subseteq \mathcal{M}_1'\cup(\mathcal{N}_1'+\varphi(A)).$$

Recall that $\scriptsize \varphi(A)=\left(
\begin{array}{cc}
 a_{11} & A_{12} \\
 A_{21} & 0 \\
 \end{array}
 \right)$ where $A_{21}\neq 0$ and $ A_{12}\neq 0$.
Let $\scriptsize R=\left(
\begin{array}{cc}
 a_{11} & A_{12} \\
 0 & 0 \\
 \end{array}
 \right)\in \mathcal{M}_1'$. Then  ${\rm rank}(R)=1$, $\mathcal{M}_1'=\mathcal{M}_1'+R$ and $\mathcal{N}_1'+\varphi(A)=\mathcal{N}_1'+R$, and hence
$\mathcal{M}_1'\cup(\mathcal{N}_1'+\varphi(A))=(\mathcal{M}_1'+R)\cup(\mathcal{N}_1'+R)$. Therefore, (\ref{vcadgjr772m}) holds.
By (\ref{cx53gdgg3wevvv6}) and (\ref{c543fuadlj970}), we have  $\varphi\left(\mathbb{D}_{\leq 1}^{2\times n}\right)\subseteq \mathcal{M}_1'$ and
$\varphi(\mathbb{N}_A)\subseteq \mathcal{M}_1'\cap(\mathcal{N}_1'+R)$. Thus, there is $\mathcal{M}'\in \{\mathcal{M}, \mathcal{N}\}$ such that
$\varphi\left(\mathbb{D}_{\leq 1}^{2\times n}\right)\subseteq \mathcal{M}'$ and $\varphi(\mathbb{N}_A)\subseteq \mathcal{M}'\cap(\mathcal{N}'+R)$
where $\{\mathcal{M}', \mathcal{N}'\}=\{\mathcal{M}, \mathcal{N}\}$.
$\qed$

Now, we prove Theorem \ref{degenerate7654} as follows.

\noindent{\bf Proof of Theorem \ref{degenerate7654}.} \  Note that the map $\varphi(X)-\varphi(0)$ is also a degenerate graph homomorphism. By the
 map $X\mapsto \varphi(X)-\varphi(0)$, we may assume with no loss of generality that $\varphi(0)=0$.
We prove this theorem only for the case of $m=2$; the case of $n=2$ is similar. From now on we assume that  $m={\rm min}\{m,n\}=2$.
Clearly, $1\leq {\rm diam}(\varphi(\mathbb{D}^{2\times n}))\leq 2$.

{\bf Case 1.} \  There exists $A\in {\mathbb{D}}_2^{2\times n}$ such that ${\rm rank}(\varphi(A))=2$.
By Lemma \ref{degenerate675235utigh},  this theorem holds.

{\bf Case 2.} \ There is no $A\in {\mathbb{D}}_2^{2\times n}$ such that ${\rm rank}(\varphi(A))=2$. Then
\begin{equation}\label{ret5fsgvd77ssf}
\mbox{${\rm rank}(\varphi(X))\leq 1$, \  for all $X\in{\mathbb{D}}^{2\times n}$.}
\end{equation}
If ${\rm diam}(\varphi(\mathbb{D}^{2\times n}))=1$, then $\varphi(\mathbb{D}^{2\times n})$ is an adjacent set, and hence this theorem  holds.
From now on we assume that  ${\rm diam}(\varphi(\mathbb{D}^{2\times n}))=2$. Then, there are two matrices $B_1,B_2\in\mathbb{D}^{2\times n}$ such that
$d(\varphi(B_1), \varphi(B_2))=2=d(B_1,B_2)$.
Let $\psi(X)=\varphi(X+B_1)-\varphi(B_1)$, $X\in{\mathbb{D}}^{2\times n}$.
Then $\psi:  {\mathbb{D}}^{2\times n}\rightarrow  {\mathbb{D}'}^{m'\times n'}$ is a graph homomorphism with $\psi(0)=0$  and  ${\rm rank}(\psi(B_2-B_1))=2$.

{\em Subcase 2.1.} \ $\psi$ is degenerate.
By Lemma \ref{degenerate675235utigh}, there are two fixed maximal sets $\mathcal{M}$ and $\mathcal{N}$ of different
types containing $0$,  such that
$\psi(\mathbb{D}^{2\times n})\subseteq (\mathcal{M}+R_1)\cup(\mathcal{N}+R_1)$,
where $R_1\in {\mathbb{D}'}^{m'\times n'}$ is fixed and ${\rm rank}(R_1)=1$.
Hence $\varphi(\mathbb{D}^{2\times n})\subseteq (\mathcal{M}+R_1+\varphi(A))\cup(\mathcal{N}+R_1+\varphi(B_1))$ and
 this theorem  holds.

{\em Subcase 2.2.} \ $\psi$ is non-degenerate. By Lemma \ref{s4245549bcnfte},
  there exist  invertible matrices  $T_1,T_2$ over $\mathbb{D}'$ and invertible matrices $T_3,T_4$ over $\mathbb{D}$,  such that
either
\begin{equation}\label{b6teyg646fd}
\mbox{$\psi(T_3^{-1}\mathcal{M}_iT_4^{-1})\subseteq T_1\mathcal{M}_i'T_2$ \ with \ $\psi(T_3^{-1}\mathcal{N}_iT_4^{-1})\subseteq T_1\mathcal{N}_i'T_2$ $(i=1,2)$,}
\end{equation}
or
\begin{equation}\label{b66557ghghf6}
\mbox{$\psi(T_3^{-1}\mathcal{M}_iT_4^{-1})\subseteq T_1\mathcal{N}_i'T_2$ \ with \ $\psi(T_3^{-1}\mathcal{N}_iT_4^{-1})\subseteq T_1\mathcal{M}_i'T_2$ $(i=1,2)$.}
\end{equation}
Without loss of generality, we may assume that  (\ref{b6teyg646fd}) holds and all $T_1, T_2,T_3,T_4$ are identity matrices. Then
$\psi(\mathcal{M}_i)\subseteq\mathcal{M}_i'$ \  with \ $\psi(\mathcal{N}_i)\subseteq \mathcal{N}_i'$, $i=1,2$.
Since $\mathcal{M}_i\cap\mathcal{N}_j=\mathbb{D}E_{ij}$ and  $\mathcal{M}_i'\cap\mathcal{N}_j'=\mathbb{D}'E_{ij}'$, $1\leq i,j\leq 2$,
we can let $\psi(xE_{ij})=x^{\sigma_{ij}}E_{ij}'$,  $x\in \mathbb{D}$,
where $\sigma_{ij}: \mathbb{D}\rightarrow \mathbb{D}'$ is an injective map with $0^{\sigma_{ij}}=0$.
By the definition of $\psi$, we obtain that
$$\mbox{$\varphi(xE_{ij}+A)=x^{\sigma_{ij}}E_{ij}'+\varphi(B_1)$, \ $x\in \mathbb{D}$, $i, j=1,2$.}$$
Using (\ref{ret5fsgvd77ssf}), we have ${\rm rank}\left(x^{\sigma_{ij}}E_{ij}'+\varphi(B_1)\right)\leq 1$,  $x\in \mathbb{D}$, $i,j=1,2$.
Applying Lemma \ref{Matrix-PID5-4bb}, it is easy to verify that $\varphi(B_1)=0$. Thus,
${\rm rank}(\varphi(B_2))=d(\varphi(B_2), \varphi(B_1))=2$, a contradiction to (\ref{ret5fsgvd77ssf}).
Therefore, Subcase 2.2 does not happen.
$\qed$

By Theorem \ref{distinguishing-non01}, the following corollary is obvious.

\begin{cor}\label{disishing-non-cor01} \ Let $\mathbb{D}, \mathbb{D}'$ be division rings  with $|\mathbb{D}|\geq 4$,
 and let $m,n,m', n'\geq 2$ be integers with $m',n'\geq {\rm min}\{m,n\}$. Suppose  $\varphi:{\mathbb{D}}^{m\times n}\rightarrow  {\mathbb{D}'}^{m'\times n'}$
is a graph homomorphism with $\varphi(0)=0$,  and there exists $A_0\in {\mathbb{D}}^{m\times n}$ such that ${\rm rank}(\varphi(A_0))={\rm min}\{m,n\}$.
If $\varphi\left(\mathbb{D}_{\leq 1}^{m\times n} \right)$ or $\varphi(\mathbb{B}_{A_0})$ is not an adjacent set, then $\varphi$ is non-degenerate.
\end{cor}

The following result is a generalization of \cite[Theorem 1.1]{Huang-Semrl-2014} (which is due to Huang and \v{S}emrl).

\begin{cor}\label{non-degenerate-ddD}  Let $\mathbb{D}, \mathbb{D}'$ be division rings with $|\mathbb{D}|\geq 4$,
and let  $m', n'\geq 2$. Suppose
$\varphi:  {\mathbb{D}}^{2\times 2}\rightarrow  {\mathbb{D}'}^{m'\times n'}$
is a graph homomorphism. Then either $\varphi$ is a distance preserving map (which is of the form either (\ref{3654665cnc0c2}) or (\ref{CX3mmbb660c2})),
or there are two fixed maximal sets $\mathcal{M}$ and $\mathcal{N}$ of different types in ${\mathbb{D}'}^{m'\times n'}$,
such that $0\in \mathcal{M}\cap\mathcal{N}$ and
\begin{equation}\label{fg328kjhyti00}
\varphi(\mathbb{D}^{2\times 2})\subseteq (\mathcal{M}+R)\cup(\mathcal{N}+R),
\end{equation}
where $R\in{\mathbb{D}'}^{m'\times n'}$ is fixed.
\end{cor}
\proof
Assume that ${\rm diam}(\varphi(\mathbb{D}^{2\times 2}))=1$. Then $\varphi(\mathbb{D}^{2\times 2})$ is an adjacent set, and hence  (\ref{fg328kjhyti00}) holds.
From now on we assume that  ${\rm diam}(\varphi(\mathbb{D}^{2\times 2}))=2$.  Then
 there are two matrices $A_0,B_0\in\mathbb{D}^{2\times 2}$ such that
 $d(\varphi(A_0),\varphi(B_0))={\rm rank}(\varphi(B_0)-\varphi(A_0))=2$.
Let $\psi(X)=\varphi(X+A_0)-\varphi(A_0)$, $X\in{\mathbb{D}}^{2\times 2}$. Then $\psi:  {\mathbb{D}}^{2\times 2}\rightarrow  {\mathbb{D}'}^{m'\times n'}$
is a graph homomorphism with $\psi(0)=0$  and  ${\rm rank}(\psi(B_0-A_0))=2$. By Corollary \ref{stronglydegenerate2},
 either $\psi$ is a distance preserving map or  $\psi\left(\mathbb{D}_{\leq 1}^{2\times 2} \right)$ is an adjacent set.

  Assume that $\psi$ is a distance preserving map. Then $\varphi$ is also a distance preserving map.
 By Corollary \ref{non-degenerate-cc02}, $\varphi$ is of the form either (\ref{3654665cnc0c2}) or (\ref{CX3mmbb660c2}).
Now, we assume that $\psi\left(\mathbb{D}_{\leq 1}^{2\times 2} \right)$ is an adjacent set. By Theorem \ref{degenerate7654} or Lemma \ref{degenerate675235utigh},
there are two fixed maximal sets $\mathcal{M}$ and $\mathcal{N}$ of different types in ${\mathbb{D}'}^{m'\times n'}$,
 such that $0\in \mathcal{M}\cap\mathcal{N}$ and (\ref{fg328kjhyti00}) holds.
$\qed$

Finally, we discuss the case of finite fields. For the case of finite fields, we have the following better results.

\begin{thm}\label{finitefieldsooooa} \ Let  $\mathbb{D}, \mathbb{D}'$ be two finite fields with  $|\mathbb{D}|>|\mathbb{D}'|$, and let
$m,n,m', n'\geq2$ be integers.  Then every graph homomorphism from $\mathbb{D}^{m\times n}$ to ${\mathbb{D}'}^{m'\times n'}$ is a (vertex) colouring.
 \end{thm}
\proof
 Suppose  $\varphi:  \mathbb{D}^{m\times n}\rightarrow  {\mathbb{D}'}^{m'\times n'}$ is a graph homomorphism. Without loss of generality, we may assume that
 $\varphi(0)=0$ and $\varphi(\mathcal{M}_1)\subseteq \mathcal{M}_1'$.
 If  $\ell$ is a line in $AG(\mathcal{M}_1)$, then $|\mathbb{D}|>|\mathbb{D}'|$ implies that $\varphi(\ell)$ is not contained in any line in $AG(\mathcal{M}_1')$,
 and hence  $\varphi(\ell)$ contains at least three  noncollinear points in $AG(\mathcal{M}_1')$.

 Let $\scriptsize A=\left(
                \begin{array}{c}
                  \alpha_1 \\
                  \vdots \\
                  \alpha_m \\
                \end{array}
              \right)\in \mathbb{D}_1^{m\times n}$. Without loss of generality, we  assume  $\alpha_m\neq 0$.  Then $\alpha_i=k_i\alpha_m$, $i=1,\ldots, m-1$.
We have $\scriptsize A\sim \left(
                \begin{array}{c}
                  \lambda \alpha_m \\
                  0 \\
                \end{array}
              \right)$ for all $\lambda\in \mathbb{D}$. Since $\scriptsize\ell_1:=\left\{ \left(
                \begin{array}{c}
                  \lambda \alpha_m \\
                  0 \\
                \end{array}
              \right): \lambda\in \mathbb{D} \right\}$ is a line in  $AG(\mathcal{M}_1)$, it follows from above result that $\varphi(A)$ is adjacent with
three  noncollinear points in $AG(\mathcal{M}_1')$. By Lemma \ref{gdkj7ssb33}, $\varphi(A)\in \mathcal{M}_1'$ for all $A\in \mathbb{D}_1^{m\times n}$. Therefore, we obtain
$\varphi\left(\mathbb{D}_{\leq 1}^{m\times n}\right)\subseteq  \mathcal{M}_1'$.

Suppose that $\varphi\left(\mathbb{D}_{\leq \,k-1}^{m\times n}\right)\subseteq  \mathcal{M}_1'$ where $2\leq k\leq m$. Let $B\in \mathbb{D}_k^{m\times n}$.
Note that $P(\mathbb{D}_{\leq \,k-1}^{m\times n})Q=\mathbb{D}_{\leq \,k-1}^{m\times n}$ for any $P\in GL_m(\mathbb{D}$ and $Q\in GL_n(\mathbb{D}$.
Without loss of generality, we can assume that $\scriptsize B=\left(
\begin{array}{c}
                B_1\\
                  0 \\
                \end{array}
              \right)$ where $\scriptsize B_1=\left(
                \begin{array}{c}
                  \alpha_1 \\
                  \vdots \\
                  \alpha_k \\
                \end{array}
              \right)\in \mathbb{D}_k^{k\times n}$. Put
$\scriptsize C_1=\left(
                \begin{array}{c}
                  \alpha_1 \\
                  \vdots \\
                  \alpha_{k-1} \\
                \end{array}
              \right)$. Then
$\scriptsize B\sim \left(
        \begin{array}{c}
          C_1+{\scriptsize\left(
                \begin{array}{c}
                  \lambda \alpha_k \\
                  0 \\
                \end{array}
              \right)}
           \\
          0 \\
        \end{array}
      \right)$ for all $\lambda\in \mathbb{D}$.
Let
$\scriptsize \ell_2=\left\{\left(
        \begin{array}{c}
          C_1+{\scriptsize\left(
                \begin{array}{c}
                  \lambda \alpha_k \\
                  0 \\
                \end{array}
              \right)}
           \\
          0 \\
        \end{array}
      \right): \lambda\in\mathbb{D}  \right\}\subset \mathbb{D}_{\leq \,k-1}^{m\times n}$.
Then $\ell_2$ is a line in some affine geometry on a maximal set in $\mathbb{D}^{m\times n}$. By the induction hypothesis, we have $\varphi(\ell_2)\subseteq \mathcal{M}_1'$.
Since $|\mathbb{D}|>|\mathbb{D}'|$,  $\varphi(\ell_2)$ contains  at least three  noncollinear points in $AG(\mathcal{M}_1')$. By Lemma \ref{gdkj7ssb33} and $\varphi(B)\sim Y$ for any $Y\in\varphi(\ell_2)$,
we obtain   $\varphi(B)\in \mathcal{M}_1'$ for any $B\in \mathbb{D}_k^{m\times n}$. Hence $\varphi\left(\mathbb{D}_{\leq \, k}^{m\times n}\right)\subseteq  \mathcal{M}_1'$.
 Taking $k=m$, we get $\varphi\left(\mathbb{D}^{m\times n}\right)\subseteq\mathcal{M}_1'$.
 $\qed$

Let $\lceil x \rceil$ denote the smallest integer at least as large as $x$.

\begin{thm}\label{finitefieldsoooob} \ Let  $\mathbb{D}, \mathbb{D}'$ be two finite fields with $4\leq |\mathbb{D}|\leq|\mathbb{D}'|\leq(|\mathbb{D}|-1)\lceil (|\mathbb{D}|+1)/2 \rceil$, and let
$m,n,m', n'\geq2$ be integers.  Suppose  $\varphi:  \mathbb{D}^{m\times n}\rightarrow  {\mathbb{D}'}^{m'\times n'}$
is a degenerate graph homomorphism. Then there are two fixed maximal sets $\mathcal{M}$ and $\mathcal{N}$ of different types containing $0$ in ${\mathbb{D}'}^{m'\times n'}$,
 such that
\begin{equation}\label{tre543bcn877m}
\mbox{$\varphi(\mathbb{D}^{m\times n})\subseteq (\mathcal{M}+R)\cup(\mathcal{N}+R)$,}
 \end{equation}
 where $R\in {\mathbb{D}'}^{m'\times n'}$ is fixed.
 \end{thm}
\proof
When ${\rm min}\{m,n\}=2$, this theorem holds by Theorem \ref{degenerate7654}. From now on we assume that $m, n>2$.
 By Theorem \ref{degenerate7654},
 there are two fixed maximal sets $\mathcal{M}$ and $\mathcal{N}$ of different types containing $0$ in ${\mathbb{D}'}^{m'\times n'}$,  such that
$\scriptsize\varphi\left(
                \begin{array}{c}
                  \mathbb{D}^{2\times n} \\
                  0 \\
                \end{array}
              \right)
\subseteq (\mathcal{M}+R)\cup(\mathcal{N}+R)$,  where $R\in {\mathbb{D}'}^{m'\times n'}$ is fixed.
Suppose  that $3\leq k\leq m$ and
\begin{equation}\label{hfg45978bcada9}
\mbox{$\varphi\left(
                \begin{array}{c}
                  \mathbb{D}^{(k-1)\times n} \\
                  0 \\
                \end{array}
              \right)
\subseteq (\mathcal{M}+R)\cup(\mathcal{N}+R)$.}
\end{equation}
 We prove  $\scriptsize\varphi\left(
                \begin{array}{c}
                  \mathbb{D}^{k\times n} \\
                  0 \\
                \end{array}
              \right)
\subseteq (\mathcal{M}+R)\cup(\mathcal{N}+R)$ as follows.

Let $\scriptsize A=\left(
\begin{array}{c}
                A_1\\
                  0 \\
                \end{array}
              \right)\in \mathbb{D}^{m\times n}$, where $\scriptsize A_1=\left(
                \begin{array}{c}
                  \alpha_1 \\
                  \vdots \\
                  \alpha_k \\
                \end{array}
              \right)\in \mathbb{D}^{k\times n}$ with $\alpha_k\neq 0$. Put  $\scriptsize B_1=\left(
                \begin{array}{c}
                  \alpha_1 \\
                  \vdots \\
                  \alpha_{k-1} \\
                \end{array}
              \right)$. Then
\begin{equation}\label{FS4TETB43dfg}
A\sim \left(
        \begin{array}{c}
          B_1+{\scriptsize \lambda\left(
                \begin{array}{c}
                  \alpha_k \\
                  0 \\
                  0 \\
                \end{array}
              \right)}
           \\
          0 \\
        \end{array}
      \right), \ \   A\sim \left(
        \begin{array}{c}
          B_1+{\scriptsize \lambda\left(
                \begin{array}{c}
                    0_{1,n} \\
                \alpha_k \\
                  0 \\
                \end{array}
              \right)}
           \\
          0 \\
        \end{array}
      \right),  \ \  A\sim \left(
        \begin{array}{c}
          B_1+{\scriptsize\lambda\left(
                \begin{array}{c}
                  d\alpha_k \\
                  \alpha_k \\
                  0 \\
                \end{array}
              \right)}
           \\
          0 \\
        \end{array}
      \right), \ \ \lambda\in \mathbb{D}, d\in \mathbb{D}^*.
\end{equation}
Let
$$\ell_1=\left\{\left(
        \begin{array}{c}
          B_1+{\small \lambda\left(
                \begin{array}{c}
                  \alpha_k \\
                  0 \\
                  0 \\
                \end{array}
              \right)}
           \\
          0 \\
        \end{array}
      \right): \lambda\in\mathbb{D}  \right\}, \ \ \  \ell_2=\left\{\left(
        \begin{array}{c}
          B_1+{\small \lambda\left(
                \begin{array}{c}
                   0_{1,n} \\
                   \alpha_k \\
                  0 \\
                \end{array}
              \right)}
           \\
          0 \\
        \end{array}
      \right): \lambda\in\mathbb{D}  \right\},$$
$$ \ell_d'=\left\{\left(
        \begin{array}{c}
          B_1+{\small\lambda\left(
                \begin{array}{c}
           d\alpha_k \\
            \alpha_k \\
                  0 \\
                \end{array}
              \right)}
           \\
          0 \\
        \end{array}
      \right): \lambda\in\mathbb{D}  \right\}, \ d\in\mathbb{D}^*.$$
Then $\ell_1, \ell_2, \ell_d'$ ($d\in \mathbb{D}^*$) are $|\mathbb{D}|+1$ distinct lines in  affine geometries on $|\mathbb{D}|+1$  maximal sets in $\mathbb{D}^{m\times n}$.
Moreover, $\ell_1\cap \ell_2=\ell_i\cap \ell_d'=\ell_{d_1}'\cap\ell_{d_2}'=\scriptsize \left(
                \begin{array}{c}
               B_1\\
                  0 \\
                \end{array}
              \right)$ for $i=1,2$ and $d,d_1,d_2\in\mathbb{D}^*$ with $d_1\neq d_2$.
Clearly, $\varphi(\ell_i)$ or $\ell_d'$ is contained in a maximal set in ${\mathbb{D}'}^{m'\times n'}$, $i=1,2$, $d\in\mathbb{D}^*$. Since $\scriptsize \ell_i, \ell_d'\subseteq \left(
                \begin{array}{c}
                  \mathbb{D}^{(k-1)\times n} \\
                  0 \\
                \end{array}
              \right)$,
it follows from (\ref{hfg45978bcada9}) that
$\varphi(\ell_i)\subseteq \mathcal{M}+R$ or $\varphi(\ell_i)\subseteq \mathcal{N}+R$, $i=1,2$; $\varphi(\ell_d')\subseteq \mathcal{M}+R$ or $\varphi(\ell_d')\subseteq \mathcal{N}+R$, $d\in \mathbb{D}^*$.

Let either $\mathcal{M}'=\mathcal{M}+R$ or $\mathcal{M}'=\mathcal{N}+R$.
Suppose that there is some $k$ or $d$ such that
$\varphi(\ell_k)$ or $\ell_d'$ is not contained in a line in $AG(\mathcal{M}')$.
Then, from (\ref{FS4TETB43dfg})  $\varphi(A)$ is adjacent with  three  noncollinear points in $AG(\mathcal{M}')$.  By Lemma \ref{gdkj7ssb33}, we get  $\varphi(A)\in \mathcal{M}'$. Hence
$\varphi(A)\in (\mathcal{M}+R)\cup(\mathcal{N}+R)$.

Now, we assume that $\varphi(\ell_k)$ [resp. $\ell_d'$] is contained in a line in $AG(\mathcal{M}')$ for  $k=1,2$ [resp. all $d\in \mathbb{D}^*$].
Let $X, Y\in \ell_1\cup\ell_2\cup_{d\in \mathbb{D}^*}\ell_d'$ with $X\neq Y$. Then
 $X\sim Y$ and hence $\varphi(X)\sim\varphi(Y)$ and $\varphi(X)\neq\varphi(Y)$. By (\ref{FS4TETB43dfg}), we have $\varphi(A)\sim\varphi(X)$ and $\varphi(A)\sim\varphi(Y)$.
 Clearly, $\mathcal{M}+R$ or $\mathcal{N}+R$
 contains $\lceil (|\mathbb{D}|+1)/2 \rceil$ elements of the set $\left\{\varphi(\ell_1), \varphi(\ell_2), \varphi(\ell_d'), d\in \mathbb{D}^*\right\}$.
 Without loss of generality, we  assume that $\mathcal{M}+R$ contains $\lceil (|\mathbb{D}|+1)/2 \rceil$ elements of the set
 $\left\{\varphi(\ell_1), \varphi(\ell_2), \varphi(\ell_d'), d\in \mathbb{D}^*\right\}$. Then,  $\mathcal{M}+R$ contains $(|\mathbb{D}|-1)\lceil (|\mathbb{D}|+1)/2 \rceil +1$
 distinct points in $\varphi(\ell_1)\cup\varphi(\ell_2)\cup_{d\in \mathbb{D}^*}\varphi(\ell_d')$. Since $|\mathbb{D}'|\leq(|\mathbb{D}|-1)\lceil (|\mathbb{D}|+1)/2 \rceil$,
 every line in  $AG(\mathcal{M}+R)$ contains at most $(|\mathbb{D}|-1)\lceil (|\mathbb{D}|+1)/2 \rceil$ points. Thus
 $AG(\mathcal{M}+R)$ contains at least three  noncollinear points.  It follows from (\ref{FS4TETB43dfg})  that
  $\varphi(A)$ are  adjacent with three  noncollinear points in $AG(\mathcal{M}+R)$. Applying Lemma \ref{gdkj7ssb33}, we obtain  $\varphi(A)\in\mathcal{M}+R$.

 Therefore, we always have $\varphi(A)\in (\mathcal{M}+R)\cup(\mathcal{N}+R)$ for any $\scriptsize A\in\left(
                \begin{array}{c}
                  \mathbb{D}^{k\times n} \\
                  0 \\
                \end{array}
              \right)$.
  Then
  $$\varphi\left(
                \begin{array}{c}
                  \mathbb{D}^{k\times n} \\
                  0 \\
                \end{array}
              \right)
\subseteq (\mathcal{M}+R)\cup(\mathcal{N}+R).$$
 Taking $k=m$, we get (\ref{tre543bcn877m}).
 $\qed$

\begin{rem}\label{rem-6trr3vbxej} \ Let  $n,m, p, q\geq 2$ be integers. By  \cite[Theorems 1.2 and 1.3]{Pazzis}, it is easy to see that
every  degenerate graph homomorphism $\varphi: \mathbb{F}_q^{m\times n}\rightarrow  \mathbb{F}_q^{p\times q}$ is a  (vertex) colouring.
When  $\mathbb{D}=\mathbb{D}'=\mathbb{F}_2$ or $\mathbb{F}_3$,
Both Theorem \ref{non-degenerate-b} and Corollary \ref{non-degenerate-ddD} still hold (cf. \cite{Pazzis, Huang-Li,Huang-Semrl-2014}).
However, when $q\leq 3$ with $q < q'$, it is an open problem to characterize the non-degenerate graph homomorphisms from $\mathbb{F}_q^{n\times n}$ to $\mathbb{F}_{q'}^{p\times q}$.
\end{rem}

\section*{Acknowledgments}

This work was supported by the National Natural Science Foundation of China (Project 11371072), and
the Scientific Research Fund of Hunan Provincial Education Department 16C0037.






\footnotesize




\end{document}